\documentclass[a4paper, twoside, 11pt]{article}

\usepackage[a4paper, margin=1in]{geometry}

\usepackage{amssymb, amsfonts, amsmath, amsthm}
\usepackage{hyperref} 
\usepackage{bm} 
\usepackage{graphicx} 
\usepackage[ruled, linesnumbered]{algorithm2e}
\usepackage{subcaption} 
\usepackage[title,titletoc]{appendix} 
\usepackage{multirow}
\usepackage{indentfirst} 
\usepackage{cleveref} 

\theoremstyle{plain}

\newtheorem{theorem}{Theorem}[section]

\theoremstyle{definition}
\newtheorem{definition}{Definition}[section]

\theoremstyle{remark}

\numberwithin{equation}{section}
\numberwithin{figure}{section}
\numberwithin{table}{section}

\def\R{\mathbb{R}}

\def\d{\mathrm{d}}

\def\diag{\mathrm{diag}}

\def\B#1{\left\{#1\right\}} 
\def\abs#1{\left\lvert#1\right\rvert}

\allowdisplaybreaks
\title{
	A positivity preserving and entropy stable nodal discontinuous Galerkin scheme for ideal MHD
}
\author{
	Yue~Wu\thanks{Division of Applied Mathematics, Brown University, Providence, RI 02912, USA (yue\_wu3@brown.edu).} \
	and Chi-Wang~Shu\thanks{Division of Applied Mathematics, Brown University, Providence, RI 02912, USA (chi-wang\_shu@brown.edu).}
}
\date{} 

\begin{document}

\maketitle


\textbf{Abstract:} Numerically solving magnetohydrodynamic (MHD) equations faces many challenges: avoiding divergence error, maintaining positivity, and satisfying entropy conditions. Among discontinuous Galerkin (DG) schemes, there has been a modal version that is locally divergence-free and positivity preserving and a nodal version that is semi-discretely entropy stable. In this work, we develop a DG scheme that combines the advantages of these two and solves all the three challenges. The key ingredients that bring these two schemes together are an HLL numerical flux with entropy stable signal speed estimates and a locally divergence-free projection. To handle problems with strong shocks, the essentially oscillation-free damping is applied. Various numerical experiments verify the accuracy and robustness of our method.



\section{Introduction}


The ideal magnetohydrodynamic (MHD) equations describe the motion of perfectly conducting quasi-neutral plasmas without viscosity or resistivity. It has wide application in scientific and engineering areas including geophysics, astrophysics, and fusion technology. Mathematically, the ideal MHD equations belong to the category of nonlinear hyperbolic conservation laws with the inherent magnetic divergence-free condition as an involution. For nonlinear hyperbolic conservation laws, it is well known that discontinuities can emerge even with smooth initial and boundary data, and complex structures including shocks, contacts, and rarefactions will occur thereby. As a result, weak solutions in the distributional sense have to be considered. However, they are generally not unique, and additional entropy conditions are required to find the physically relevant and unique solution thereby \cite{Hyperbolic_Dafermos_4th}, which states that the mathematical entropy remains constant in smooth regions and decreases at shocks. To ensure fidelity of the numerical solution, many schemes \cite{10.1017/S0962492902000156, 10.1137/110836961, 10.4208/csiam-am.2020-0003} are proposed with an emphasis on entropy stability.

One crucial aspect that distinguishes the MHD equations from other common hyperbolic systems is the divergence-free condition. Physically, it indicates that no magnetic monopoles ever exist, and mathematically, it will hold at future time if satisfied by the initial condition so as to be called an involution \cite{10.1007/BF00280911}. Such an involution is not merely a consequence, like in the Maxwell equations, but also a factor that influences the entropy stability of the solution as pointed out by \cite{10.1007/0-387-38034-5_4}. Existing studies have shown that violating this condition will deteriorate the numerical solution \cite{10.1016/0021-99918090079-0, 10.1006/jcph.2000.6519, 10.2140/camcos.2013.8.1}, and various numerical methods have been designed to tackle this problem, including the projection method \cite{10.1016/0021-99918090079-0}, the constrained transport method \cite{10.1086/166684, 10.1016/j.jcp.2004.11.016} and its unstaggered version \cite{10.1137/050627022, 10.1016/j.jcp.2011.02.009}, the generalized Lagrange multiplier method \cite{10.1006/jcph.2001.6961, 10.1016/j.jcp.2018.03.002}, the locally divergence-free (LDF) method \cite{10.1007/s10915-004-4146-4}, and the globally divergence-free method \cite{10.1007/s10915-018-0750-6, 10.1007/s10915-020-01289-8, 10.2140/camcos.2021.16.59}. Besides using the original conservative system, a non-conservative source term has been proposed independently by Godunov \cite{Godunov1972} to entropically symmetrize the system and by Powell \cite{Powell1994NASA, 10.1006/jcph.1999.6299} to remove the rank deficiency of Jacobian matrices and to satisfy Galilean invariance, which is later called the Godunov--Powell (GP) source term. It allows the second law of thermodynamics to be satisfied even when the divergence of magnetic field is nonzero, both analytically \cite{10.1007/978-3-642-58535-7_5} and in numerical schemes \cite{10.1016/j.jcp.2015.09.055, 10.1137/15M1013626, 10.1016/j.jcp.2017.10.043, 10.1016/j.jcp.2018.03.002}, and enjoy numerical stability even with non-zero divergence \cite{10.1006/jcph.1999.6299, 10.1137/18M1168042}.

High-order explicit shock capturing methods for hyperbolic conservation laws that are capable of achieving high-order accuracy in smooth regions and maintain stability near discontinuities have been widely studied and applied over the past decades, including the finite difference and finite volume weighted essentially non-oscillatory (WENO) \cite{10.1006/jcph.1996.0130, 10.1017/S0962492920000057} and the discontinuous Galerkin (DG) \cite{10.1006/jcph.1998.5892, 10.1007/978-3-642-59721-3} methods. The discontinuous Galerkin method is known for its arbitrary \(h/p\) adaptivity, geometric flexibility, extreme locality and parallel scalability. However, for nonlinear hyperbolic conservation laws such as the MHD equations, three major problems arise regarding robustness of the computation. The first one is the suppression of oscillation generated by the Gibbs phenomenon in the vicinity of discontinuities. A widely used cure is to apply slope limiters adopted from finite volume schemes, such as the total variation diminishing (TVD), the total variation bounded (TVB) \cite{0021999189901836}, and the WENO \cite{10.1137/S1064827503425298, 10.1016/j.jcp.2012.08.028} limiters, and another choice is some artificial damping that can be incorporated into the semi-discrete formulation like the essentially oscillation-free family \cite{10.1137/21M140835X, 10.1090/mcom/3998, 10.1016/j.jcp.2024.113498}. The second one is preserving admissibility of variables, which is also called the positivity preserving (PP) property in the context of MHD, i.e., preserving the positivity of density and pressure, so as to avoid simulation from breaking down. Rigorous analysis has been carried out and provably PP schemes have been designed in the series work of \cite{10.1137/18M1168017, 10.1137/18M1168042, 10.1007/s00211-019-01042-w}, where four main ingredients are necessary: locally divergence-free magnetic field, HLL numerical flux with suitably chosen signal speed estimates, compatible discretization of the GP source term on cell interfaces, and the Zhang--Shu limiter \cite{10.1016/j.jcp.2010.08.016}. The third one is about the entropy stability of the scheme. Semi-discretely entropy stable (ES) nodal DG schemes have been proposed by Chen and Shu (on Lobatto nodes) \cite{10.1016/j.jcp.2017.05.025}, Chan (the hybridized SBP operators) \cite{10.1016/j.jcp.2018.02.033, 10.1137/18M1209234}, Crean (the global SBP operators) \cite{10.1016/j.jcp.2017.12.015} based on the SBP-SAT \cite{10.1016/j.compfluid.2014.02.016} flux differencing framework, and Abgrall (artificial dissipation) \cite{10.1016/j.jcp.2018.06.031}. The readers are referred to \cite{10.4208/csiam-am.2020-0003} for a systematic review, and these methods are applied to the MHD equations in \cite{10.1016/j.jcp.2015.09.055, 10.1016/j.jcp.2017.10.043, 10.1016/j.jcp.2018.06.027, 10.1016/j.jcp.2022.111851}.

As pointed out by \cite{10.1016/j.jcp.2025.113911}, the PP property of the updated cell averages relies on the LDF property of the original magnetic field. Existing flux differencing type nodal semi-discretely ES schemes are not provably PP because of the divergence and the choice of interface numerical flux. Our work designs an algorithm that incorporates the LDF-PP framework into the semi-discretely ES DG scheme. We study the semi-discrete entropy production by the LF and HLL numerical fluxes, propose a bound for the signal speed estimates that makes the interface flux entropy stable, and build an essentially oscillation-free DG scheme that achieves both PP and ES.

The rest of the paper is organized as follows. In \Cref{se:pre}, we introduce the MHD equations with the GP source term and its entropy analysis. In \Cref{se:nm}, we construct a semi-discretely ES LF/HLL numerical flux, and describe the whole PP ES DG scheme. Numerical examples are provided in \Cref{se:nr}. Finally, we conclude in \Cref{se:con}.

\section{Preliminaries}\label{se:pre}

\subsection{The MHD equations with the Godunov--Powell source term}

The MHD equations with the GP source term in Einstein summation form \cite{jameson2006eigenvalues} are
\begin{equation}\label{eq: MHD equations (Einstein form)}
	\begin{cases}
		\frac{\partial \rho}{\partial t} + \frac{\partial}{\partial x_j} (\rho u_j) = 0, \\
		\frac{\partial (\rho u_i)}{\partial t} + \frac{\partial}{\partial x_j} (\rho u_i u_j - B_i B_j) + \frac{\partial}{\partial x_i} (p + \frac{\abs{\bm B}^2}{2}) + B_i \frac{\partial B_j}{\partial x_j} = 0, \\
		\frac{\partial B_i}{\partial t} + \frac{\partial}{\partial x_j} (B_i u_j - u_i B_j) + u_i \frac{\partial B_j}{\partial x_j} = 0, \\
		\frac{\partial E}{\partial t} + \frac{\partial}{\partial x_j} ( (E + p + \frac{\abs{B}^2}{2}) u_j - u_i B_i B_j ) + u_i B_i \frac{\partial B_j}{\partial x_j} = 0.
	\end{cases}
\end{equation}
They can also be rewritten into the following hyperbolic conservation law form augmented by a non-conservative source term as
\begin{equation}\label{eq: MHD equations (combined form)}
	\frac{\partial \bm w}{\partial t} + \nabla \cdot \bm f(\bm w) + \bm S(\bm w) (\nabla \cdot \bm B) = \bm 0,
\end{equation}
where the variables and flux functions are
\begin{equation*}
	\bm w = \begin{pmatrix}
		\rho \\
		\rho \bm u \\
		\bm B \\
		E
	\end{pmatrix} \in \R^{8 \times 1}, \quad
	\bm S(\bm w) = \begin{pmatrix}
		0 \\
		\bm B \\
		\bm u \\
		\bm u \cdot \bm B
	\end{pmatrix} \in \R^{8 \times 1},
\end{equation*}
\begin{equation*}
	\bm f = (\bm f_1, \bm f_2, \bm f_3)
	= \begin{pmatrix}
		\rho \bm u^T \\
		\rho \bm u \otimes \bm u - \bm B \otimes \bm B + (p + \frac{\abs{\bm B}^2}{2}) \bm I \\
		\bm B \otimes \bm u - \bm u \otimes \bm B \\
		(E + p + \frac{\abs{B}^2}{2}) \bm u^T - (\bm u \cdot \bm B) \bm B^T
	\end{pmatrix} \in \R^{8 \times 3},
\end{equation*}
\begin{equation*}
	\bm f_i = \begin{pmatrix}
		\rho u_i \\
		\rho u_i \bm u - B_i \bm B + (p + \frac{\abs{\bm B}^2}{2}) \bm e_i \\
		u_i \bm B - B_i \bm u \\
		(E + p + \frac{\abs{B}^2}{2}) u_i - (\bm u \cdot \bm B) B_i
	\end{pmatrix} \in \R^{8 \times 1}.
\end{equation*}
Here, \(\rho\) is the density, \(\bm u = (u_1, u_2, u_3)^T\) is the fluid velocity, \(\bm B = (B_1, B_2, B_3)^T\) is the magnetic field strength, \(E = \rho e + \frac{\rho}{2} \abs{\bm u}^2 + \frac{1}{2} \abs{\bm B}^2\) is the total energy, \(e\) is the specific internal energy, and \(p = p(\rho, e)\) is the pressure determined by the equation of state (EOS). In order to close the system, an equation of state has to be provided. We consider the \(\gamma\)-law EOS for ideal polytropic gas,
\begin{equation*}
	p = \rho e (\gamma - 1),
\end{equation*}
with \(\gamma = \frac{c_p}{c_v} \in (1, \frac{5}{3}]\) being the specific heat capacity ratio which is constant.

The third term on the left hand side of \eqref{eq: MHD equations (combined form)} is the GP source term. Taking divergence on the third equation of \eqref{eq: MHD equations (Einstein form)} gives the evolution equation of the magnetic divergence,
\begin{equation}
	\frac{\partial (\nabla \cdot \bm B)}{\partial t} + \nabla \cdot (\bm u (\nabla \cdot \bm B)) = 0,
\end{equation}
which states that the magnetic divergence is convected as a passive scalar \cite{Powell1994NASA} and coincides with the involution that \(\nabla \cdot \bm B = 0\) will hold for all \(t \geq 0\) if it holds at \(t = 0\).

\subsection{Entropy analysis}

Weak solutions to MHD equations are not necessarily unique. Additional partial differential inequalities, i.e., the entropy conditions, are needed to pick the physically relevant solution. For scalar equations, it is known that imposing entropy inequality picks the unique entropy solution among weak solutions \cite{Hyperbolic_Dafermos_4th}. For this specific system, we are concerned with the second law of thermodynamics. Let \(s(\bm w) = \ln(p \rho^{-\gamma})\) be the thermodynamic entropy. Then the mathematical entropy pair for the system is
\begin{equation}
	U(\bm w) = \frac{-\rho s}{\gamma - 1}, \quad
	\bm F(\bm w) = \frac{-\rho s \bm u^T}{\gamma - 1} \in \R^{1 \times 2},
\end{equation}
that satisfies
\begin{equation}
	\frac{\partial F_\alpha}{\partial \bm w} = \frac{\partial U}{\partial \bm w} \frac{\partial \bm f_\alpha}{\partial \bm w} + \phi(\bm v(\bm w)) \frac{\partial B_\alpha}{\partial \bm w}.
\end{equation}
Here, the entropy variable is
\begin{equation}
	\bm v(\bm w) = (\nabla_{\bm w} U)^T = (\frac{\gamma-s}{\gamma-1} - \beta \abs{\bm u}^2, 2 \beta \bm u^T, 2 \beta \bm B^T, -2\beta)^T \in \R^{8 \times 1}, \quad \beta = \frac{\rho}{2 p}.
\end{equation}
The coefficient variable in the Godunov--Powell source term is given by
\begin{equation}
	\phi(\bm v) = 2 \beta (\bm u \cdot \bm B) = -\frac{v_2 v_5 + v_3 v_6 + v_4 v_7}{v_8},
\end{equation}
such that
\begin{equation}
	\phi'(\bm v) \bm v = \phi(\bm v), \quad
	\phi'(\bm v) = S(\bm w(\bm v)) \in \R^{1 \times 8}.
\end{equation}
The corresponding potential function and potential flux are
\begin{equation}
	\varphi(\bm v) = \bm v^T \bm w - U = \rho + \beta \abs{\bm B}^2, \quad
	\bm\psi(\bm v) = \bm v^T \bm f + \phi(\bm v) \bm B^T - \bm F = \varphi \bm u^T \in \R^{1 \times 3},
\end{equation}
which satisfy
\begin{equation}
	\nabla_{\bm v}\varphi = \bm w^T, \quad
	\frac{\partial \bm\psi}{\partial \bm v} = \bm f^T + \bm B \otimes \phi'(\bm v).
\end{equation}
Then the system can be symmetrized into
\begin{equation}
	\varphi''(\bm v) \frac{\partial \bm v}{\partial t} + \sum_{m=1}^2 (\psi_m''(\bm v) - B_m \phi''(\bm v)) \frac{\partial \bm v}{\partial x_m} = 0,
\end{equation}
and the corresponding entropy inequality is
\begin{equation}
	\frac{\partial U(\bm w)}{\partial t} + \nabla \cdot \bm F(\bm w) \leq 0.
\end{equation}
The last inequality describes dissipation of mathematical entropy along with evolution, which is desirable also in numerical solutions.

\section{Numerical methods}\label{se:nm}

In this section, we present the positivity preserving and entropy stable nodal discontinuous Galerkin (PP-ES-DG) scheme for the ideal MHD equations with the Godunov--Powell source term. Following a brief review of the widely studied entropy stable nodal DG framework for this system \cite{10.1016/j.jcp.2017.10.043}, we present the key design of an entropy stable HLL numerical flux that incorporates entropy stability and positivity preservation. Then the whole algorithmic flow is described that integrates the essentially oscillation-free limiter and the PP limiter into the whole framework.

\subsection{The entropy stable nodal discontinuous Galerkin framework}

We introduce the key entropy related flux used in the derivation before presenting the numerical schemes.

\begin{definition}[entropy conservative flux]
	We call a symmetric and consistent flux pair \((\bm f_S, \bm B_S)\) \textit{entropy conservative} (EC) if
	\begin{equation*}
		\begin{aligned}
			\forall \bm w^-, \bm w^+, \quad
			& (\bm v(\bm w^+) - \bm v(\bm w^-))^T \bm f_S(\bm w^-, \bm w^+) + (\phi(\bm v(\bm w^+)) - \phi(\bm v(\bm w^-))) \bm B_S(\bm w^-, \bm w^+)^T \\
			& = \bm\psi(\bm v(\bm w^+))^T - \bm\psi(\bm v(\bm w^-))^T.
		\end{aligned}
	\end{equation*}
\end{definition}
Examples of the EC flux \(\bm f_S\) include the KEPEC flux \cite[Appendix B]{10.1016/j.jcp.2015.09.055} for the MHD system with the Janhunen source term \cite{10.1006/jcph.2000.6479}, the Chandrashekar flux \cite[Section 3.2]{10.1137/15M1013626} for the MHD system with the Godunov--Powell source term and the Hindenlang flux \cite{hindenlang2019new} (see \cite[Appendix A]{10.1016/j.jcp.2023.112607} for the formula) for the GLM-MHD system with the Godunov--Powell source term (that does not produce divergence error in the volume flux differencing). All these examples use the simplest arithmetic average magnetic field volume flux \(\bm B_S(\bm w^-, \bm w^+) = \frac{\bm B(\bm w^-) + \bm B(\bm w^+)}{2}\).

\begin{definition}[entropy stable flux]
	We call a consistent numerical flux pair \((\widehat{\bm f}, \widehat{\bm B})\) \textit{entropy stable} (ES) if
	\begin{equation*}
		\begin{aligned}
			\forall \bm w^-, \bm w^+, \quad
			& (\bm v(\bm w^+) - \bm v(\bm w^-))^T \widehat{\bm f}(\bm w^-, \bm w^+) + (\phi(\bm v(\bm w^+)) - \phi(\bm v(\bm w^-))) \widehat{\bm B}(\bm w^-, \bm w^+)^T \\
			& \leq \bm\psi(\bm v(\bm w^+))^T - \bm\psi(\bm v(\bm w^-))^T.
		\end{aligned}
	\end{equation*}
\end{definition}

For simplicity, we only present the 1D case for this subsection. Extension to 2D can be done by applying the 1D procedure dimension-by-dimension. The 1D system is
\begin{equation}\label{eq: MHD equations (1D)}
	\frac{\partial \bm w}{\partial t} + \frac{\partial \bm f_1(\bm w)}{\partial x} + \bm S(\bm w) \frac{\partial B_1}{\partial x} = \bm 0.
\end{equation}
We discretize the 1D domain \([a,b]\) into \(N\) uniform non-overlapping elements
\begin{equation*}
	I_i = [x_{i-\frac{1}{2}}, x_{i+\frac{1}{2}}], \quad
	1 \leq i \leq N,
\end{equation*}
where
\begin{equation*}
	x_{i+\frac{1}{2}} = a + i h, \quad h = \frac{b-a}{N_x}.
\end{equation*}
Let \(k \geq 1\) be the polynomial order. The DG space is
\begin{equation}
	\mathcal{V}_h^k = \B{\bm w_h: \bm w_h|_{I_{i}} \in \mathcal{P}^k(I_{i}), \forall I_{i}}.
\end{equation}
Here, \(\mathcal{P}^k(I_{i})\) denotes the space of polynomials of maximum order \(k\) on element \(I_{i}\). Let \(\B{x_{i}^{\mathrm{LGL}}}_{i=0}^{k}\) and \(\B{\omega_{i}}_{i=0}^{k}\) be the \((k+1)\)-point Legendre--Gauss--Lobatto (LGL) nodes and weights on the reference cell \([-1,1]\). Let \(\B{l_i(x)}_{i=0}^k\) be the LGL nodal basis where \(l_i(x_j^\mathrm{LGL}) = \delta_{ij}\). Define the differentiation matrix \(\bm D \in \R^{(k+1) \times (k+1)}\) (with index starting from zero) such that
\begin{equation*}
	l_j'(x) = \sum_{i=0}^k D_{ij} l_i(x).
\end{equation*}
By the algebraic degree of exactness of the LGL quadrature, it is straight forward to verify the summation-by-parts (SBP) property,
\begin{equation}\label{eq: SBP property}
	\bm W \bm D + \bm D^T \bm W = \bm B,
\end{equation}
where
\begin{equation*}
	\bm W = \diag\left(\omega_0^\mathrm{LGL}, \omega_1^\mathrm{LGL}, \cdots, \omega_k^\mathrm{LGL}\right),
\end{equation*}
\begin{equation*}
	\bm B = \diag\left(-1, 0, \cdots, 0, 1\right) = \diag\left(\tau_0, \cdots, \tau_k\right).
\end{equation*}
On each element \(I_i\), the numerical solution \(\bm w_h\) is represented by its nodal values on LGL nodes, \(\bm w_h(x_{i \alpha})\), where \(x_{i \alpha}\) is the affine-mapped LGL nodes on \(I_i\), i.e., \(x_{i \alpha} = x_{i} + \frac{h}{2} x_{\alpha}^\mathrm{LGL}\). The semi-discrete scheme of \eqref{eq: MHD equations (1D)} on a single element \(I_i\) for nodal value at \(x_{i \alpha}\) is
\begin{equation}\label{eq: 1D scheme}
	\frac{h}{2} \frac{\d}{\d t} \bm w_h(x_{i \alpha}) + 2 \sum_{\beta = 0}^{k} D_{\alpha \beta} \bm f_{1 S}(\bm w_h(x_{i \alpha}), \bm w_h(x_{i \beta})) + \sum_{\beta = 0}^{k} D_{\alpha \beta} \bm S(\bm w_h(x_{i \alpha})) B_1(x_{i \beta}) = \frac{\tau_\alpha}{\omega_\alpha} \bm F_{1 \ast}^{i \alpha}.
\end{equation}
Here, \(\bm f_{1 S}\) is the Chandrashekar EC flux \cite[Section 3.2]{10.1137/15M1013626}, and the star vector is
\begin{equation}
	\begin{aligned}
		& \bm F_{1 \ast}^{i 0} = \bm f_1(\bm w_h(x_{i-\frac{1}{2}}^{+})) - \widehat{\bm f_1}_{i-\frac{1}{2}} + \bm S(\bm w_h(x_{i-\frac{1}{2}}^{+})) (B_1(\bm w_h(x_{i-\frac{1}{2}}^{+})) - \widehat{B_1}_{i-\frac{1}{2}}), \\
		& \bm F_{1 \ast}^{i k} = \bm f_1(\bm w_h(x_{i+\frac{1}{2}}^{-})) - \widehat{\bm f_1}_{i+\frac{1}{2}} + \bm S(\bm w_h(x_{i+\frac{1}{2}}^{-})) (B_1(\bm w_h(x_{i+\frac{1}{2}}^{-})) - \widehat{B_1}_{i+\frac{1}{2}}), \\
		& \bm F_{1 \ast}^{i \alpha} = \bm 0, \quad 1 \leq \alpha \leq k-1.
	\end{aligned}
\end{equation}

\begin{theorem}
	The semi-discrete scheme is entropy conservative (or stable) if the interface numerical flux pair \((\widehat{\bm f}, \widehat{\bm B})\) is entropy conservative (or stable).
\end{theorem}

\subsection{Entropy stable HLL flux}

In \cite{10.1016/j.jcp.2017.05.025, 10.1016/j.jcp.2025.113911}, ES numerical flux is obtained by directly adding penalization of jump to the EC flux. However, this choice generally cannot guarantee the positivity of cell average after a single Euler forward step. As an alternative, we seek to directly use the HLL flux so that PP is guaranteed and ES can be achieved with our carefully derived signal speed estimates. Recall the PP numerical flux of Wu and Shu \cite{10.1007/s00211-019-01042-w} and its compatible normal magnetic flux. Unlike for the Euler equations where there exists a so-called two-rarefaction estimate \cite{10.1007/BF01414629} that provably bounds the physical signal speeds \cite{10.1016/j.jcp.2016.05.054}, there has not been any reliable one for the Godunov--Powell MHD system, so the entropy stability of the HLL flux is not studied yet. For the Janhunen MHD system, some results have been obtained for the Janhunen MHD system \cite{10.1007/s00211-007-0108-8, 10.1007/s00211-010-0289-4}, but the discretization of the source term is intrinsically different.

The semi-discrete entropy stability of the HLL surface discretization is equivalent to
\begin{equation}
	-S_L S_R \geq \frac{(S_R - S_L)(F_{1,R} - F_{1,L}) - (S_R \bm v_R - S_L \bm v_L)^T (\bm f_{1,R} - \bm f_{1,L}) - (S_R \phi_R - S_L \phi_L) (B_{1,R} - B_{1,L})}{(\bm v_R - \bm v_L)^T (\bm w_R - \bm w_L)}.
\end{equation}
Here, we always assume \(S_L \leq 0 \leq S_R\) and \(S_L < S_R\) for simplicity. The denominator to the right of the inequality is symmetric positive definite by strict convexity of the entropy. Let
\begin{subequations}\label{eq: AL and AR forms}
	\begin{equation}
		A_L := \frac{(F_{1,R} - F_{1,L}) - \bm v_L^T (\bm f_{1,R} - \bm f_{1,L}) - \phi_L (B_{1,R} - B_{1,L})}{(\bm v_R - \bm v_L)^T (\bm w_R - \bm w_L)}
		=: \frac{b_L(\bm w_L, \bm w_R)}{a(\bm w_L, \bm w_R)}
	\end{equation}
	and
	\begin{equation}
		A_R := \frac{(F_{1,R} - F_{1,L}) - \bm v_R^T (\bm f_{1,R} - \bm f_{1,L}) - \phi_R (B_{1,R} - B_{1,L})}{(\bm v_R - \bm v_L)^T (\bm w_R - \bm w_L)}
		=: \frac{b_R(\bm w_L, \bm w_R)}{a(\bm w_L, \bm w_R)}
	\end{equation}
\end{subequations}
be the linear coefficient of \(-S_L\) and \(S_R\) in to the right of the inequality, respectively, and let \(b(\bm w_L, \bm w_R) := b_L(\bm w_L, \bm w_R) + b_R(\bm w_L, \bm w_R)\). The condition can be rewritten as
\begin{equation}
	(-S_L - A_R) (S_R - A_L) \geq A_L A_R.
\end{equation}
We found the following ways to satisfy the conditions.
\begin{itemize}
	\item For the HLL case, require \(S_L \leq -A_R^+ - \sqrt{A_L^+ A_R^+}\) and \(S_R \geq A_L^+ + \sqrt{A_L^+ A_R^+}\).
	\item For the LF case, require \(-S_L = S_R \geq (A_L + A_R)^+\).
\end{itemize}
Since \(a(\cdot, \cdot)\) will be close to \(0^+\) as \(w_L\) and \(w_R\) are close to each other, we analyze the behavior near this only singular point. By chain rule,
\begin{equation*}
	\frac{\partial}{\partial \bm w_R} a(\bm w_L, \bm w_R)
	= (\bm v_R - \bm v_L)^T + (\bm w_R - \bm w_L)^T \frac{\partial \bm v_R}{\partial \bm w_R},
\end{equation*}
\begin{equation*}
	\frac{\partial^2}{\partial \bm w_R^2} a(\bm w_L, \bm w_R)|_{\bm w_R = \bm w_L = \bm w}
	= \frac{\partial \bm v}{\partial \bm w} + (\frac{\partial \bm v}{\partial \bm w})^T
	= 2 \frac{\partial \bm v}{\partial \bm w}
	= 2 \frac{\partial^2 U}{\partial \bm w^2}
	\succ 0,
\end{equation*}
\begin{equation*}
	\frac{\partial}{\partial \bm w_L} b_L(\bm w_L, \bm w_R)
	= -\left[(\bm f_{1,R} - \bm f_{1,L})^T + (B_{1,R} - B_{1,L}) \phi'(\bm v_L)\right] \frac{\partial \bm v_L}{\partial \bm w_L},
\end{equation*}
\begin{equation*}
	\frac{\partial}{\partial \bm w_R} b_L(\bm w_L, \bm w_R)
	= (\bm v_R - \bm v_L)^T \frac{\partial \bm f_{1,R}}{\partial \bm w_R} + (\phi_R - \phi_L) \frac{\partial B_{1,R}}{\partial \bm w_R},
\end{equation*}
\begin{equation*}
	\frac{\partial}{\partial \bm w_L} b_R(\bm w_L, \bm w_R)
	= (\bm v_R - \bm v_L)^T \frac{\partial \bm f_{1,L}}{\partial \bm w_L} + (\phi_R - \phi_L) \frac{\partial B_{1,L}}{\partial \bm w_L},
\end{equation*}
\begin{equation*}
	\frac{\partial}{\partial \bm w_R} b_R(\bm w_L, \bm w_R)
	= -\left[(\bm f_{1,R} - \bm f_{1,L})^T + (B_{1,R} - B_{1,L}) \phi'(\bm v_R)\right] \frac{\partial \bm v_R}{\partial \bm w_R},
\end{equation*}
\begin{equation*}
	\begin{aligned}
		& \frac{\partial^2}{\partial \bm w_L^2} b_L(\bm w_L, \bm w_R)|_{\bm w_L = \bm w_R = \bm w}
		= - \frac{\partial^2}{\partial \bm w_R^2} b_R(\bm w_L, \bm w_R)|_{\bm w_R = \bm w_L = \bm w}
		= (\frac{\partial \bm v}{\partial \bm w})^T (\frac{\partial \bm f_1}{\partial \bm w} + \phi'(\bm v)^T \frac{\partial B_x}{\partial \bm w}) \\
		= & \frac{\partial^2 U}{\partial \bm w^2} \bm J_1
		= \frac{\partial^2 F_1}{\partial \bm w^2} - \sum_{k=1}^{8} v^{(k)} \nabla_{\bm w}^2 f_{1,k},
	\end{aligned}
\end{equation*}
\begin{equation*}
	\begin{aligned}
		& \frac{\partial^2}{\partial \bm w_R^2} b_L(\bm w_L, \bm w_R)|_{\bm w_R = \bm w_L = \bm w}
		= -\frac{\partial^2}{\partial \bm w_L^2} b_R(\bm w_L, \bm w_R)|_{\bm w_L = \bm w_R = \bm w}
		= \left[(\frac{\partial \bm f_1}{\partial \bm w})^T + (\frac{\partial B_x}{\partial \bm w})^T \phi'(\bm v)\right] \frac{\partial \bm v}{\partial \bm w} \\
		= & \bm J_1^T \frac{\partial^2 U}{\partial \bm w^2}
		= \frac{\partial^2 F_1}{\partial \bm w^2} - \sum_{k=1}^{8} v^{(k)} \nabla_{\bm w}^2 f_{1,k}.
	\end{aligned}
\end{equation*}
Hence, both \(A_L\) and \(A_R\) are actually well-defined and smooth within the whole admissible set. Moreover, by noting that \(\frac{\partial^2 U}{\partial \bm w^2} \bm J_1\) is symmetric and so is \(\sqrt{\frac{\partial^2 U}{\partial \bm w^2}} \bm J_1 \sqrt{\frac{\partial^2 U}{\partial \bm w^2}}^{-1}\), we have the local bound
\begin{equation*}
	\limsup_{\bm w_L \to \bm w} A_L(\bm w_L, \bm w)
	= \limsup_{\bm w_R \to \bm w} A_L(\bm w, \bm w_R)
	= \max_{\bm h \neq \bm 0} \frac{\bm h^T \frac{\partial^2 U}{\partial \bm w^2} \bm J_1 \bm h}{2 \bm h^T \frac{\partial^2 U}{\partial \bm w^2} \bm h}
	= \max_{\bm h \neq \bm 0} \frac{\bm h^T \sqrt{\frac{\partial^2 U}{\partial \bm w^2}} \bm J_1 \sqrt{\frac{\partial^2 U}{\partial \bm w^2}}^{-1} \bm h}{2 \bm h^T \bm h}
	= \lambda_8(\bm w),
\end{equation*}
\begin{equation*}
	\liminf_{\bm w_L \to \bm w} A_L(\bm w_L, \bm w)
	= \liminf_{\bm w_R \to \bm w} A_L(\bm w, \bm w_R)
	= \lambda_1(\bm w),
\end{equation*}
\begin{equation*}
	\limsup_{\bm w_L \to \bm w} A_R(\bm w_L, \bm w)
	= \limsup_{\bm w_R \to \bm w} A_R(\bm w, \bm w_R)
	= -\lambda_1(\bm w),
\end{equation*}
\begin{equation*}
	\liminf_{\bm w_L \to \bm w} A_R(\bm w_L, \bm w)
	= \liminf_{\bm w_R \to \bm w} A_R(\bm w, \bm w_R)
	= -\lambda_8(\bm w),
\end{equation*}
\begin{equation*}
	\lim_{\bm w_L \to \bm w} (A_L + A_R)(\bm w_L, \bm w)
	= \lim_{\bm w_R \to \bm w} (A_L + A_R)(\bm w, \bm w_R)
	= 0.
\end{equation*}
Here, \(\lambda_1 \leq \lambda_2 \leq \cdots \leq \lambda_8\) are the eigenvalues of \(\bm J_1\). These estimates further show that when \(\bm w_L\) and \(\bm w_R\) are close to each other, i.e., in smooth regions, the HLL lower bound here will not much exceed the spectral radius and the LF lower bound here will even be nearly zero.

We remark that these bounds discussed above are only lower bounds for the sake of entropy stability. They should be used together with PP lower bounds \cite{10.1137/18M1168042, 10.1007/s00211-019-01042-w} to ensure stable computation.

\subsection{Damping and limiting of the numerical solution}\label{sse: damping and limiting}

Damping or limiters are needed to either suppress non-physical oscillation generated by strong shocks or to maintain admissibility of the numerical solution. A modified version of the OEDG damping \cite{10.1090/mcom/3998} is employed for suppressing oscillation, and the PP limiter of Wu and Shu \cite{10.1137/18M1168042} is used for maintaining positivity. One issue that the nodal DG scheme will not produce a numerical solution that is locally divergence free, which is a prerequisite for the PP limiter. The cure is to apply the LDF projection right after a new stage is computed. The resulting limiting procedure is in the following order: (1) apply the LDF projection, (2) apply the OEDG damping, (3) apply the PP limiter\footnote{The PP limiter should be conducted upon the LGL nodes rather than some optimal convex decomposition nodes, because positive nodal values are required in computing the volume EC flux.} if needed. The OEDG damping is a convex limiter so that it will not deteriorate the LDF property, and both the LDF projection and the OEDG damping preserves the cell average, hence the PP limiter is guaranteed to work when the vanilla updated cell average is positive. Moreover, the PP limiter is also a convex limiter and hence the final numerical solution is also locally divergence-free.

We consider the 2D scheme with uniform grid discretization of the rectangular domain \(\Omega = [a, b] \times [c, d]\),
\begin{equation*}
	I_{ij} = [x_{i-\frac{1}{2}}, x_{i+\frac{1}{2}}] \times [y_{j-\frac{1}{2}}, y_{j+\frac{1}{2}}], \quad
	1 \leq i \leq N_x, 1 \leq j \leq N_y,
\end{equation*}
where
\begin{equation*}
	\begin{aligned}
		& x_{i+\frac{1}{2}} = a + i h_x, \quad h_x = \frac{b-a}{N_x}, \\
		& y_{j+\frac{1}{2}} = c + j h_y, \quad h_y = \frac{d-c}{N_y}.
	\end{aligned}
\end{equation*}

\subsubsection{The LDF projection}

The LDF projection is an element-wise post-processing introduced in \cite{10.4208/cicp.180515.230616a}. From now on, we only discuss the operations on a single physical element \(I_{ij}\).

Let \(K = [-1,1]^2\) be the reference element with coordinate \(\widehat{\bm x}\) and let \(\bm x\) be the coordinate of \(I_{ij}\). Define the push-forward mapping \(\bm T: K \mapsto I_{ij}\) s.t. \(\bm x = \bm T(\widehat{\bm x})\). In our setting, its Jacobian is a constant matrix
\begin{equation}
	\bm J = \frac{\partial \bm x}{\partial \widehat{\bm x}}
	= \diag\left(\frac{h_x}{2}, \frac{h_y}{2}\right).
\end{equation}
On the reference element \(K\), define the polynomial spaces
\begin{equation}
	\begin{aligned}
		& \widehat{\mathcal{V}} := \mathcal{Q}^{k}(K), \\
		& \widehat{\mathcal{W}} := \widehat{\mathcal{V}}^2, \\
		& \widehat{\mathcal{M}} := \mathcal{Q}^{k-1,k}(K) + \mathcal{Q}^{k,k-1}(K).
	\end{aligned}
\end{equation}
Here, \(\widehat{\mathcal{V}}\) and \(\widehat{\mathcal{W}}\) are the approximation spaces for scalar and 2D-vector variables, and \(\widehat{\mathcal{M}}\) is the divergence space of \(\widehat{\mathcal{W}}\). Following \cite[Definition 9.8]{10.1007/978-3-030-56341-7}, we use the pull-back transformation \(\psi^g(v) := v \circ \bm T\) to map scalar functions on \(I_{ij}\) to those on \(K\) and the contravariant Piola transformation \(\psi^d(\bm v) := \det(\bm J) \bm J^{-1} (\bm v \circ \bm T)\) to map 2D-vector functions on \(I_{ij}\) to those on \(K\). Since \(\bm J\) is diagonal and constant, the physical approximating spaces can be rewritten as
\begin{equation}
	\begin{aligned}
		& \mathcal{V} = \mathcal{Q}^{k}(I_{ij}) = (\psi^g)^{-1}(\widehat{\mathcal{V}}), \\
		& \mathcal{W} = [\mathcal{V}]^2 = (\psi^d)^{-1}(\widehat{\mathcal{W}}), \\
		& \mathcal{M} = (\psi^g)^{-1}(\widehat{\mathcal{M}}).
	\end{aligned}
\end{equation}
We use the variational formulation that seeks the LDF projected solution \(\bm B_h \in \mathcal{W}\) from an input \(\bm B_{h, 0} \in \mathcal{W}\),
\begin{equation}\label{eq: physical variation form}
	\begin{aligned}
		& \text{find \((\bm B_h, p_h) \in \mathcal{W} \times \mathcal{M}\) s.t. \(\forall (\bm v_h, q_h) \in \mathcal{W} \times \mathcal{M}\),} \\
		& \begin{cases}
			\int_{I_{ij}} (\bm J^{-1} \bm B_h)^T (\bm J^{-1} \bm v_h) \det(\bm J) \d \bm x + \int_{I_{ij}} p_h \nabla \cdot \bm v_h \d \bm x = \int_{I_{ij}} (\bm J^{-1} \bm B_{h, 0})^T (\bm J^{-1} \bm v_h) \det(\bm J) \d \bm x, \\
			\int_{I_{ij}} \nabla \cdot \bm B_h q_h \d \bm x = 0.
		\end{cases}
	\end{aligned}
\end{equation}
Denote the Piola mapped functions by
\begin{equation*}
	\widehat{\bm B_h} = \psi^d(\bm B_h), \quad
	\widehat{\bm B_{h, 0}} = \psi^d(\bm B_{h, 0}), \quad
	\widehat{p_h} = \psi^g(p_h).
\end{equation*}
Then \eqref{eq: physical variation form} is \textit{equivalent} to the following one on the reference cell \(K\):
\begin{equation}\label{eq: reference variation form}
	\begin{aligned}
		& \text{find \((\widehat{\bm B_h}, \widehat{p_h}) \in \widehat{\mathcal{W}} \times \widehat{\mathcal{M}}\) s.t. \(\forall (\widehat{\bm v_h}, \widehat{q_h}) \in \widehat{\mathcal{W}} \times \widehat{\mathcal{M}}\),} \\
		& \begin{cases}
			\int_{K} \widehat{\bm B_h}^T \widehat{\bm v_h} \d \widehat{\bm x} + \int_{K} \widehat{p_h} \nabla \cdot \widehat{\bm v_h} \d \widehat{\bm x} = \int_{K} \widehat{\bm B_{h, 0}}^T \widehat{\bm v_h} \d \widehat{\bm x}, \\
			\int_{K} \nabla \cdot \widehat{\bm B_h} \widehat{q_h} \d \widehat{\bm x} = 0.
		\end{cases}
	\end{aligned}
\end{equation}
The well-posedness of this formulation \eqref{eq: reference variation form} is easy to check by the Lagrange multiplier theory. Owing to the commuting property \cite[Corollary 9.9]{10.1007/978-3-030-56341-7}, \(\widehat{\bm B_h}\) is divergence-free implies that \(\bm B_h\) is also divergence-free. The reference-cell formulation \eqref{eq: reference variation form} does not contain the Jacobian \(\bm J\), so its matrices can be pre-computed and stored. We use Legendre modal coefficient representations for \(\widehat{\bm B_{h, 0}}\) and \(\widehat{\bm B_h}\), which is easy to transfer from and to nodal values at quadrature nodes. The LDF projection algorithm flow is described in \Cref{algo: LDF projection}. Here, \(\bm V\) is the Vandermonde matrix of 1D Legendre basis at LGL nodes, i.e., \(V_{\alpha \beta} = p_\beta(x_\alpha^\mathrm{LGL})\), where \(p_\beta\) is the \(\beta\)-th Legendre polynomial. The Legendre modal coefficient-based solver modal\_LDF\_ref for \eqref{eq: reference variation form} is simply constructed by forming the linear saddle-point system and solving with the method of Schur complement. Intermediate matrices can be saved before starting the time stepping.

\begin{algorithm}[htbp]
	\caption{LDF projection for cell \(I_{ij}\)}
	\label{algo: LDF projection}

	\KwIn{values of \(\bm B_{h, 0}\) at \((k+1)^2\) LGL nodes of \(I_{ij}\) in matrix form: \([B_{1, h, 0}]_{\alpha \beta} = B_{1,h,0}(x_i + \frac{h_x}{2}x_\alpha^\mathrm{LGL}, y_j + \frac{h_y}{2}x_\beta^\mathrm{LGL})\) and \([B_{2, h, 0}]_{\alpha \beta} = B_{2,h,0}(x_i + \frac{h_x}{2}x_\alpha^\mathrm{LGL}, y_j + \frac{h_y}{2}x_\beta^\mathrm{LGL})\)}

	\KwOut{values of divergence-free \(\bm B_{h}\) at \((k+1)^2\) LGL nodes of \(I_{ij}\) in matrix form: \([B_{1, h}]_{\alpha \beta} = B_{1,h}(x_i + \frac{h_x}{2}x_\alpha^\mathrm{LGL}, y_j + \frac{h_y}{2}x_\beta^\mathrm{LGL})\) and \([B_{2, h}]_{\alpha \beta} = B_{2,h}(x_i + \frac{h_x}{2}x_\alpha^\mathrm{LGL}, y_j + \frac{h_y}{2}x_\beta^\mathrm{LGL})\)}

	\tcp{Compute the Legendre coefficient matrices of \(\widehat{\bm B_{h,0}}\).}
	\([\widehat{B_{1,h,0}}] \leftarrow \frac{h_y}{2} \bm V^{-1} [B_{1, h, 0}] \bm V^{-T}\), \([\widehat{B_{2,h,0}}] \leftarrow \frac{h_x}{2} \bm V^{-1} [B_{2, h, 0}] \bm V^{-T}\) \;

	\tcp{Conduct the LDF projection on the reference element \(K\) using \eqref{eq: reference variation form}.}
	\(([\widehat{B_{1,h}}], [\widehat{B_{2,h}}]) \leftarrow \text{modal\_LDF\_ref}([\widehat{B_{1,h,0}}], [\widehat{B_{2,h,0}}])\) \;

	\tcp{Transform to nodal values of \(\bm B_{h}\) in matrix form.}
	\([B_{1,h}] \leftarrow \frac{2}{h_y} \bm V [\widehat{B_{1,h}}] \bm V^T\), \([B_{2,h}] \leftarrow \frac{2}{h_x} \bm V [\widehat{B_{2,h}}] \bm V^T\) \;
\end{algorithm}

\subsubsection{The OEDG damping}

Given a time step (or a Runge--Kutta stage) \(\bm w_h^n\), the modified OEDG damping procedure for cell \(I_{ij}\) is
\begin{equation}\label{eq: OEDG damping}
	\begin{cases}
		\frac{\partial \bm w_h}{\partial t} + \sum_{m=0}^k \delta_{I_{ij}}^m(\bm w_h^n) (\bm w_h - \overline{\bm w_h^n}|_{I_{ij}}) = \bm 0, & t \in [t^n, t^n + \tau], \\
		\bm w_h|_{t=t^n} = \bm w_h^n.
	\end{cases}
\end{equation}
where \(\overline{\bm w_h}|_{I_{ij}}\) is the cell average of \(\bm w_h\) on \(I_{ij}\). This pseudo evolution with parameter \(\tau\) is added after forming any Runge--Kutta stage, and it can be solved exactly by
\begin{equation}
	\widetilde{\bm w_h^n} = \overline{\bm w_h^n}|_{I_{ij}} + e^{-\tau \sum_{m=0}^k \delta_{I_{ij}}^m(\bm w_h^n)} (\bm w_h - \overline{\bm w_h^n}|_{I_{ij}}).
\end{equation}
The coefficients \(\delta_{I_{ij}}^m\) that measure the jump of all orders of derivatives are the same as \cite[Eqn. (3.9)]{10.1090/mcom/3998}.

\subsubsection{The PP limiter}

After the LDF projection and OEDG damping, we apply the Zhang--Shu PP limiter \cite{10.1016/j.jcp.2010.08.016} to the nodal solution in each cell. The detailed steps are just the same as \cite[Subsection 3.2.2]{10.1137/18M1168042}. Since our numerical solution is represented by nodal values at LGL nodes, the limiter can directly work on these values rather than modal coefficients. The remaining question is the validity: Is the cell average after a Euler forward step still admissible?

For the 1D scheme \eqref{eq: 1D scheme} for cell \(I_{i}\), multiply with \(\omega_\alpha\) and sum over \(\alpha\) to get the evolution of cell average,
\begin{equation}\label{eq: 1D evolution of cell average}
	\begin{aligned}
		& \frac{h}{2}\frac{\d}{\d t}\overline{\bm w_h}_{I_i} + \sum_{\alpha,\beta = 0}^k \omega_\alpha \bm S(\bm w_h(x_{i \alpha})) D_{\alpha\beta} B_1(x_{i \beta}) \\
		& = \widehat{\bm f_1}_{j-\frac{1}{2}} - \widehat{\bm f_1}_{j+\frac{1}{2}} + \bm S(\bm w_h(x_{i+\frac{1}{2}}^{-})) (B_1(\bm w_h(x_{i+\frac{1}{2}}^{-})) - \widehat{B_1}_{j+\frac{1}{2}}) - \bm S(\bm w_h(x_{i-\frac{1}{2}}^{+})) (B_1(\bm w_h(x_{i-\frac{1}{2}}^{+})) - \widehat{B_1}_{j-\frac{1}{2}}).
	\end{aligned}
\end{equation}
Here, we have used the SBP property \eqref{eq: SBP property}. Since the 2D scheme is a dimension-by-dimension application of the 1D version, its evolution of cell average can be easily derived from \eqref{eq: 1D evolution of cell average},
\begin{equation}\label{eq: 2D evolution of cell average}
	\frac{h_x h_y}{4} \frac{\d}{\d t} \overline{\bm w_h}|_{I_{ij}} + \sum_{\alpha,\beta,p,q=0}^k \omega_\alpha \omega_\beta \bm S(\bm w_h(\bm x_{ij \alpha\beta})) (D_{\alpha p} B_1(\bm x_{ij p \beta}) + D_{\beta q} B_2(\bm x_{ij \alpha q})) = \mathrm{RHS}.
\end{equation}
The RHS term here contains only boundary terms and is in the same form of the derivative term in \cite[Eqn. (63)]{10.1007/s00211-019-01042-w}. Therefore, the Euler forward step will produce a positive updated cell average, according to the rigorous analysis by Wu and Shu \cite{10.1007/s00211-019-01042-w}, if the second term to the left of the equality vanish. In fact,
\begin{equation*}
	\sum_{p,q=0}^k (D_{\alpha p} B_1(\bm x_{ij p \beta}) + D_{\beta q} B_2(\bm x_{ij \alpha q}))
\end{equation*}
is exactly the divergence value at the \((\alpha,\beta)\)-node of cell \(I_{ij}\), hence this term indeed vanishes if the last step has an LDF magnetic field. When time discretization schemes that can be rewritten as convex combination of Euler forward steps are used, and the LDF projection is applied both to the initial condition and immediately after each RK stage is formed, such condition holds so that the updated cell averages are admissible, because both the OEDG damping and PP limiter are convex limiters. 

We remark that although the original Wu--Shu framework of positivity-preserving MHD schemes are based on the modal DG or finite volume framework, its core is about proving the updated cell average is still admissible, which only relies on the cell average evolution equation. In our nodal DG setting, the cell average evolution equation \eqref{eq: 1D evolution of cell average} or \eqref{eq: 2D evolution of cell average} has the same form as that in \cite{10.1007/s00211-019-01042-w} which proves to provide a positive updated cell average. The proof of positivity for updated cell averages is based on the geometric quasilinearization framework of Wu and Shu \cite{10.1137/21M1458247}.

\subsection{Time discretization}

We use the third order strong stability preserving (SSP) Runge--Kutta (RK) method \cite{10.1090/S0025-5718-98-00913-2, 10.1137/S003614450036757X} to integrate the system in time. The limiting procedure in the previous subsection is applied after each new stage is computed. We use \(L_h\) to denote the semi-discrete operator of \eqref{eq: MHD equations (combined form)} and \(P_\tau\) to denote the limiting procedure described in the previous subsection where the OEDG part is called with time step \(\tau\). The fully discrete scheme is:
\begin{equation}
	\begin{aligned}
		& \bm w_h^{n,1} = P_{\Delta t}(\bm w_h^n + \Delta t L_h(\bm w_h^n)), \\
		& \bm w_h^{n,2} = P_{\Delta t}(\frac{3}{4}\bm w_h^n + \frac{1}{4}(\bm w_h^{n,1} + \Delta t L_h(\bm w_h^{n,1}))), \\
		& \bm w_h^{n+1} = P_{\Delta t}(\frac{1}{3}\bm w_h^{n} + \frac{2}{3}(\bm w_h^{n,2} + \Delta t L_h(\bm w_h^{n,2}))).
	\end{aligned}
\end{equation}
Since the \(P_{\Delta t}\) damping and limiting operator provides a locally divergence-free and positive (on LGL nodes) numerical solution according to our analysis in \Cref{sse: damping and limiting}, the final SSP-RK result \(\bm w_h^{n+1}\) will also satisfy these constraints.

For the 2D computation, the time step size \(\Delta t\) is determined by the following CFL condition:
\begin{equation}
	\Delta t = \min_{i j \alpha \beta} \frac{\mathrm{CFL}}{\frac{\sigma_{x, ij\alpha\beta}}{h_x} + \frac{\sigma_{y, ij\alpha\beta}}{h_y}}.
\end{equation}
Here, \(\sigma_{x, ij\alpha\beta}\) and \(\sigma_{y, ij\alpha\beta}\) are the maximum wave speed estimates in \(x\)- and \(y\)-directions on the \((\alpha, \beta)\)-node of cell \(I_{ij}\). For provable preservation of cell average positivity after a single forward Euler step from a positive state, strict CFL conditions as in \cite[Eqn. (20)]{10.1137/18M1168042} (for the LF case) or \cite[Eqn. (55)]{10.1007/s00211-019-01042-w} (for the HLL case) should be used, 
which is rather restrictive.  In our implementation, a simple try-and-correct procedure, which was used in previous literature as well, is used instead: we use the standard CFL condition for linear stability, and if there are negative cell averages in an intermediate RK stage when stepping from \(t^n\) to \(t^{n+1} = t^n + \Delta t\), just half the time step \(\Delta t\) and restart the whole RK procedure from \(t^n\).  This is repeated if necessary.  The value of the theoretical proof of a CFL condition to guarantee positivity is that this procedure of halving the time step will need to be carried for at most a fixed number of times.

\section{Numerical results}\label{se:nr}


In this section, we present numerical results to validate our proposed scheme. In all the tests, we choose \(k = 2\) and \(\mathrm{CFL} = \frac{0.6}{2k+1}\) and use \(\gamma = \frac{5}{3}\) if not otherwise specified. When computing \(A_L\) and \(A_R\) in \eqref{eq: AL and AR forms}, a small positive number \(\epsilon = 10^{-8}\) is added to the denominator to avoid division by zero.

\subsection{Smooth Alfv\'en wave}

This examples describes the circularly polarized Alfv\'en wave propagating in the periodic domain \([0, \frac{1}{\cos\alpha}] \times [0, \frac{1}{\sin\alpha}]\) \cite{10.1006/jcph.2000.6519}. The initial condition is given by the primitive variables
\begin{align*}
	& \rho = 1, \quad p = 0.1, \\
	& (u_1, u_2, u_3) = (-v_\perp \sin\alpha, v_\perp \cos\alpha, 0.1 \cos(2\pi x_\parallel)), \\
	& (B_1, B_2, B_3) = (B_\parallel \cos\alpha - B_\perp \sin\alpha, B_\parallel\sin\alpha + B_\perp \cos\alpha, u_3).
\end{align*}
The auxiliary variables are \(B_\parallel = 1\), \(B_\perp = v_\perp = 0.1 \sin(2\pi x_\parallel)\), \(x_\parallel = x \cos \alpha + y \sin \alpha\). The angle is chosen as \(\alpha = \frac{\pi}{6}\). The final time is taken to be \(T = 5\). \(N \times N\) uniform cells are used for spatial discretization. The result is shown in \Cref{tb: smooth alfven}. Desired \((k+1)\)-th order convergence rate is observed.

\begin{table}[htb]
	\centering
	\caption{Convergence history for \(B_\perp\) of the smooth Alfv\'en wave example}
	\label{tb: smooth alfven}
	\begin{tabular}{llllll}
		\hline
		\(k\) & \(N\) & \(L^2\) error & order & \(L^\infty\) error & order \\
		\hline
		1 & 32  & 2.6765E-02 & - & 2.3963E-02 & - \\
		~ & 64  & 6.1594E-03 & 2.119  & 6.4522E-03 & 1.893  \\
		~ & 128  & 1.4426E-03 & 2.094  & 1.6130E-03 & 2.000  \\
		~ & 256  & 3.4588E-04 & 2.060  & 3.7271E-04 & 2.114  \\
		~ & 512  & 8.5301E-05 & 2.020  & 8.7290E-05 & 2.094  \\
		\hline
		2 & 32  & 3.2991E-04 & - & 3.5541E-04 & - \\
		~ & 64  & 2.2516E-05 & 3.873  & 2.5284E-05 & 3.813  \\
		~ & 128  & 1.4739E-06 & 3.933  & 1.8771E-06 & 3.752  \\
		~ & 256  & 9.8033E-08 & 3.910  & 1.7033E-07 & 3.462  \\
		~ & 512  & 7.2171E-09 & 3.764  & 1.7112E-08 & 3.315  \\
		\hline
		3 & 32  & 1.6173E-05 & - & 1.7444E-05 & - \\
		~ & 64  & 3.0694E-07 & 5.719  & 5.0685E-07 & 5.105  \\
		~ & 128  & 1.3146E-08 & 4.545  & 1.7136E-08 & 4.886  \\
		~ & 256  & 5.1646E-10 & 4.670  & 9.0525E-10 & 4.243 \\
		\hline
	\end{tabular}
\end{table}

\subsection{Yee--Sj\"ogreen 2D Riemann problem}

\begin{figure}[htbp]
	\centering

	\begin{subfigure}{0.3\linewidth}
		\centering
		\includegraphics[width=\linewidth]{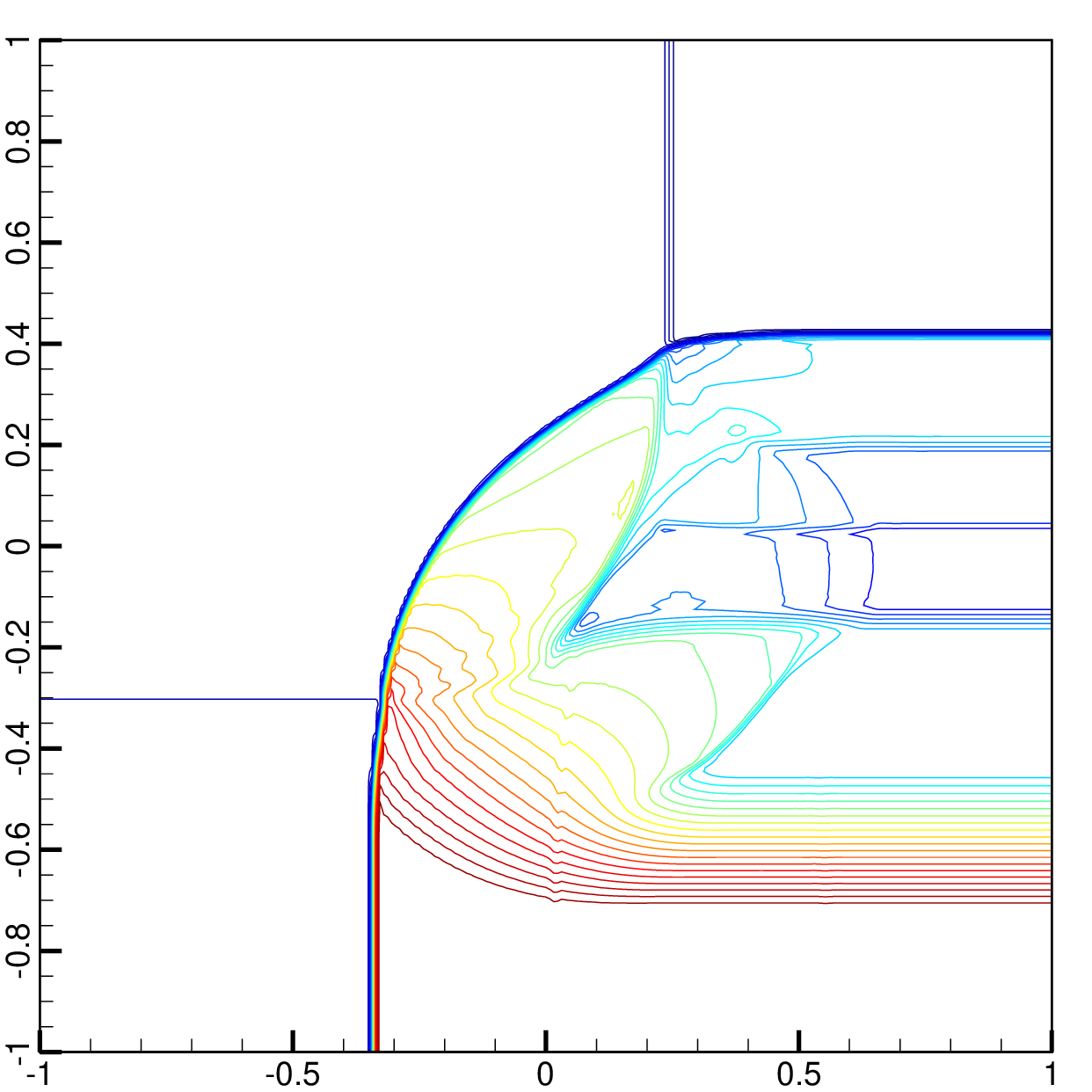}
		\caption{\(\rho\), \(256 \times 256\) grid}
	\end{subfigure}
	\begin{subfigure}{0.3\linewidth}
		\centering
		\includegraphics[width=\linewidth]{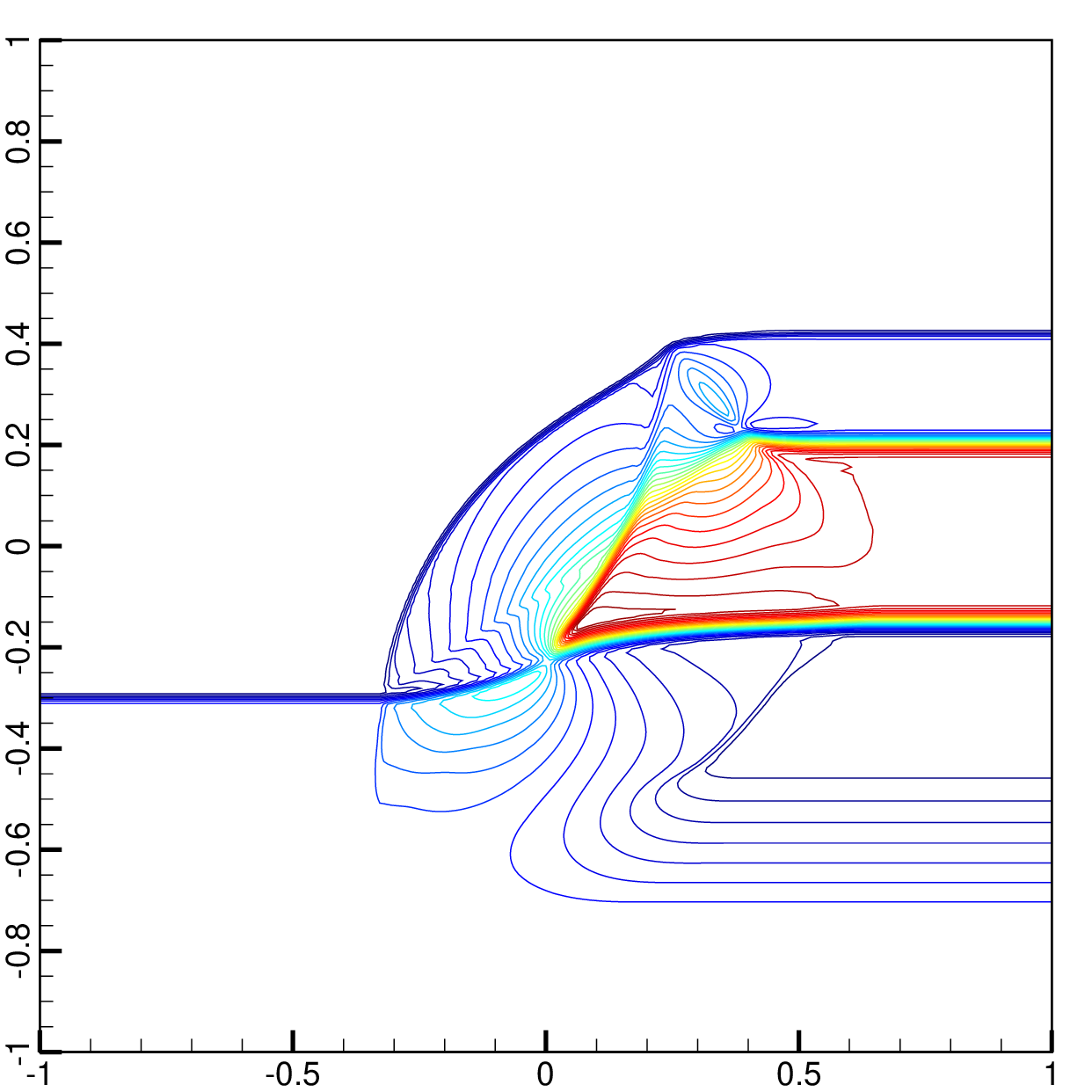}
		\caption{\(B_1\), \(256 \times 256\) grid}
	\end{subfigure}
	\begin{subfigure}{0.3\linewidth}
		\centering
		\includegraphics[width=\linewidth]{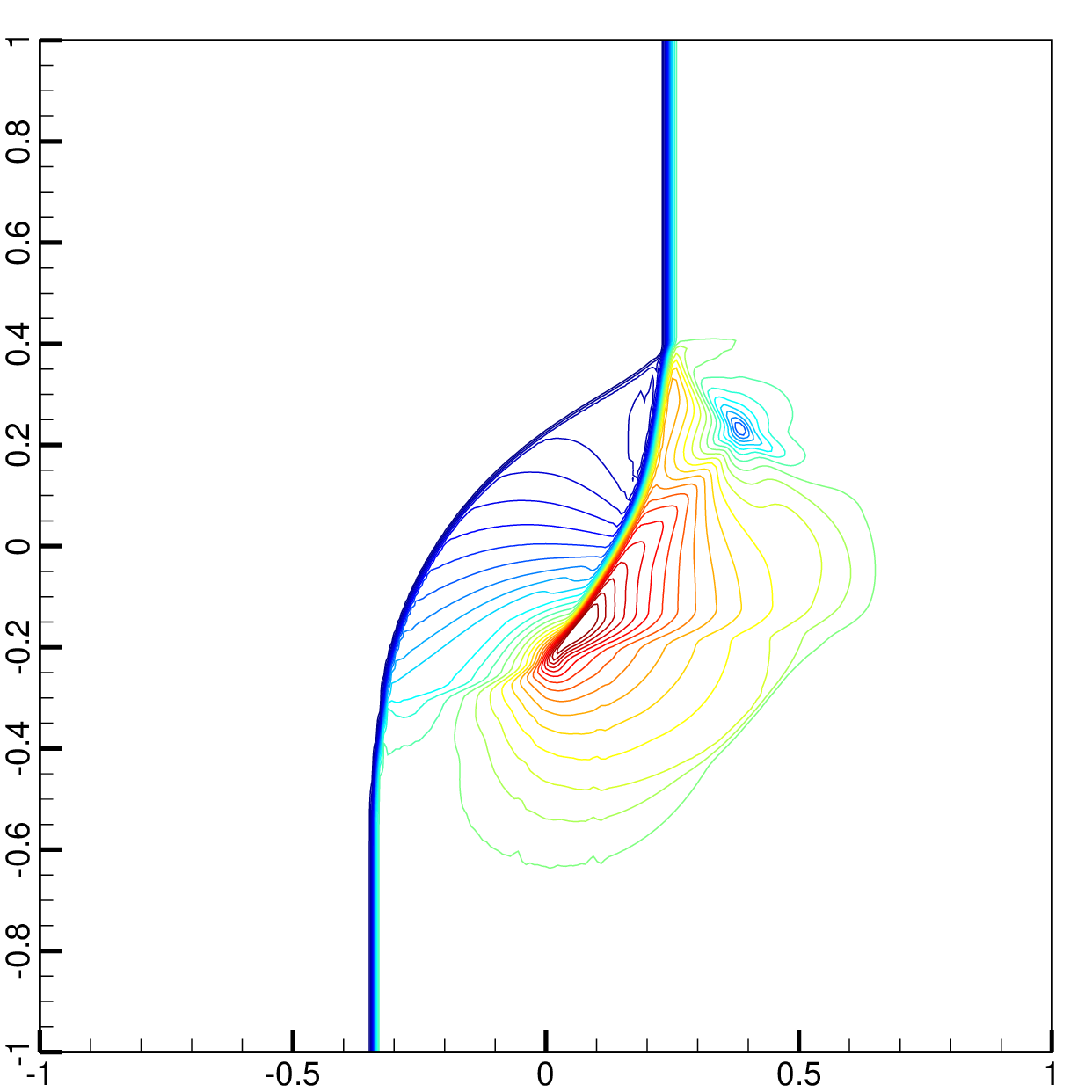}
		\caption{\(B_2\), \(256 \times 256\) grid}
	\end{subfigure}

	\begin{subfigure}{0.3\linewidth}
		\centering
		\includegraphics[width=\linewidth]{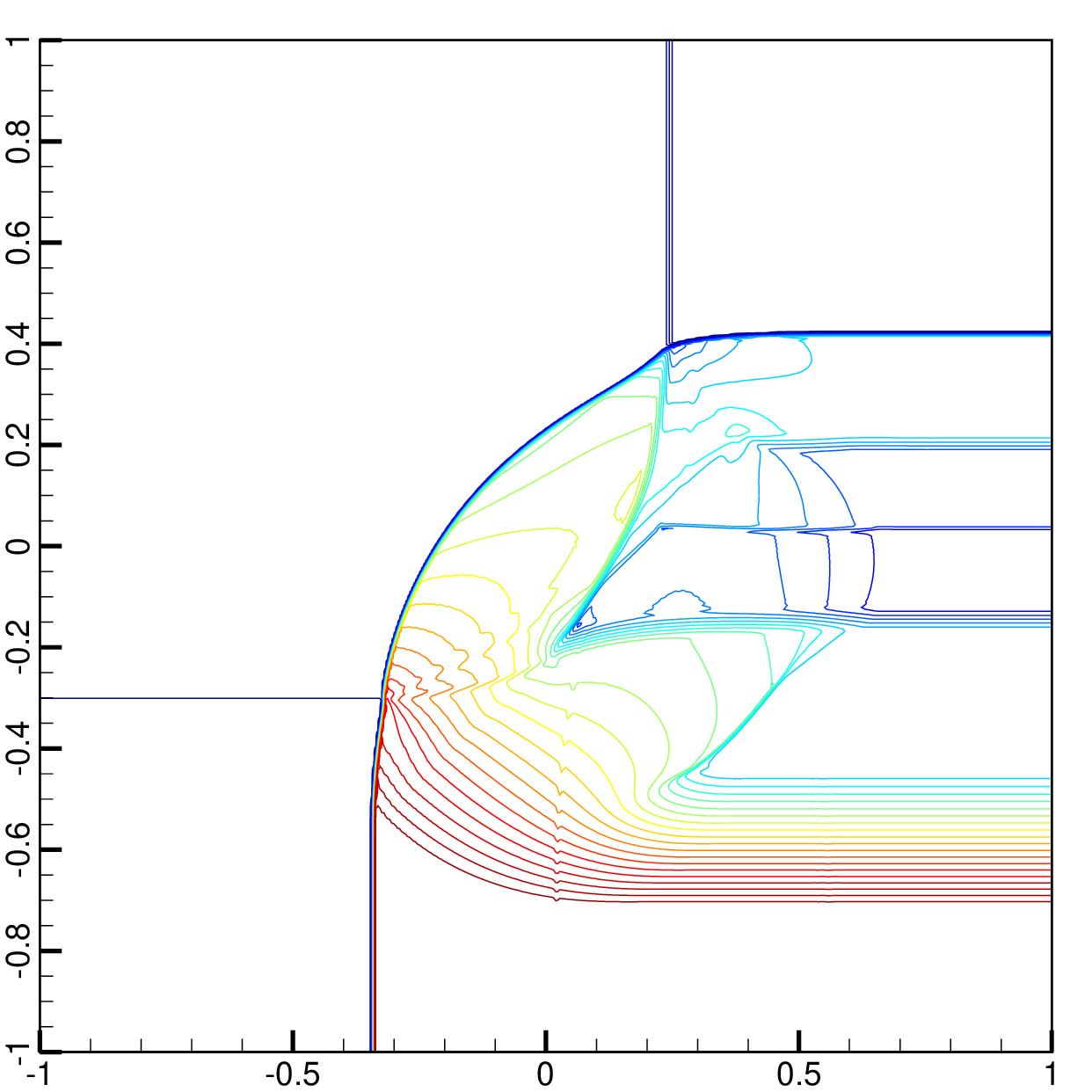}
		\caption{\(\rho\), \(512 \times 512\) grid}
	\end{subfigure}
	\begin{subfigure}{0.3\linewidth}
		\centering
		\includegraphics[width=\linewidth]{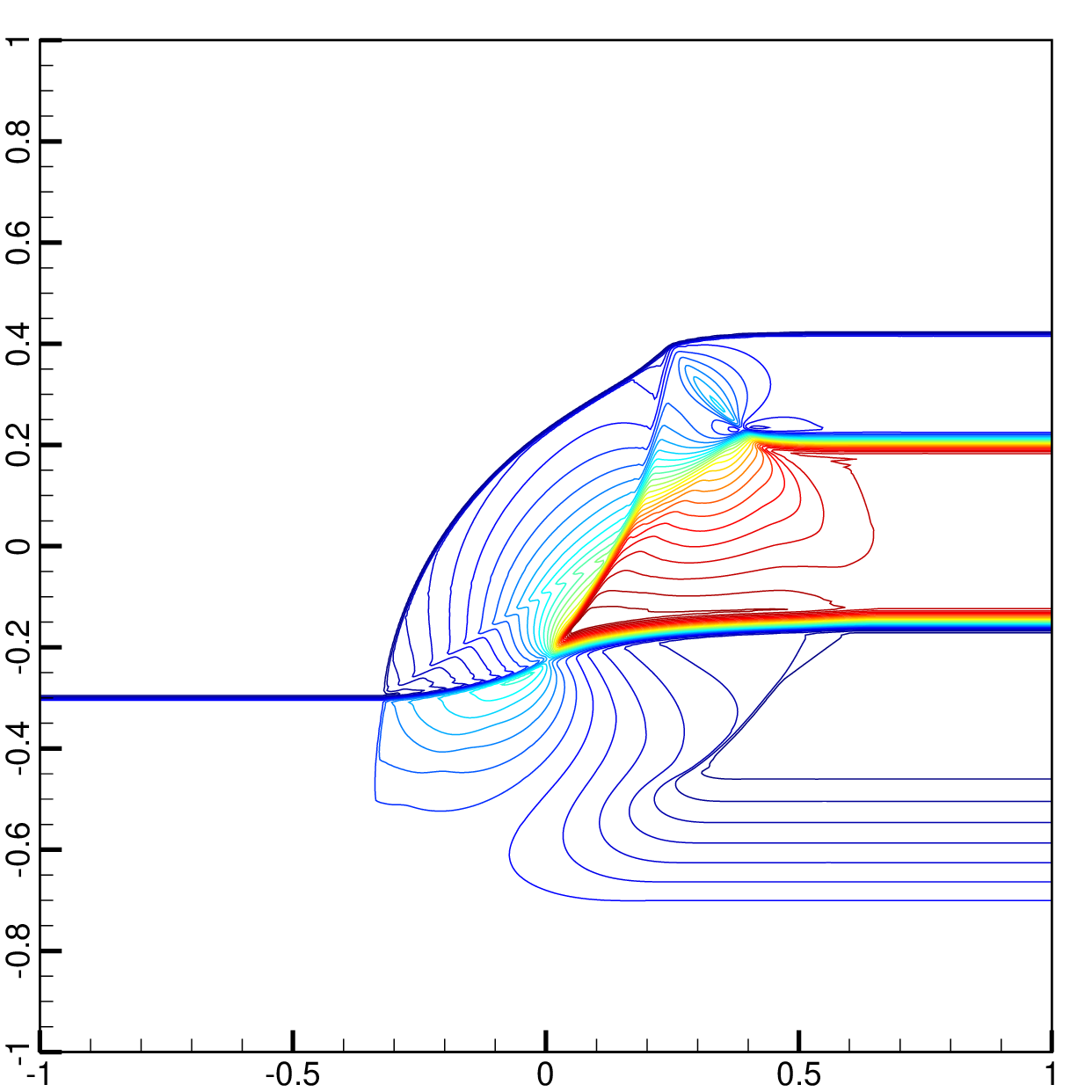}
		\caption{\(B_1\), \(512 \times 512\) grid}
	\end{subfigure}
	\begin{subfigure}{0.3\linewidth}
		\centering
		\includegraphics[width=\linewidth]{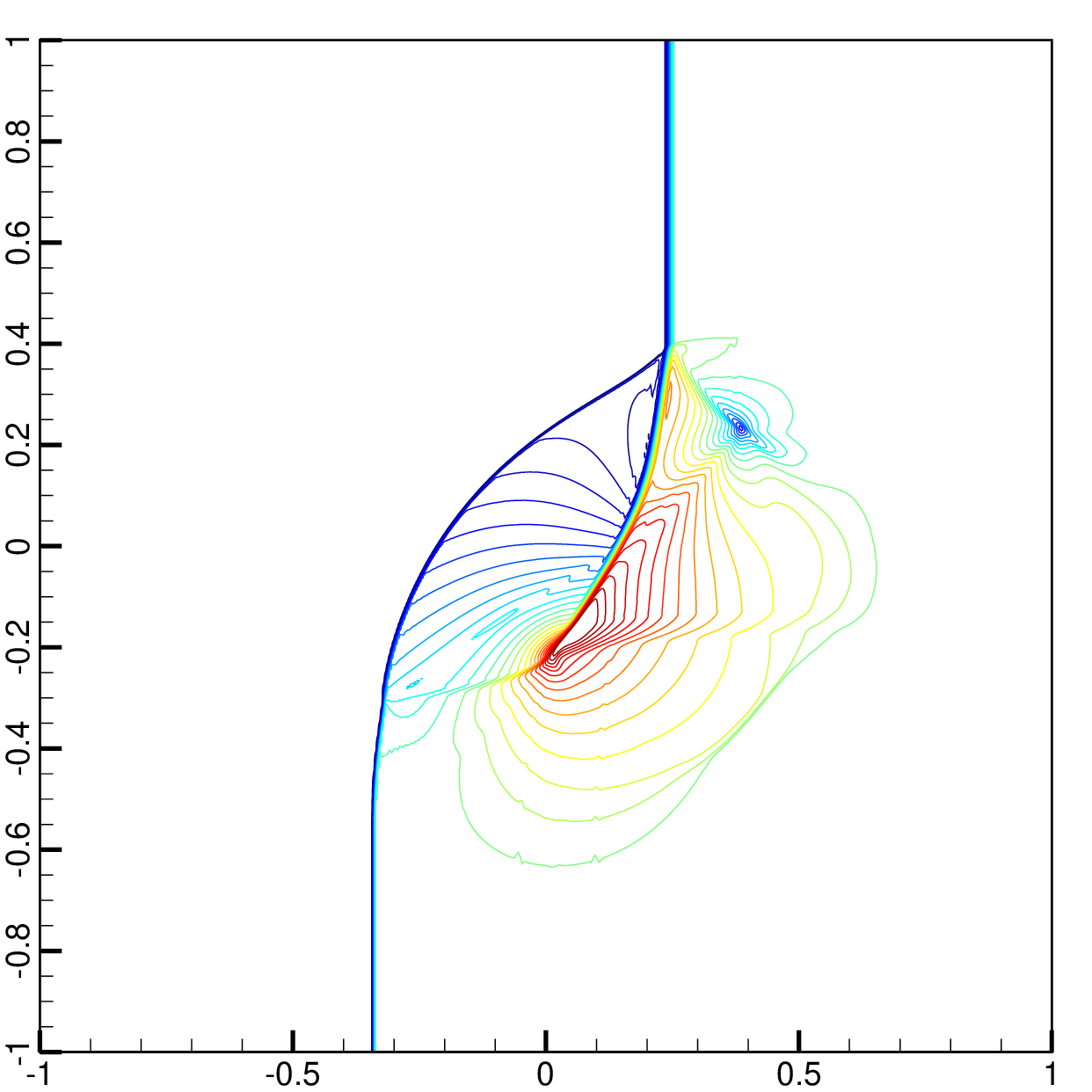}
		\caption{\(B_2\), \(512 \times 512\) grid}
	\end{subfigure}

	\caption{Contour plots of the Yee--Sj\"ogreen Riemann problem. 30 equally spaced contours for \(\rho \in [0.95, 1.85]\), \(B_1 \in [0.37, 1.24]\) and \(B_2 \in [0.55, 1.41]\).}
	\label{fig: Yee-2D-Riemann}
\end{figure}

Consider the 2D Riemann problem studied by Yee and Sj\"ogreen \cite{10.1007/3-540-27170-8_42}. The initial condition in conservative variables is listed in \Cref{tb: Yee 2D Riemann initial}. Extrapolation boundary conditions are used. We use \(256 \times 256\) and \(512 \times 512\) grids for comparison. Plots are shown in \Cref{fig: Yee-2D-Riemann}. The results are comparable to \cite{10.1006/jcph.2001.6961, 10.1007/s10915-005-9004-5, 10.1016/j.jcp.2025.113911}.

\begin{table}{htb}
	\centering
	\caption{Initial condition for the Yee--Sj\"ogreen 2D Riemann problem}
	\label{tb: Yee 2D Riemann initial}
	\begin{tabular}{c cccccccc}
		\hline
		region & \(\rho\) & \(m_1\) & \(m_2\) & \(m_3\) & \(B_1\) & \(B_2\) & \(B_3\) & \(E\) \\
		\hline
		\(x>0, y>0\) & 0.9308 & 1.4557 & -0.4633 & 0.0575 & 0.3501 & 0.9830 & 0.3050 & 5.0838 \\
		\(x<0, y>0\) & 1.0304 & 1.5774 & -1.0455 & -0.1016 & 0.3501 & 0.5078 & 0.1576 & 5.7813 \\
		\(x<0, y<0\) & 1.0000 & 1.7500 & -1.0000 & 0.0000 & 0.5642 & 0.5078 & 0.2539 & 6.0000 \\
		\(x>0, y<0\) & 1.8887 & 0.2334 & -1.7422 & 0.0733 & 0.5642 & 0.9830 & 0.4915 & 12.999 \\
		\hline
	\end{tabular}
\end{table}

\subsection{Orszag--Tang vortex}

\begin{figure}[htbp]
	\centering

	\begin{subfigure}{0.3\linewidth}
		\centering
		\includegraphics[width=\linewidth]{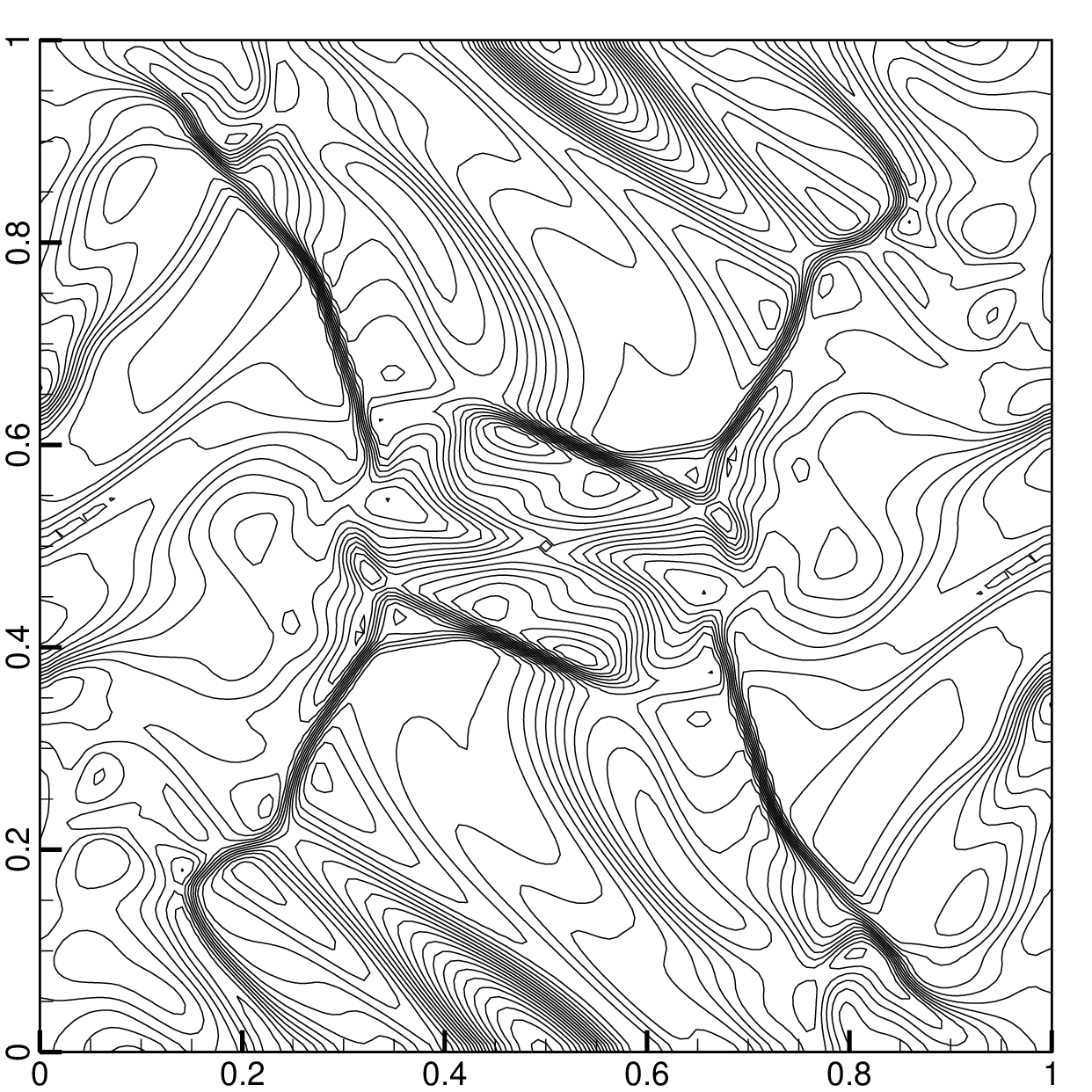}
		\caption{\(k=1\), \(N=128\)}
	\end{subfigure}
	\begin{subfigure}{0.3\linewidth}
		\centering
		\includegraphics[width=\linewidth]{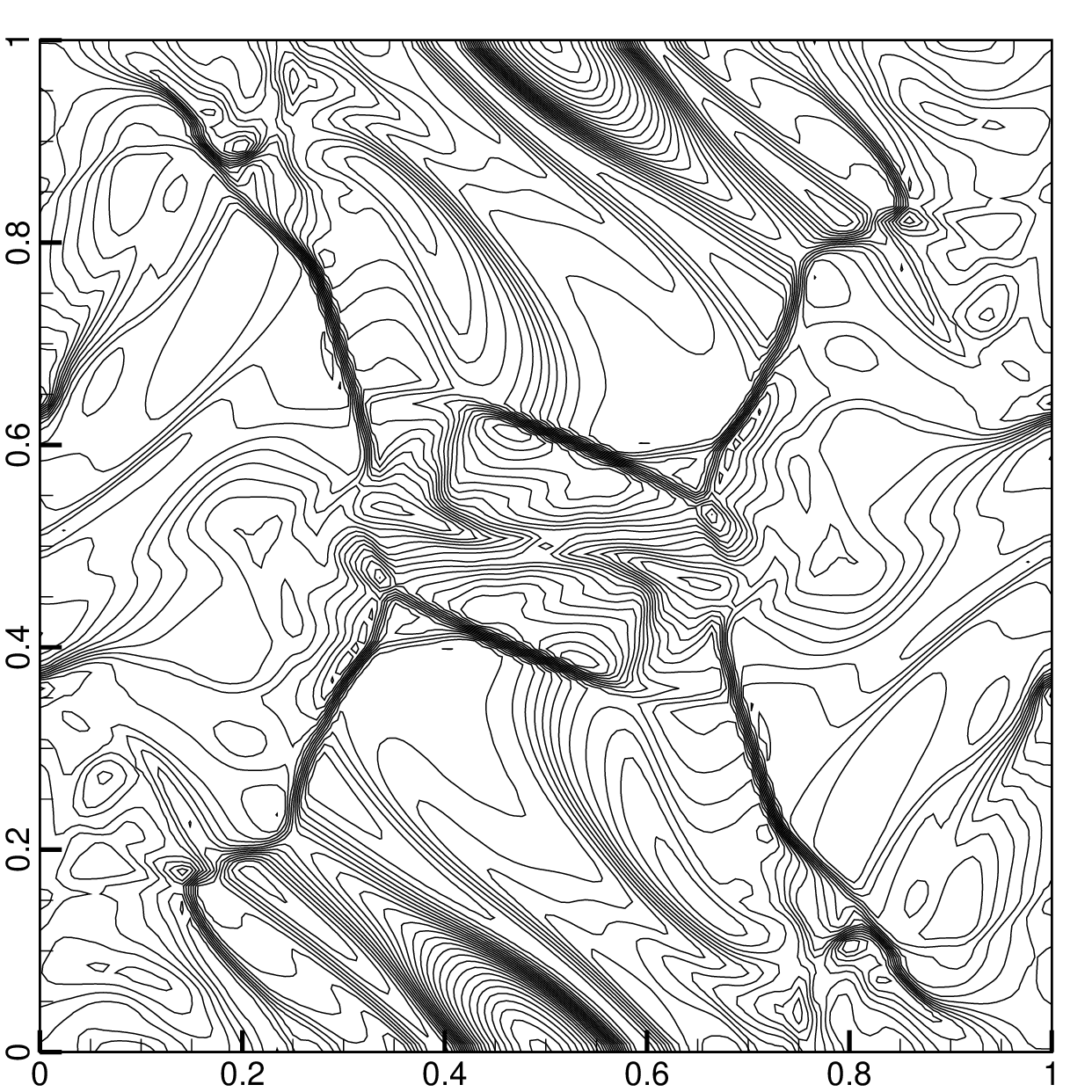}
		\caption{\(k=2\), \(N=128\)}
	\end{subfigure}
	\begin{subfigure}{0.3\linewidth}
		\centering
		\includegraphics[width=\linewidth]{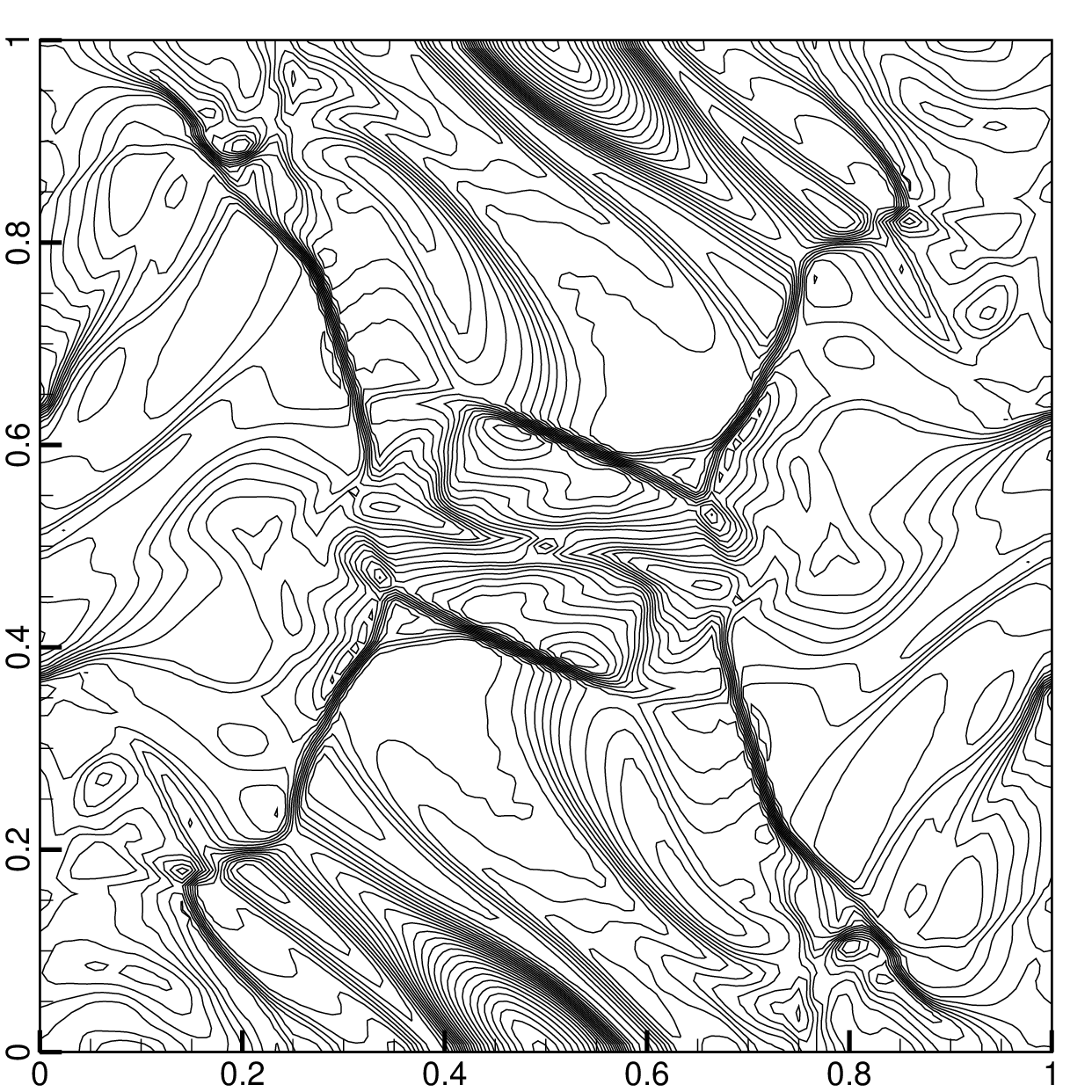}
		\caption{\(k=3\), \(N=128\)}
	\end{subfigure}

	\begin{subfigure}{0.3\linewidth}
		\centering
		\includegraphics[width=\linewidth]{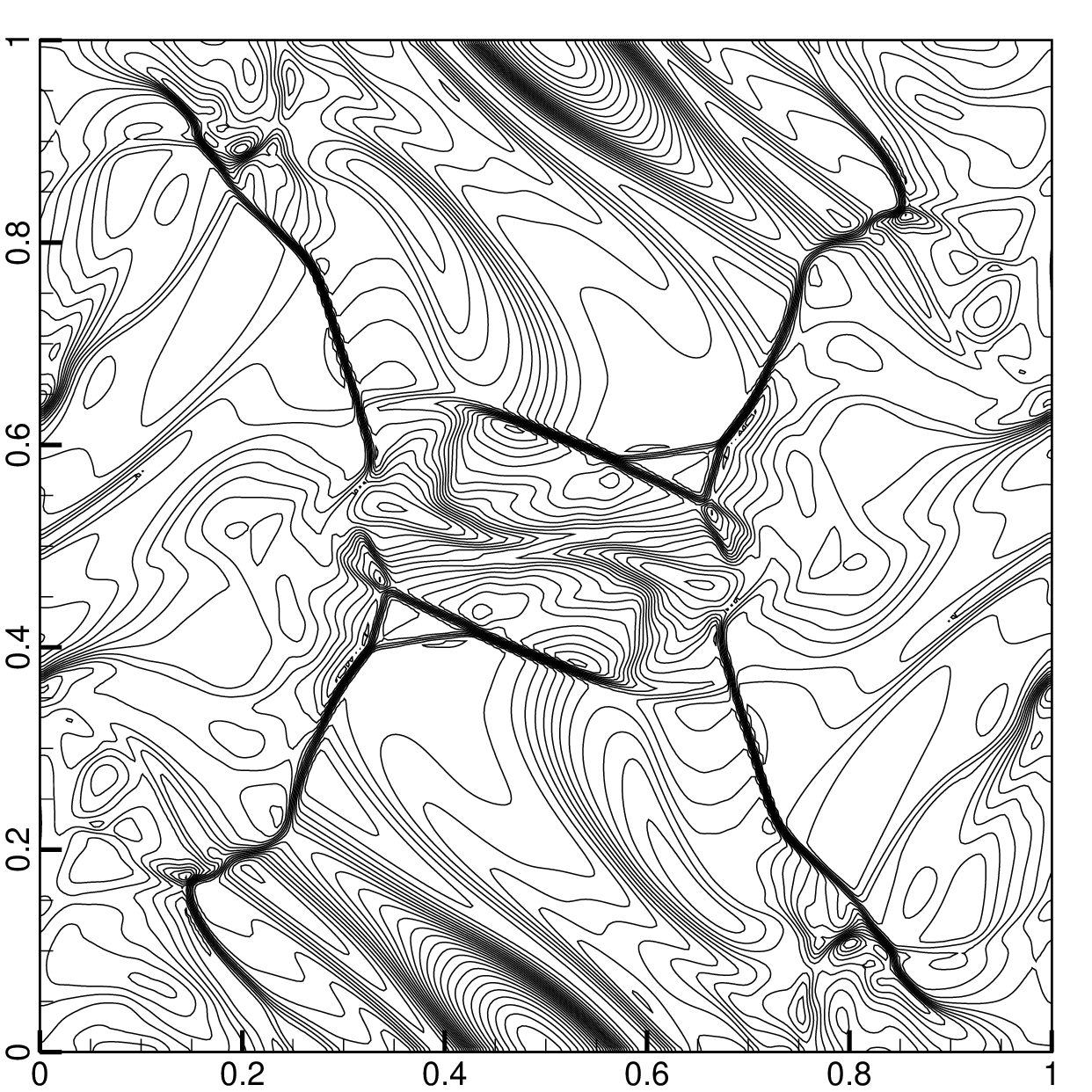}
		\caption{\(k=1\), \(N=256\)}
	\end{subfigure}
	\begin{subfigure}{0.3\linewidth}
		\centering
		\includegraphics[width=\linewidth]{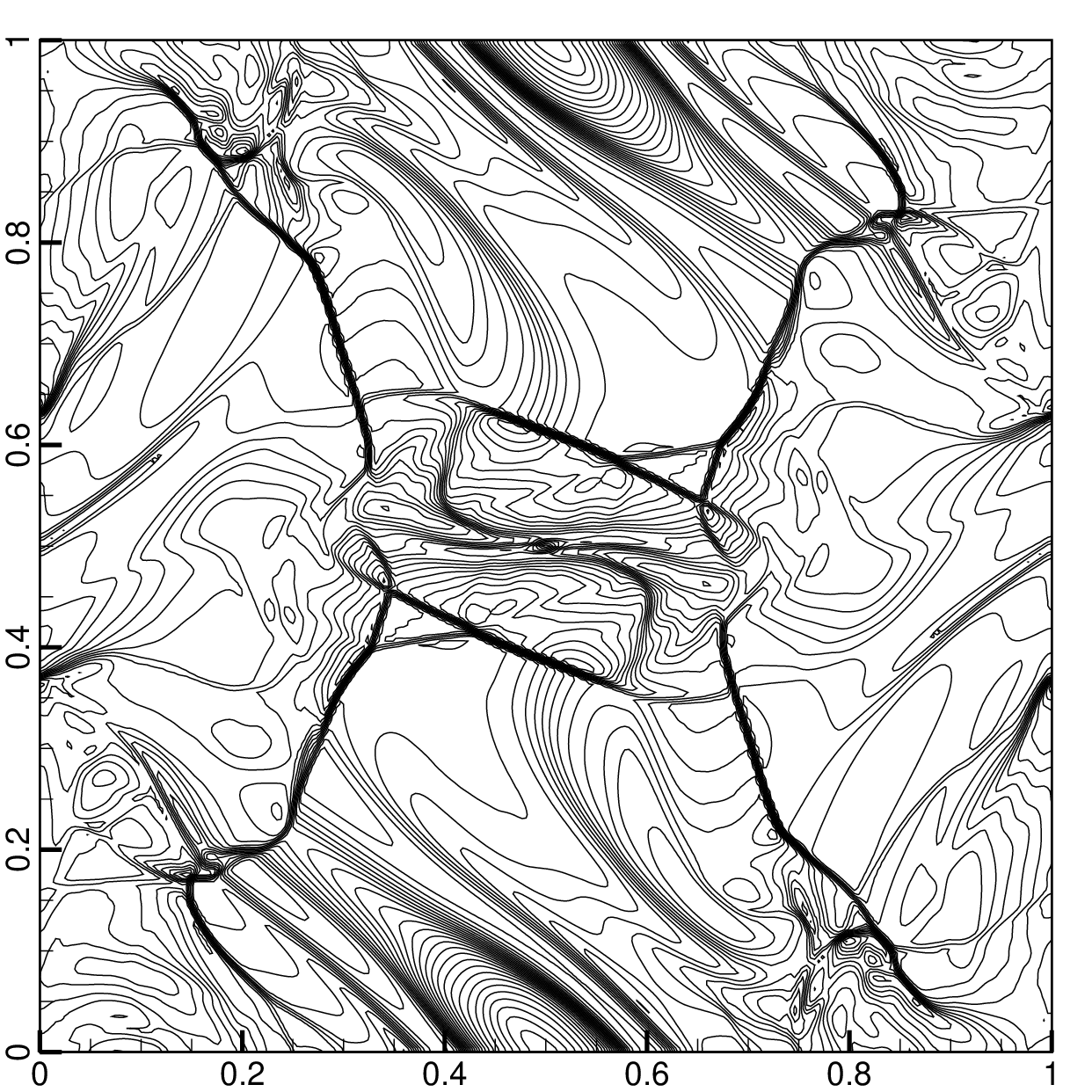}
		\caption{\(k=2\), \(N=256\)}
	\end{subfigure}
	\begin{subfigure}{0.3\linewidth}
		\centering
		\includegraphics[width=\linewidth]{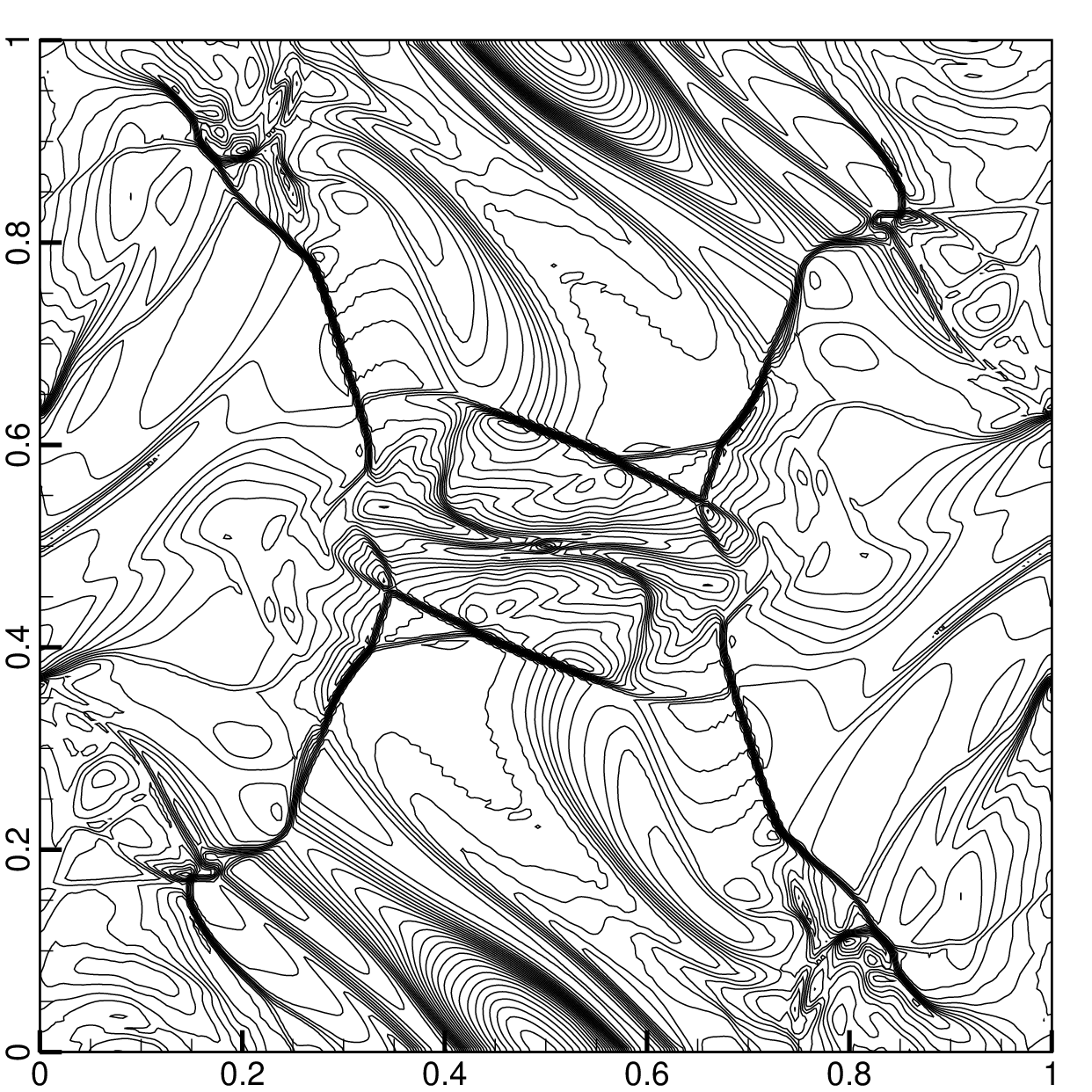}
		\caption{\(k=3\), \(N=256\)}
	\end{subfigure}

	\begin{subfigure}{0.3\linewidth}
		\centering
		\includegraphics[width=\linewidth]{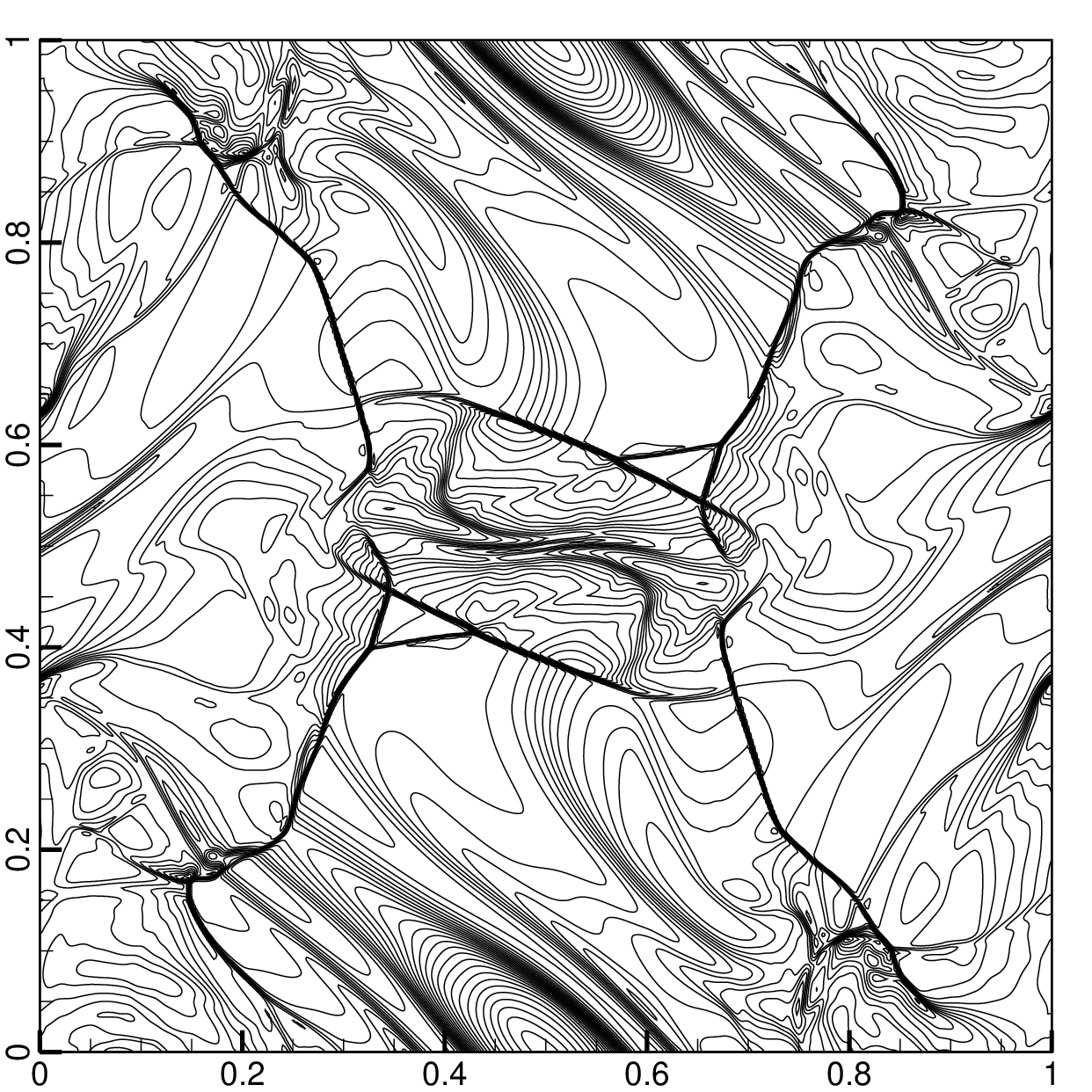}
		\caption{\(k=1\), \(N=512\)}
	\end{subfigure}
	\begin{subfigure}{0.3\linewidth}
		\centering
		\includegraphics[width=\linewidth]{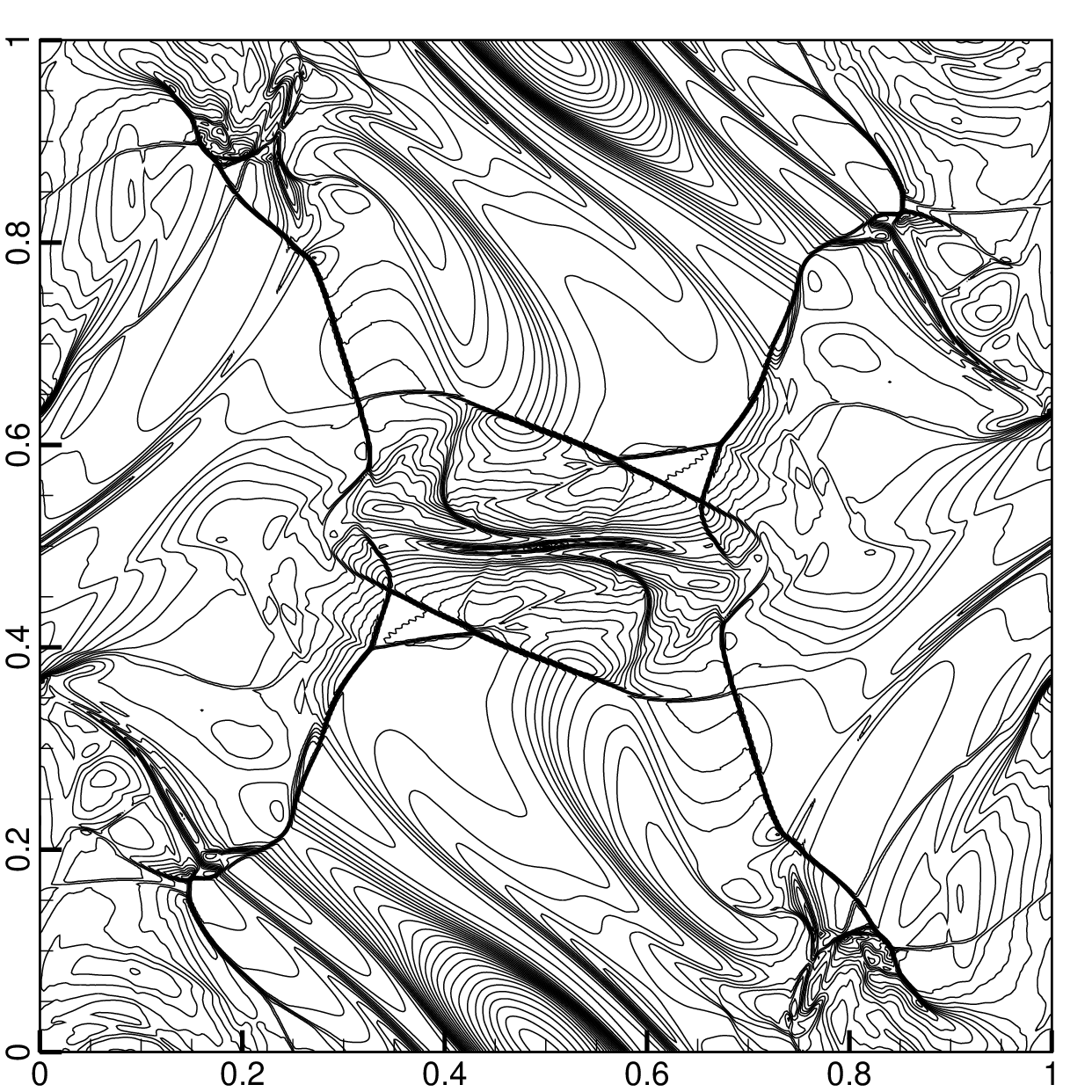}
		\caption{\(k=2\), \(N=512\)}
	\end{subfigure}
	\begin{subfigure}{0.3\linewidth}
		\centering
		\includegraphics[width=\linewidth]{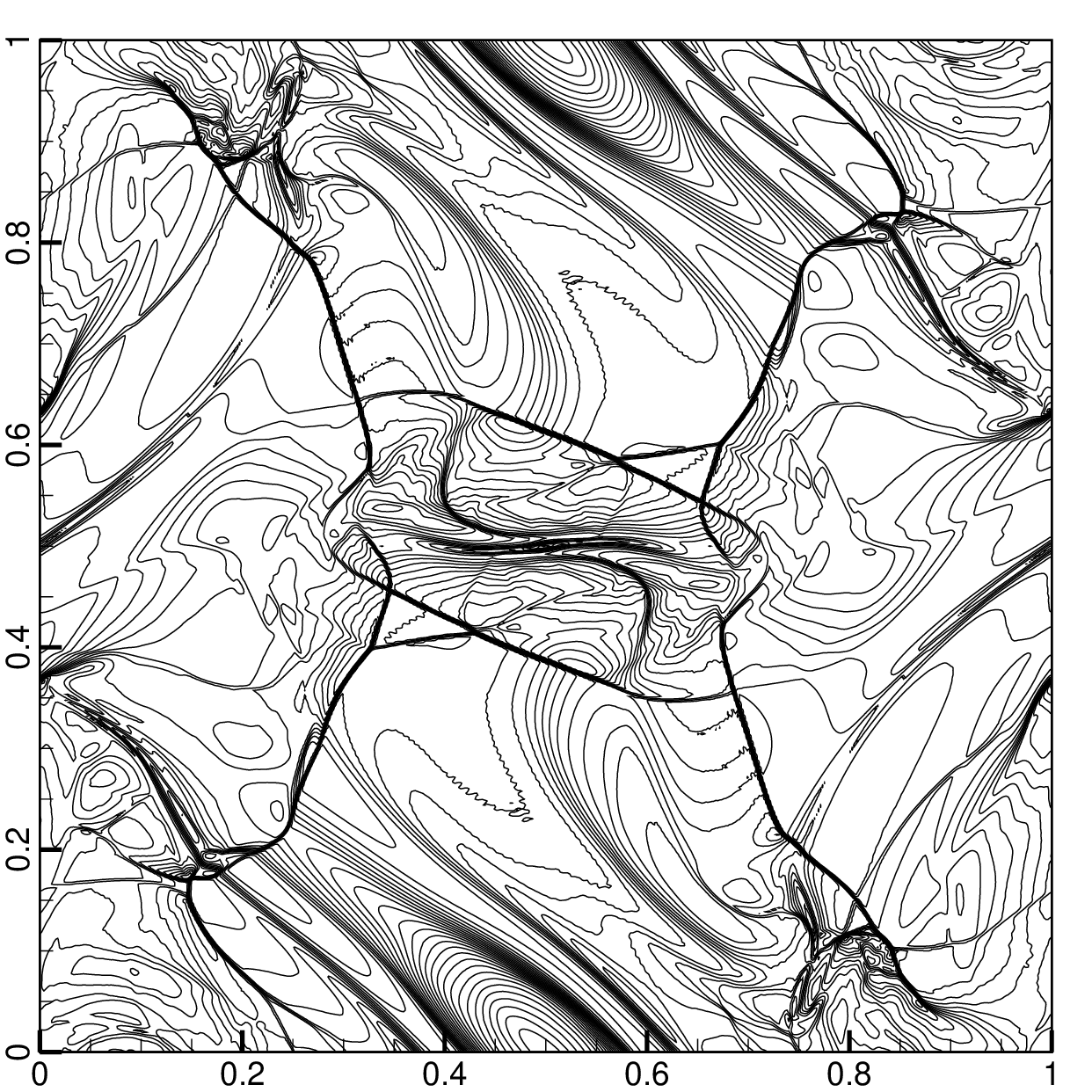}
		\caption{\(k=3\), \(N=512\)}
	\end{subfigure}

	\caption{Contour plots of Orszag--Tang vortex at \(T = 0.5\). 30 equally spaced contours for \(\rho \in [0.08, 0.5]\).}
	\label{fig: OZ rho}
\end{figure}

\begin{figure}[htbp]
	\centering
	\includegraphics[width=0.5\linewidth]{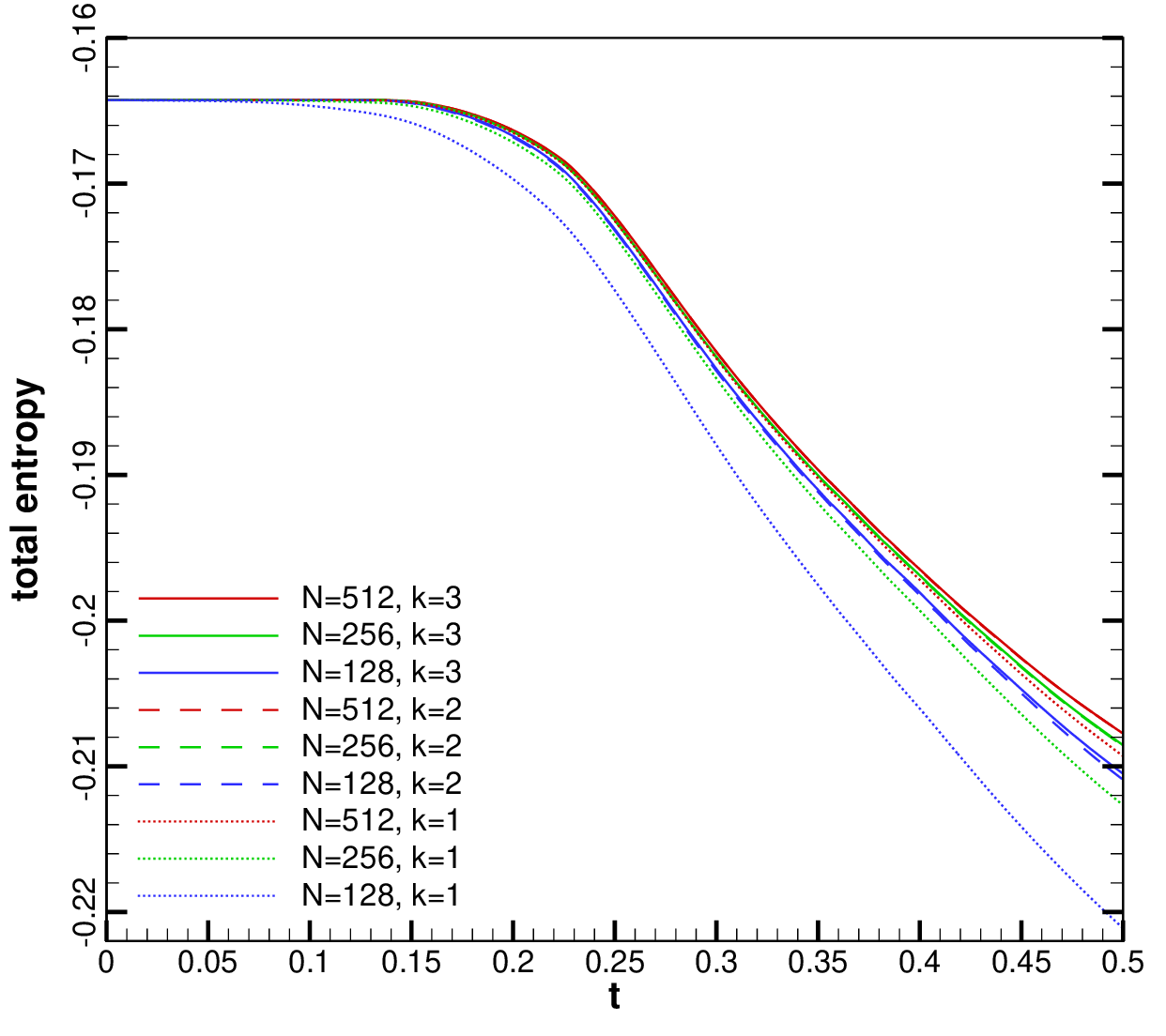}
	\caption{Evolution of entropy of Orszag--Tang vortex problem.}
	\label{fig: OZ entropy}
\end{figure}

The compressible Orszag--Tang vortex problem is a widely used 2D test problem, firstly studied by \cite{10.1063/1.859081}. The domain is \([0, 1]^2\), the final time is \(T = 0.5\), and the initial condition is given by the primitive variables
\begin{equation*}
	(\rho, u_1, u_2, u_3, B_1, B_2, B_3, p)
	= (\frac{25}{36\pi}, -\sin(2\pi y), \sin(2\pi x), 0, \frac{-\sin(2\pi y)}{\sqrt{4\pi}}, \frac{\sin(4\pi x)}{\sqrt{4\pi}}, 0, \frac{5}{12\pi}).
\end{equation*}
We use \(N \times N\) grids with \(N = 128, 256, 512\) and \(k = 1, 2, 3\) for comparison in \Cref{fig: OZ rho}. The results agree with \cite{10.2140/camcos.2021.16.59} well. \Cref{fig: OZ entropy} shows the evolution of total mathematical entropy. For all \(N\)'s and \(k\)'s, the total entropy is decreasing, which is consistent with entropy stability.

\subsection{Rotor}

\begin{figure}[htbp]
	\centering

	\begin{subfigure}{0.3\linewidth}
		\centering
		\includegraphics[width=\linewidth]{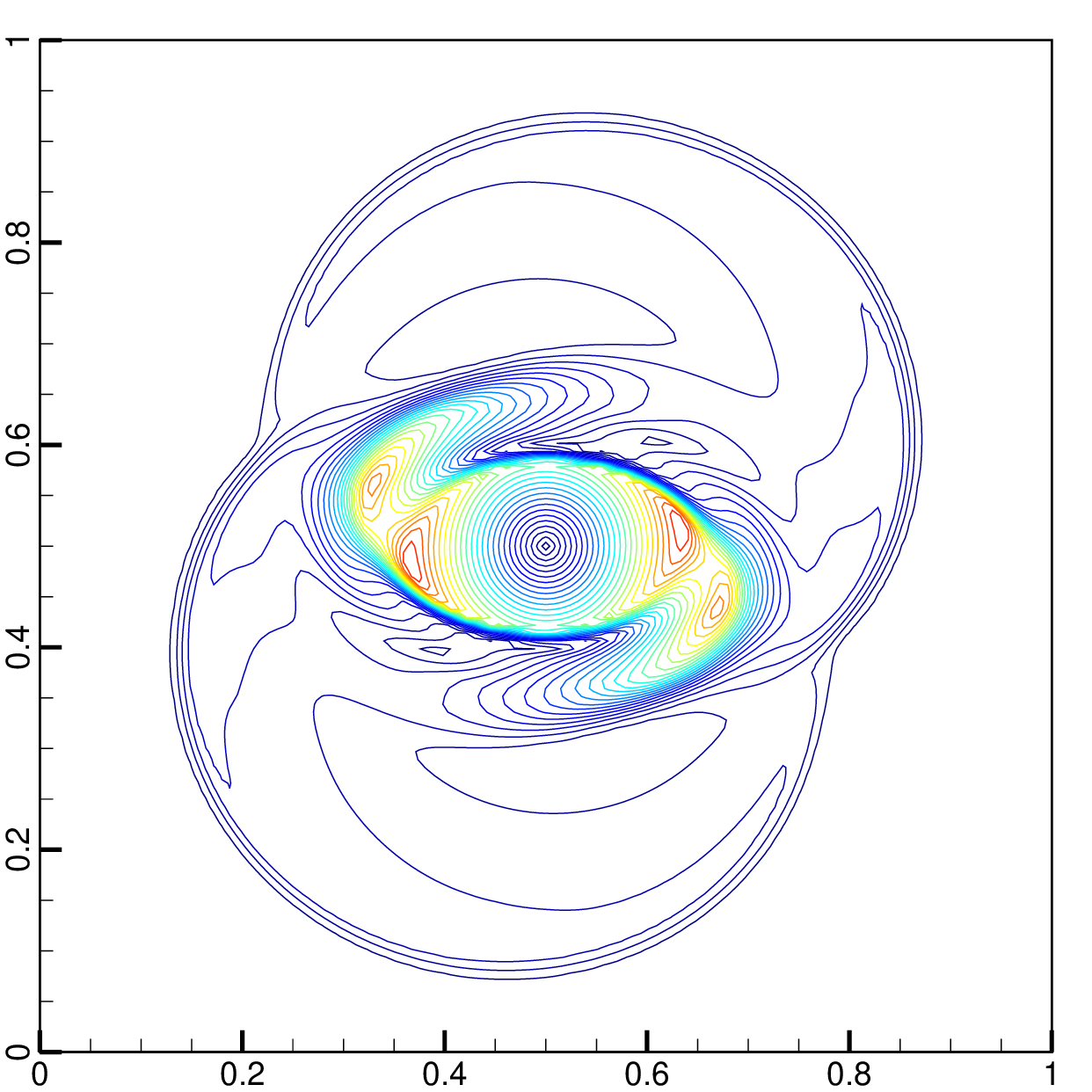}
		\caption{\(k=1\), \(N=128\)}
	\end{subfigure}
	\begin{subfigure}{0.3\linewidth}
		\centering
		\includegraphics[width=\linewidth]{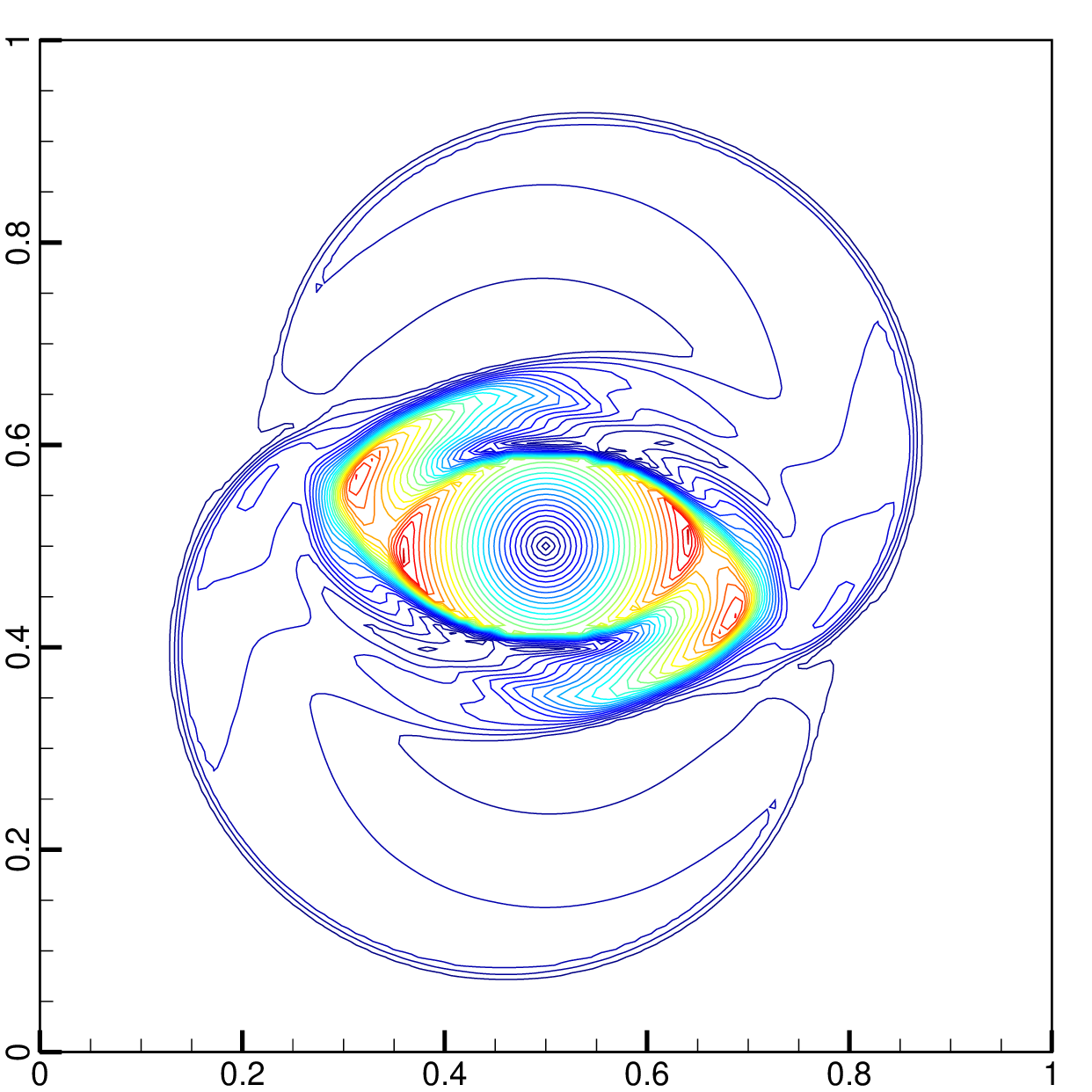}
		\caption{\(k=2\), \(N=128\)}
	\end{subfigure}
	\begin{subfigure}{0.3\linewidth}
		\centering
		\includegraphics[width=\linewidth]{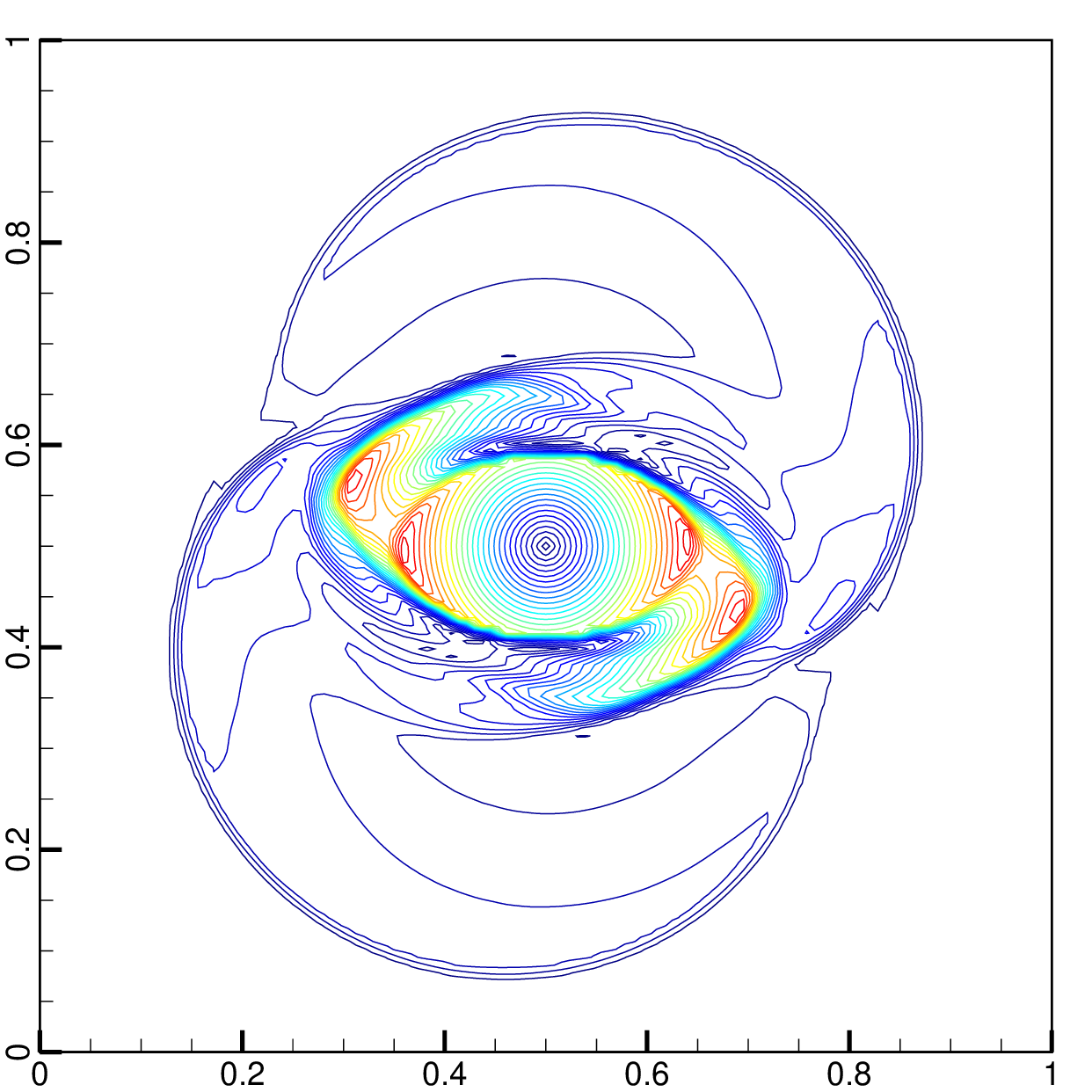}
		\caption{\(k=3\), \(N=128\)}
	\end{subfigure}

	\begin{subfigure}{0.3\linewidth}
		\centering
		\includegraphics[width=\linewidth]{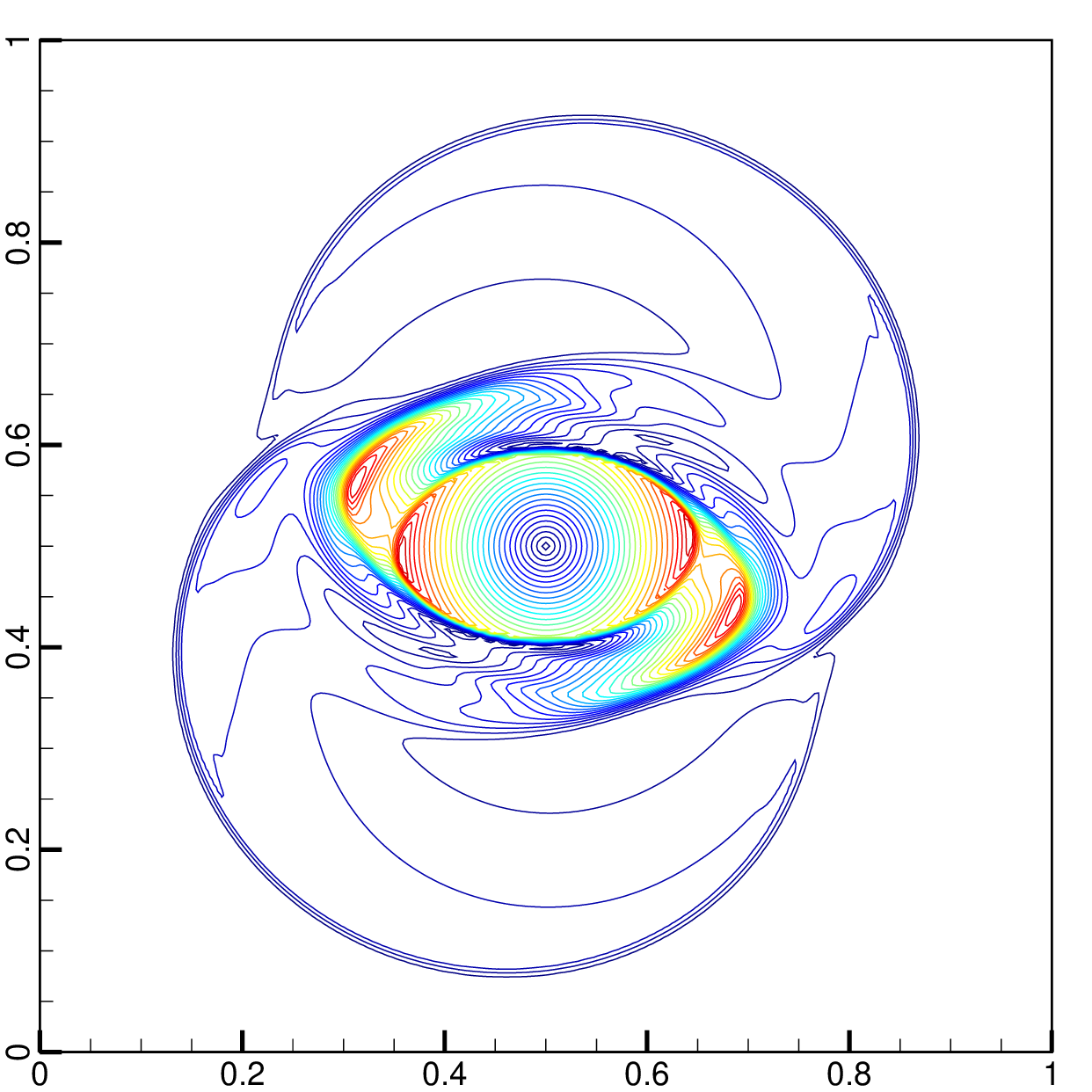}
		\caption{\(k=1\), \(N=256\)}
	\end{subfigure}
	\begin{subfigure}{0.3\linewidth}
		\centering
		\includegraphics[width=\linewidth]{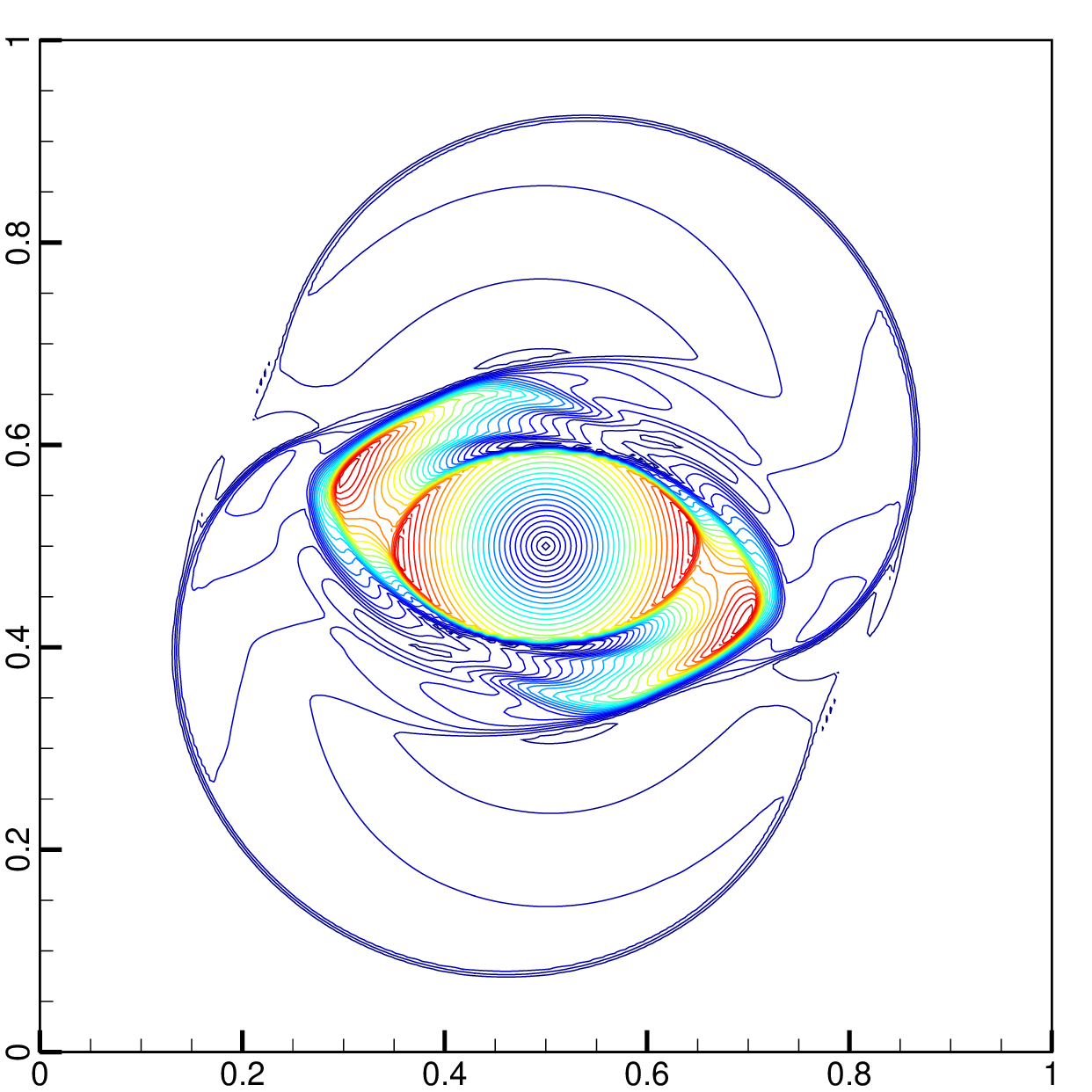}
		\caption{\(k=2\), \(N=256\)}
	\end{subfigure}
	\begin{subfigure}{0.3\linewidth}
		\centering
		\includegraphics[width=\linewidth]{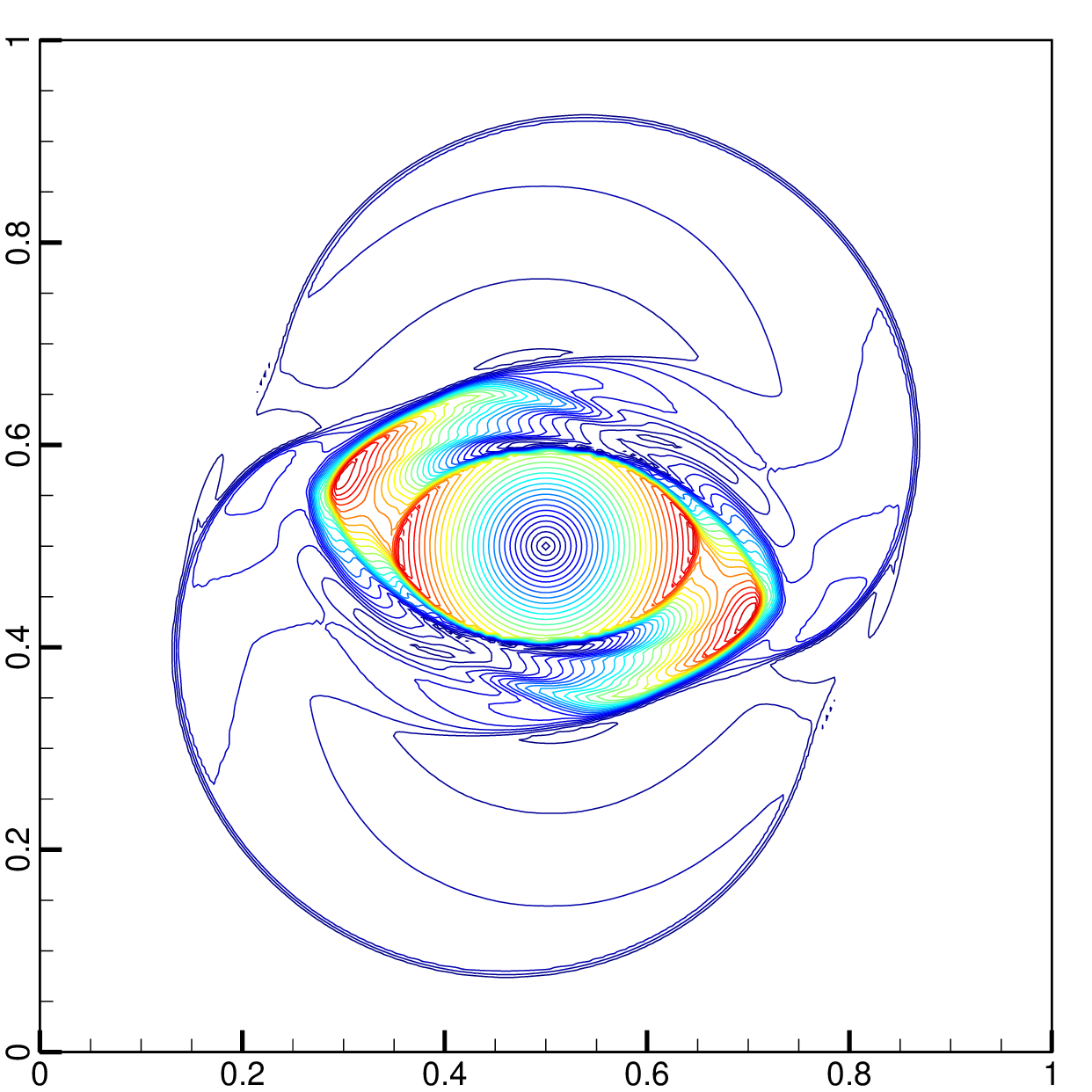}
		\caption{\(k=3\), \(N=256\)}
	\end{subfigure}

	\begin{subfigure}{0.3\linewidth}
		\centering
		\includegraphics[width=\linewidth]{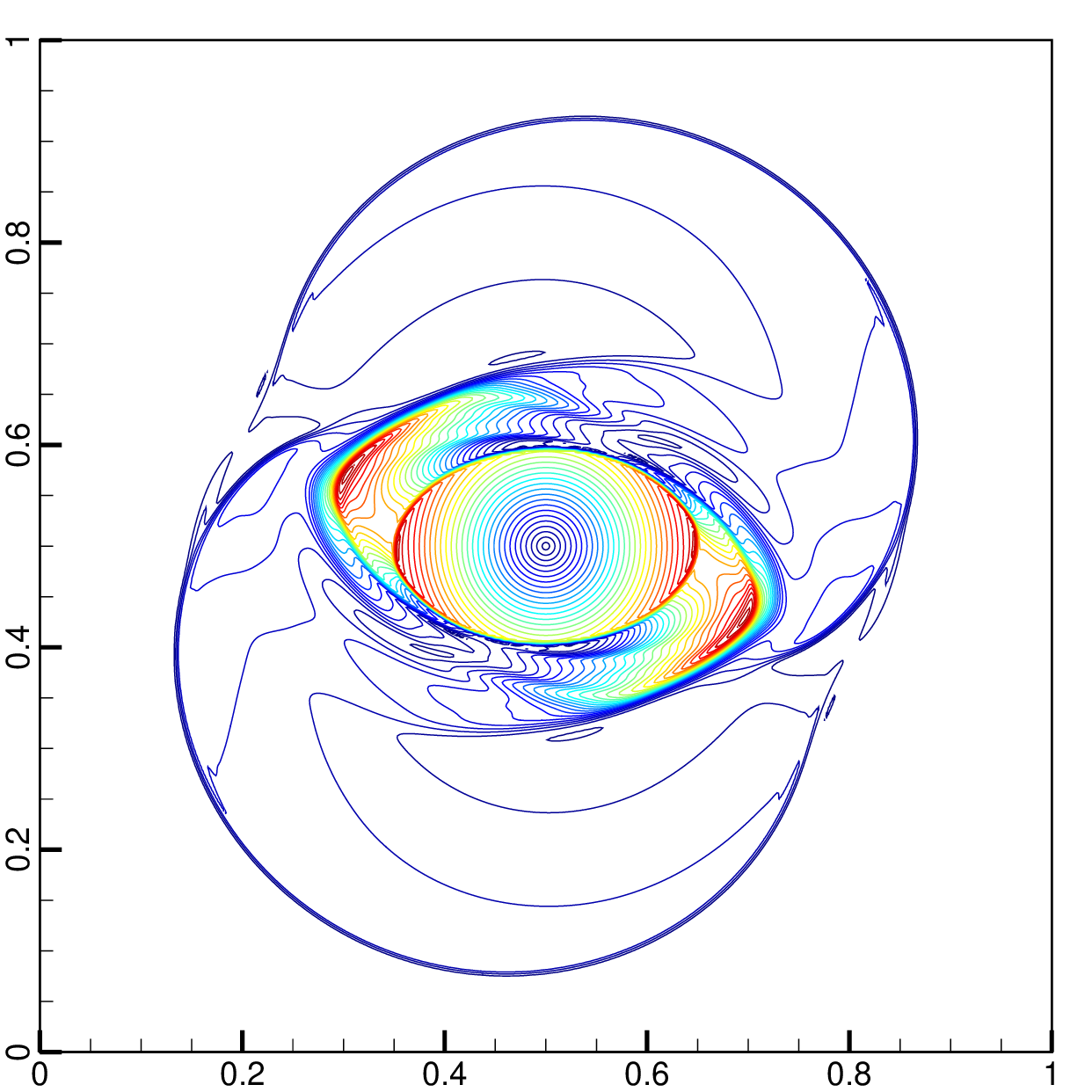}
		\caption{\(k=1\), \(N=512\)}
	\end{subfigure}
	\begin{subfigure}{0.3\linewidth}
		\centering
		\includegraphics[width=\linewidth]{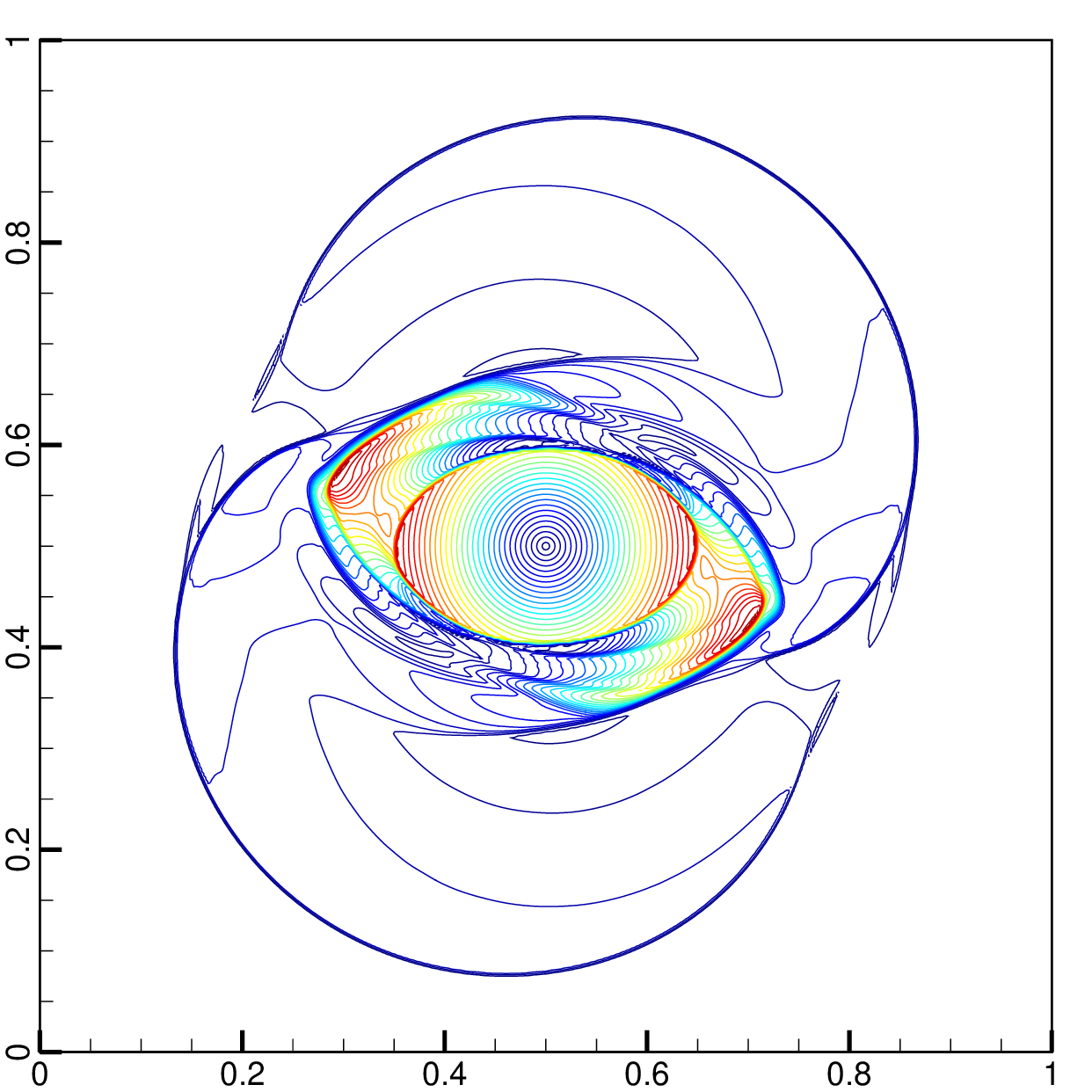}
		\caption{\(k=2\), \(N=512\)}
	\end{subfigure}
	\begin{subfigure}{0.3\linewidth}
		\centering
		\includegraphics[width=\linewidth]{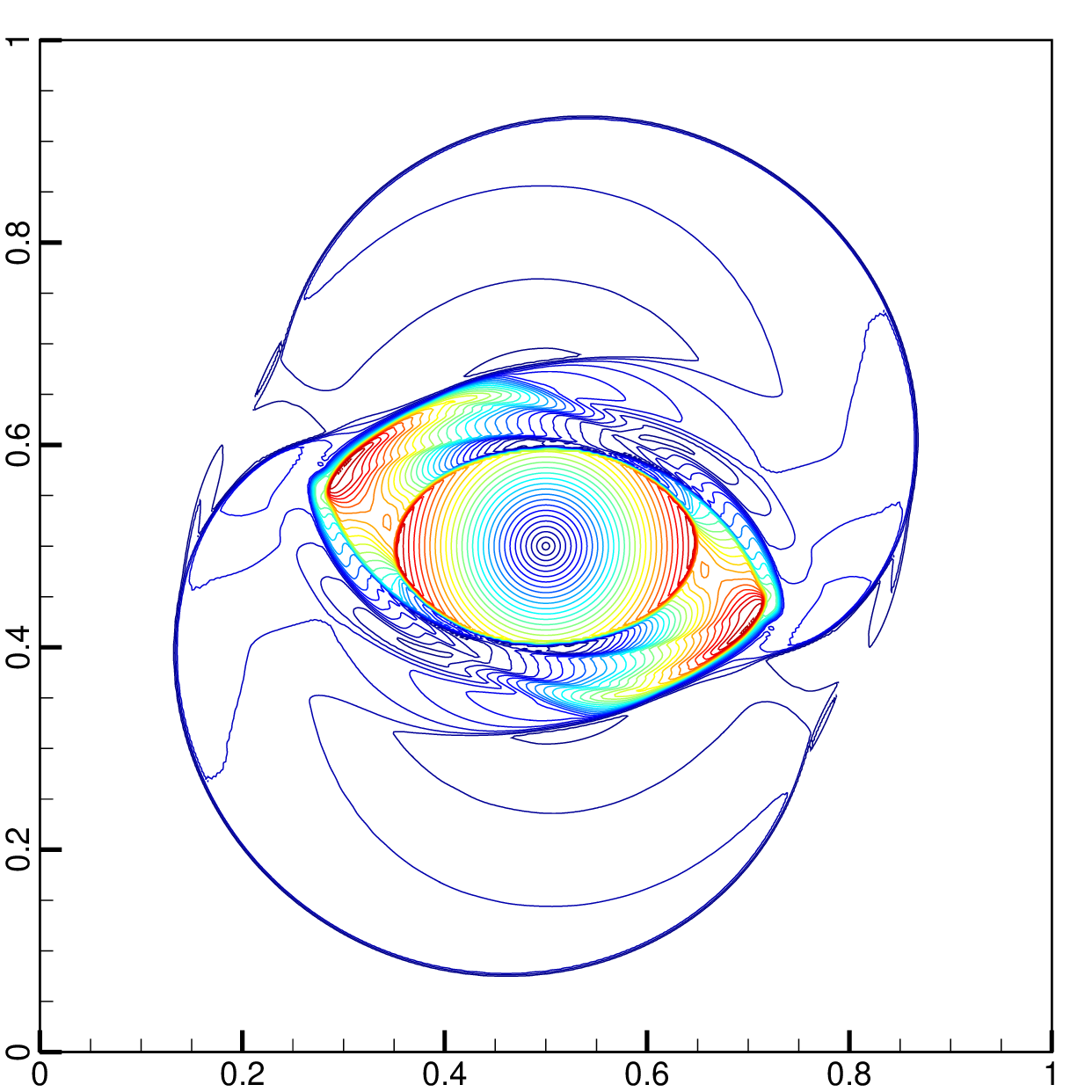}
		\caption{\(k=3\), \(N=512\)}
	\end{subfigure}

	\caption{Contour plots of Mach number for the first rotor problem at \(T = 0.15\). 30 equally spaced contours for \(\frac{\abs{\bm u}}{c} \in [0.1, 4.2]\).}
	\label{fig: RT mach}
\end{figure}

\begin{figure}[htbp]
	\centering
	\includegraphics[width=0.5\linewidth]{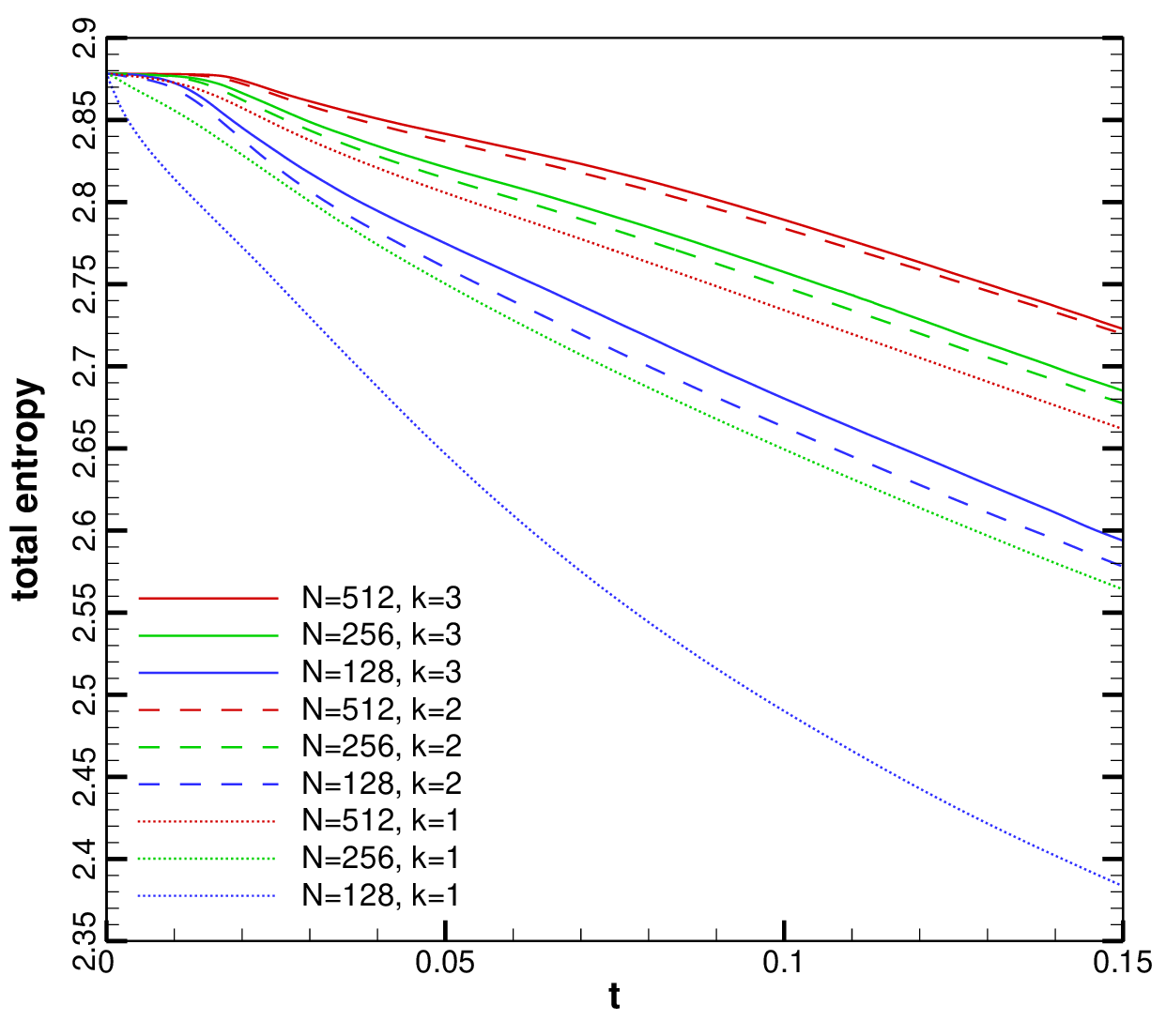}
	\caption{Evolution of entropy of the first rotor problem.}
	\label{fig: RT entropy}
\end{figure}

We consider the first rotor problem \cite{10.1006/jcph.2000.6519}, whose original version was proposed in \cite{10.1006/jcph.1998.6153}. A dense disc of fluid is rotating in the center of the region and the ambient fluid is static. A background magnetic field makes the central disc into oblate shape. The computation domain is \([0,1]^2\) with periodic boundary conditions, and the initial condition is
\begin{equation*}
	\begin{aligned}
		& r < r_0, \quad & \rho = 10, \quad \bm u = \frac{u_0}{r_0} (-(y-\frac{1}{2}), x - \frac{1}{2}, 0)^T, \\
		& r_0 < r < r_1, \quad & \rho = 1 + 9f, \quad \bm u = \frac{f u_0}{r_0} (-(y-\frac{1}{2}), x - \frac{1}{2}, 0)^T, \quad f = \frac{r_1 - r}{r_1 - r_0} \\
		& r > r_1, \quad & \rho = 1, \quad \bm u = (0, 0, 0)^T,
	\end{aligned}
\end{equation*}
and
\begin{equation*}
	p = 1, \quad \bm B = (5, 0, 0)^T.
\end{equation*}
Here, \(r = \sqrt{(x-\frac{1}{2})^2 + (y-\frac{1}{2})^2}\). The parameters are \(r_0 = 0.1\), \(r_1 = 0.115\), and \(u_0 = 2\). The adiabatic constant \(\gamma = 1.4\) is used. We compute the problem up to \(T = 0.15\). Mach plots are shown in \Cref{fig: RT mach} and the evolution of entropy is shown in \Cref{fig: RT entropy}. In this example, spurious distortions will develop at the corners of the oblate central region if the divergence of magnetic field is not well handled by the numerical scheme \cite{10.1007/s10915-004-4146-4}. We see from our results that no distortions occur, and our results agree with those in \cite{10.1007/s10915-004-4146-4, 10.1016/j.jcp.2017.10.043, 10.1007/s10915-020-01289-8, 10.2140/camcos.2021.16.59, 10.1016/j.jcp.2025.113911} well.

\subsection{Blast wave}

\begin{figure}[htbp]
	\centering

	\begin{subfigure}{0.3\linewidth}
		\centering
		\includegraphics[width=\linewidth]{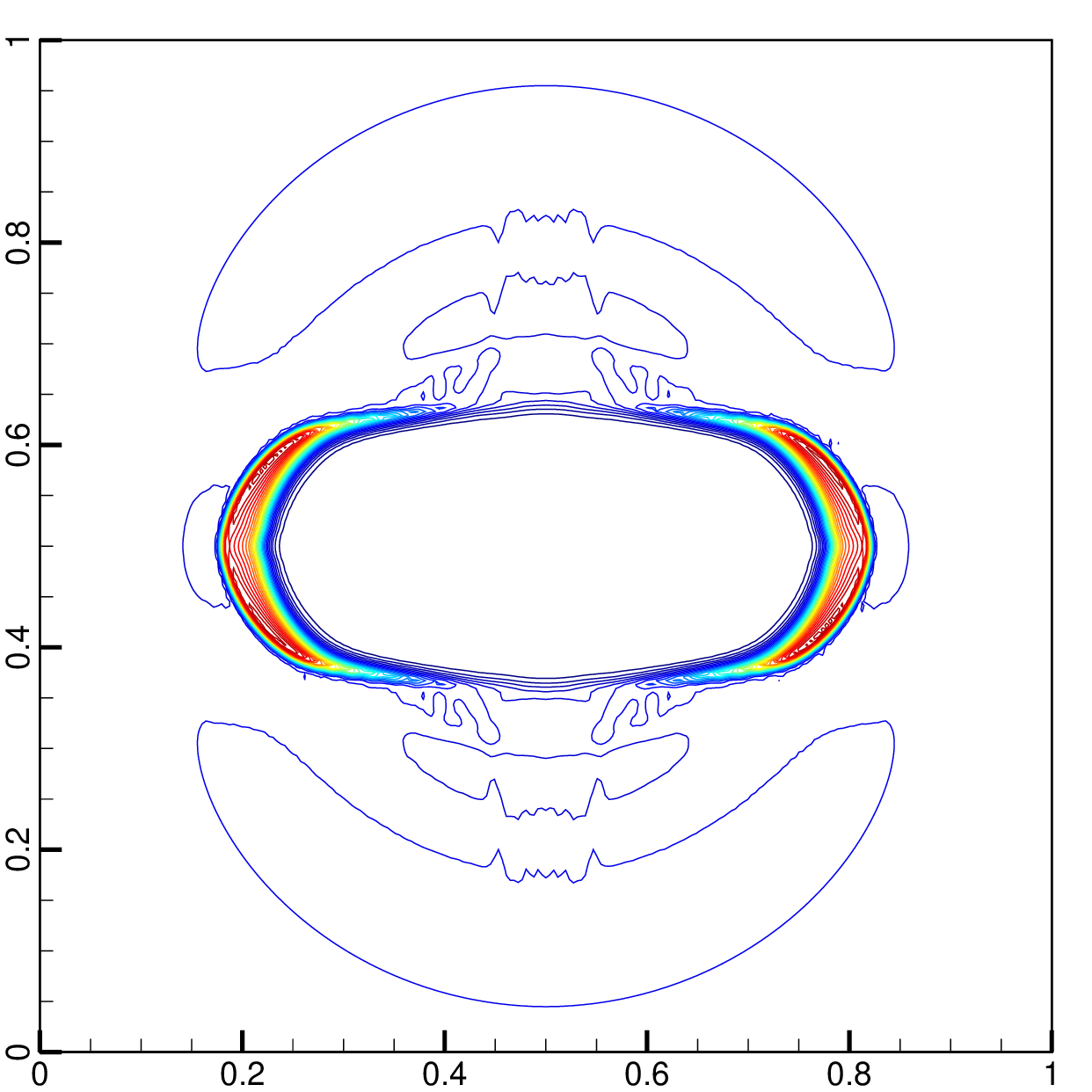}
		\caption{\(\rho\)}
	\end{subfigure}
	\begin{subfigure}{0.3\linewidth}
		\centering
		\includegraphics[width=\linewidth]{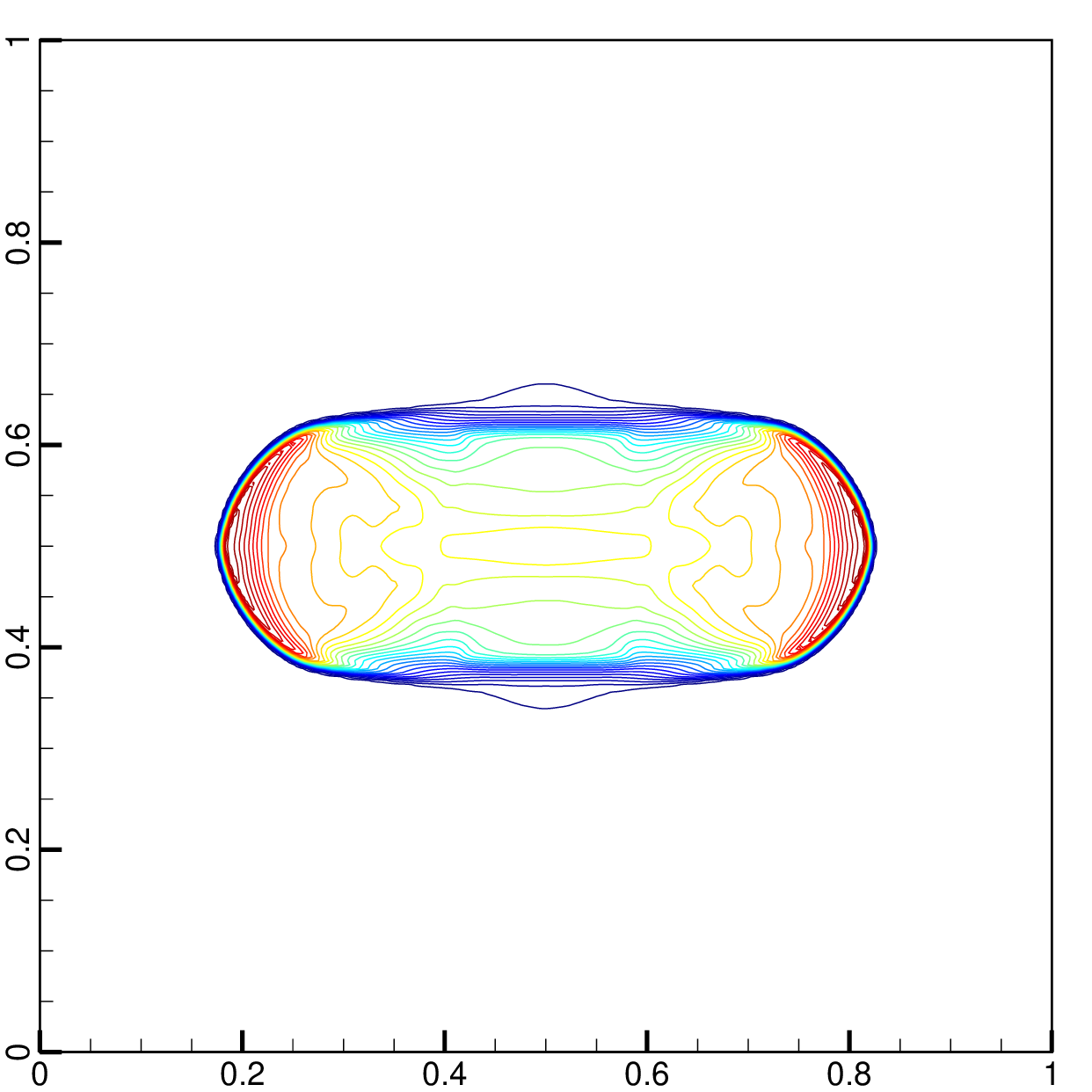}
		\caption{\(p\)}
	\end{subfigure}

	\begin{subfigure}{0.3\linewidth}
		\centering
		\includegraphics[width=\linewidth]{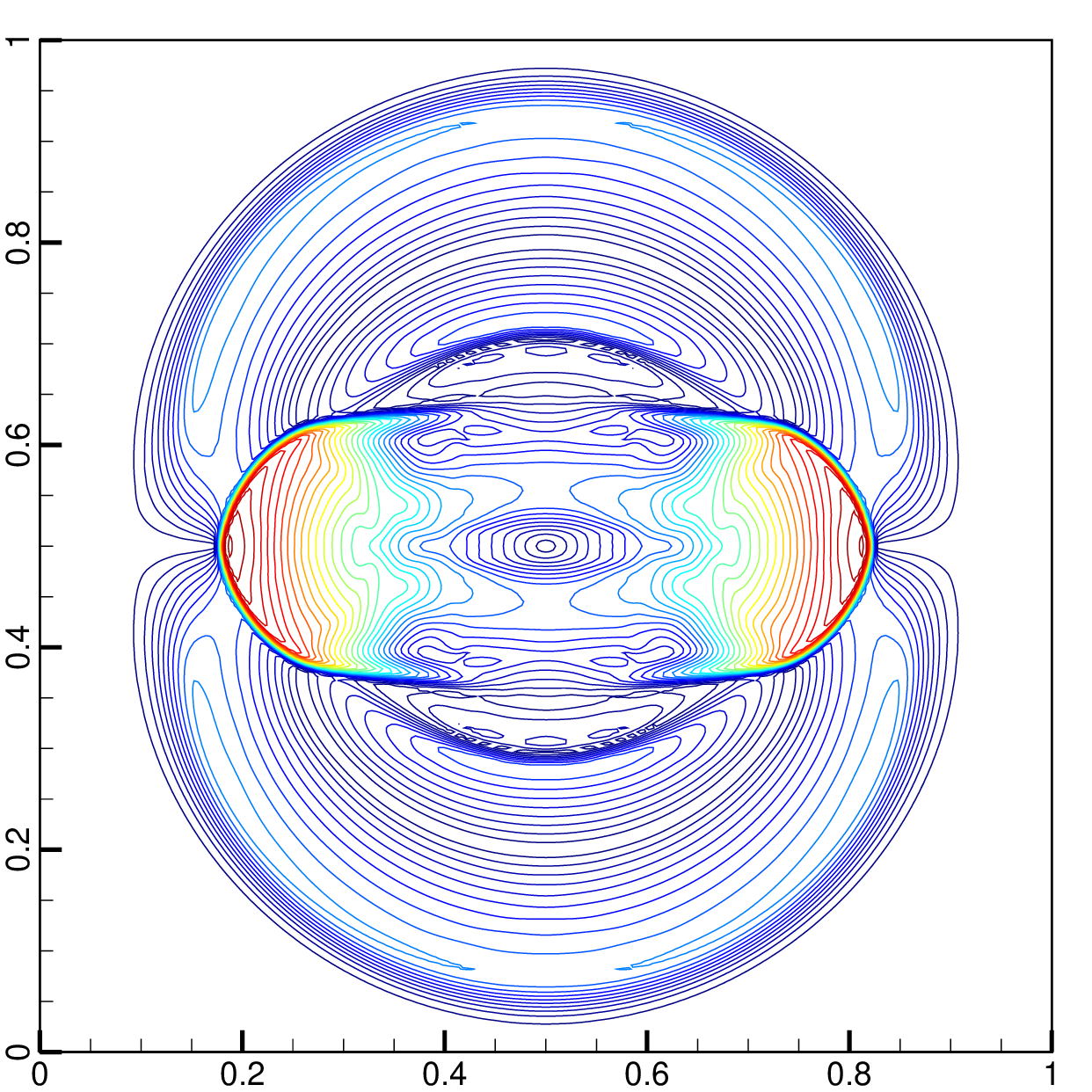}
		\caption{\(\abs{\bm u}\)}
	\end{subfigure}
	\begin{subfigure}{0.3\linewidth}
		\centering
		\includegraphics[width=\linewidth]{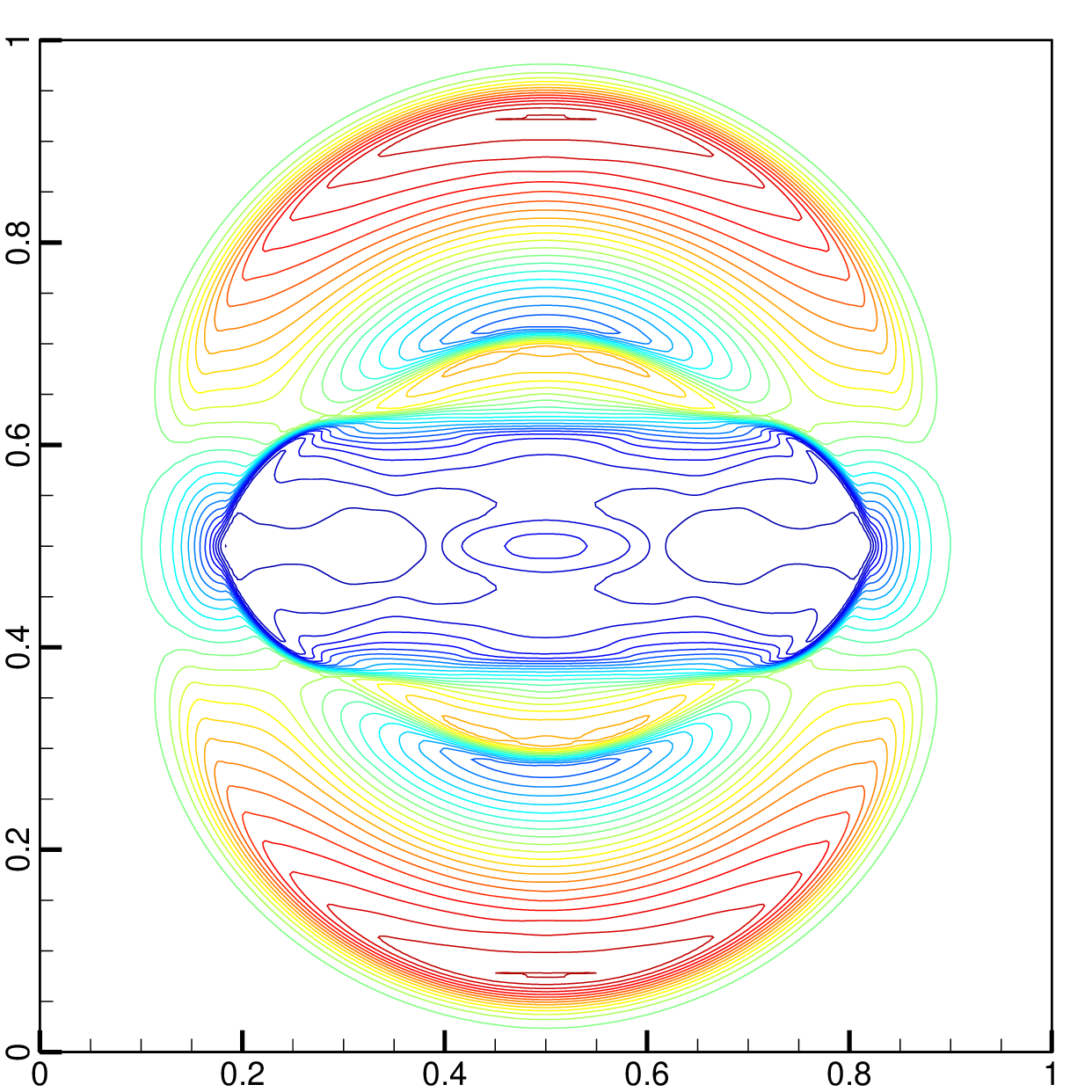}
		\caption{\(\frac{\abs{\bm B}^2}{2}\)}
	\end{subfigure}

	\caption{Contour plots of the blast wave problem. 30 equally spaced contours for \(\rho \in [0.45, 4.4]\), \(p \in [10, 250]\), \(\abs{\bm u} \in [0.5, 16]\), and \(\frac{\abs{\bm B}^2}{2} \in [200, 600]\).}
	\label{fig: blast}
\end{figure}

This problem was first proposed in \cite{10.1006/jcph.1998.6153} to test the robustness of a numerical scheme in the vicinity of low plasma beta, i.e., the ratio between the thermodynamic pressure and the magnetic pressure, and has been widely used for verification. The domain is \([0,1]^2\) with periodic boundary conditions. The initial condition consists of constant states
\begin{equation*}
	\rho = 1, \quad \bm u = (0, 0, 0)^T, \quad \bm B = (100, 0, 0)^T,
\end{equation*}
and a pressure distribution
\begin{equation*}
	p = \begin{cases}
		1000, & r < 0.1, \\
		0.1, & r > 0.1,
	\end{cases}
\end{equation*}
where \(r = \sqrt{(x-\frac{1}{2})^2 + (y-\frac{1}{2})^2}\) is the distance from the centroid. \(\gamma = 1.4\) is used for this example. \(256 \times 256\) grids are used. The density, pressure, velocity magnitude, and magnetic pressure are shown in \Cref{fig: blast}. We found our proposed scheme is robust and agree with results of \cite{10.1016/j.jcp.2012.12.019, 10.1137/18M1168042, 10.2140/camcos.2021.16.59} well.

\subsection{Cloud-shock interaction}

\begin{figure}[htbp]
	\centering

	\begin{subfigure}{0.3\linewidth}
		\centering
		\includegraphics[width=\linewidth]{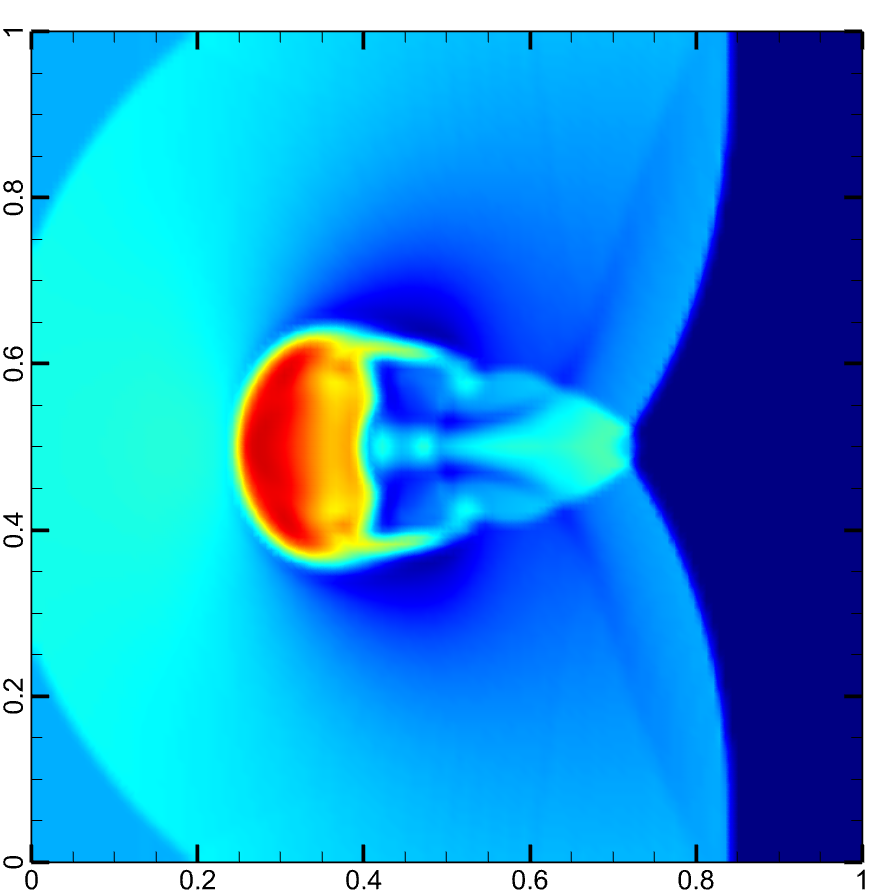}
		\caption{\(\ln\rho\), \(N=128\)}
	\end{subfigure}
	\begin{subfigure}{0.3\linewidth}
		\centering
		\includegraphics[width=\linewidth]{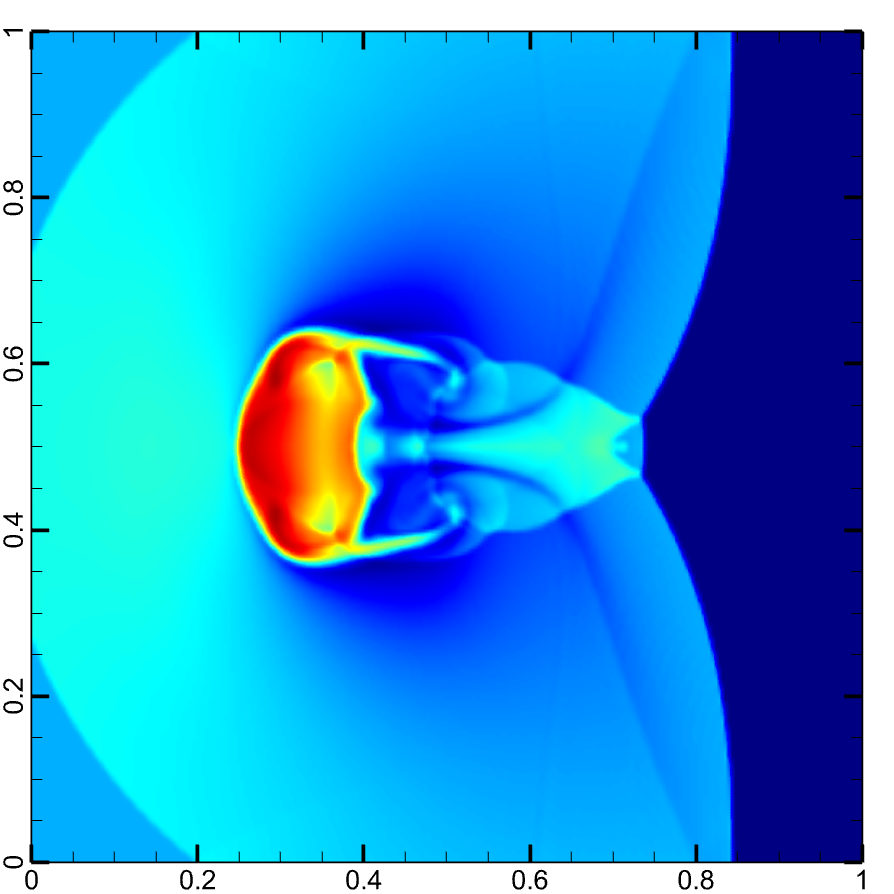}
		\caption{\(\ln\rho\), \(N=256\)}
	\end{subfigure}
	\begin{subfigure}{0.3\linewidth}
		\centering
		\includegraphics[width=\linewidth]{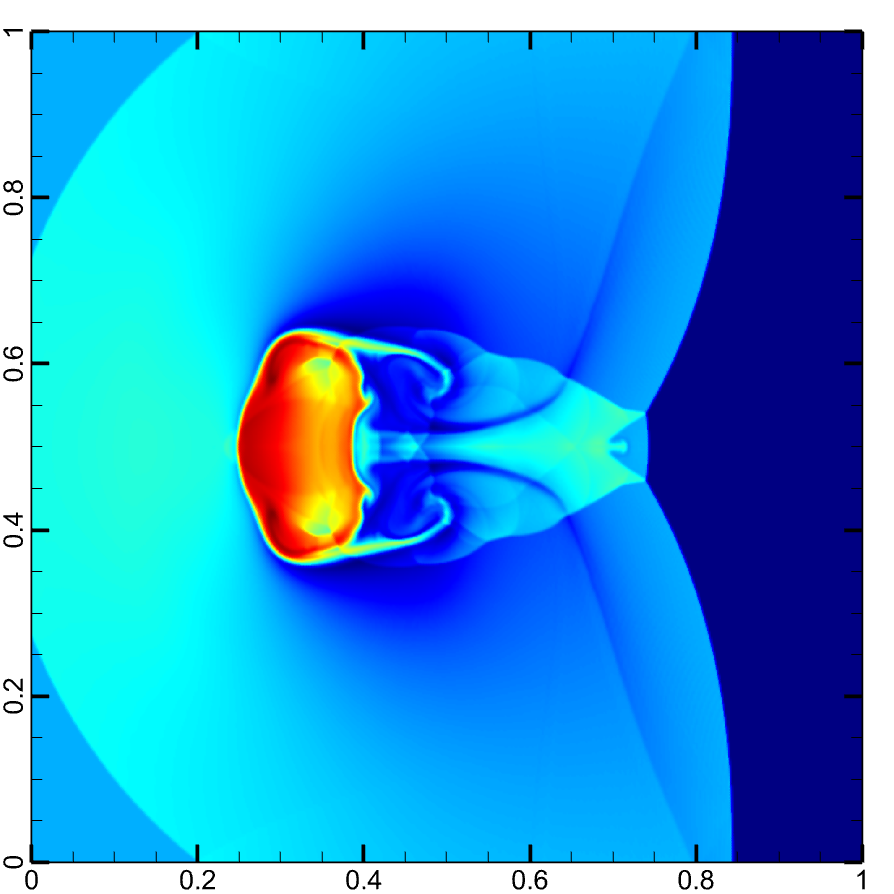}
		\caption{\(\ln\rho\), \(N=512\)}
	\end{subfigure}

	\begin{subfigure}{0.3\linewidth}
		\centering
		\includegraphics[width=\linewidth]{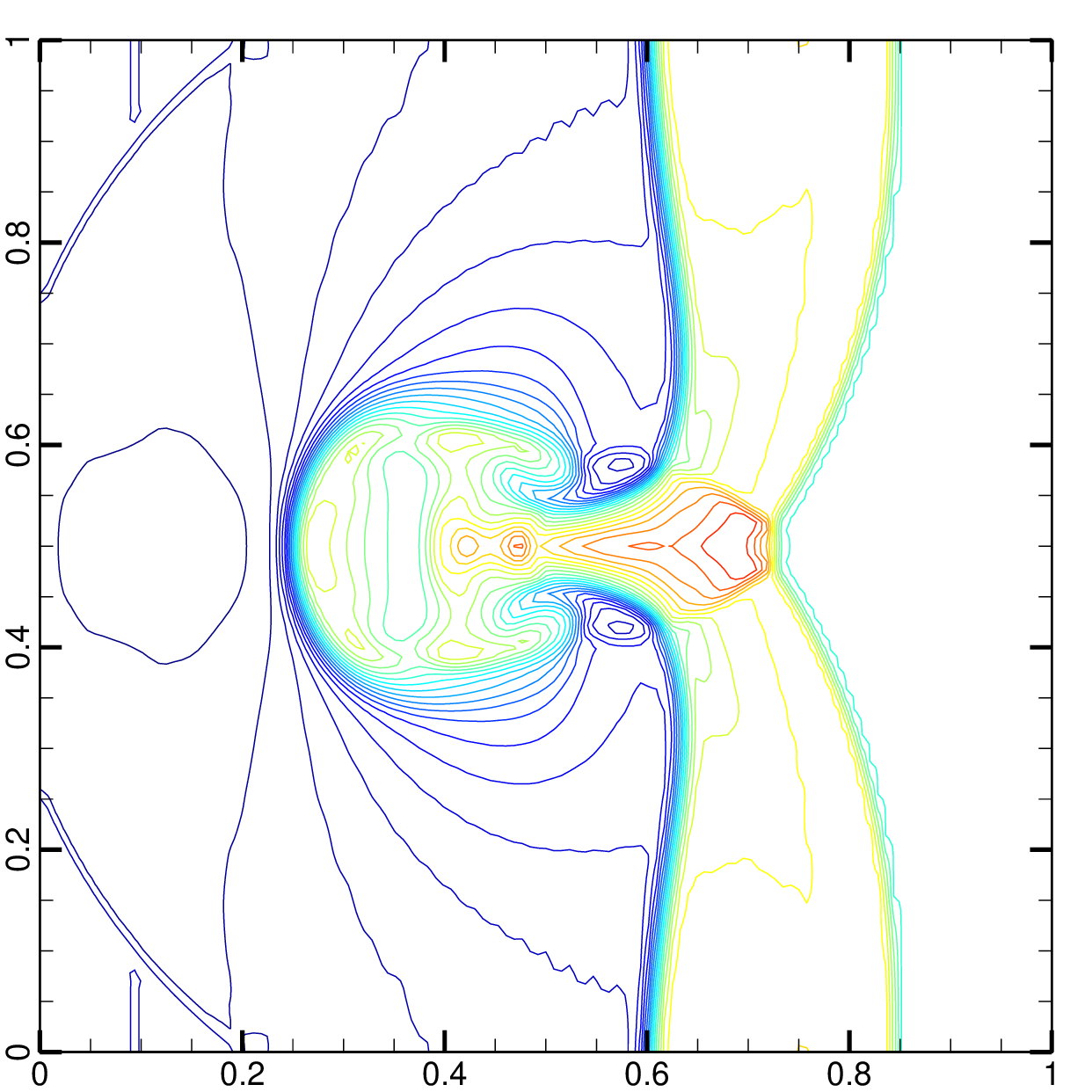}
		\caption{\(B_3\), \(N=128\)}
	\end{subfigure}
	\begin{subfigure}{0.3\linewidth}
		\centering
		\includegraphics[width=\linewidth]{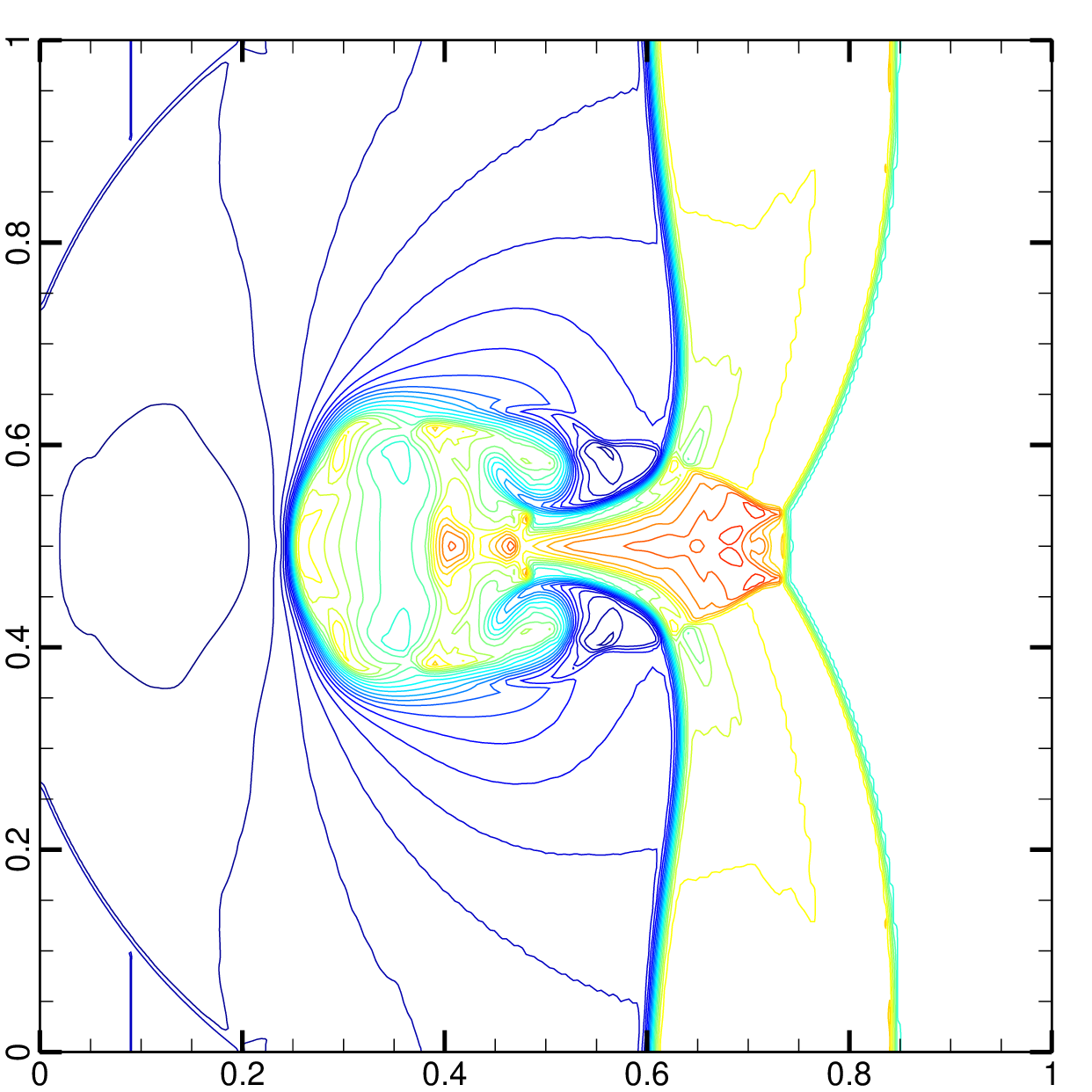}
		\caption{\(B_3\), \(N=256\)}
	\end{subfigure}
	\begin{subfigure}{0.3\linewidth}
		\centering
		\includegraphics[width=\linewidth]{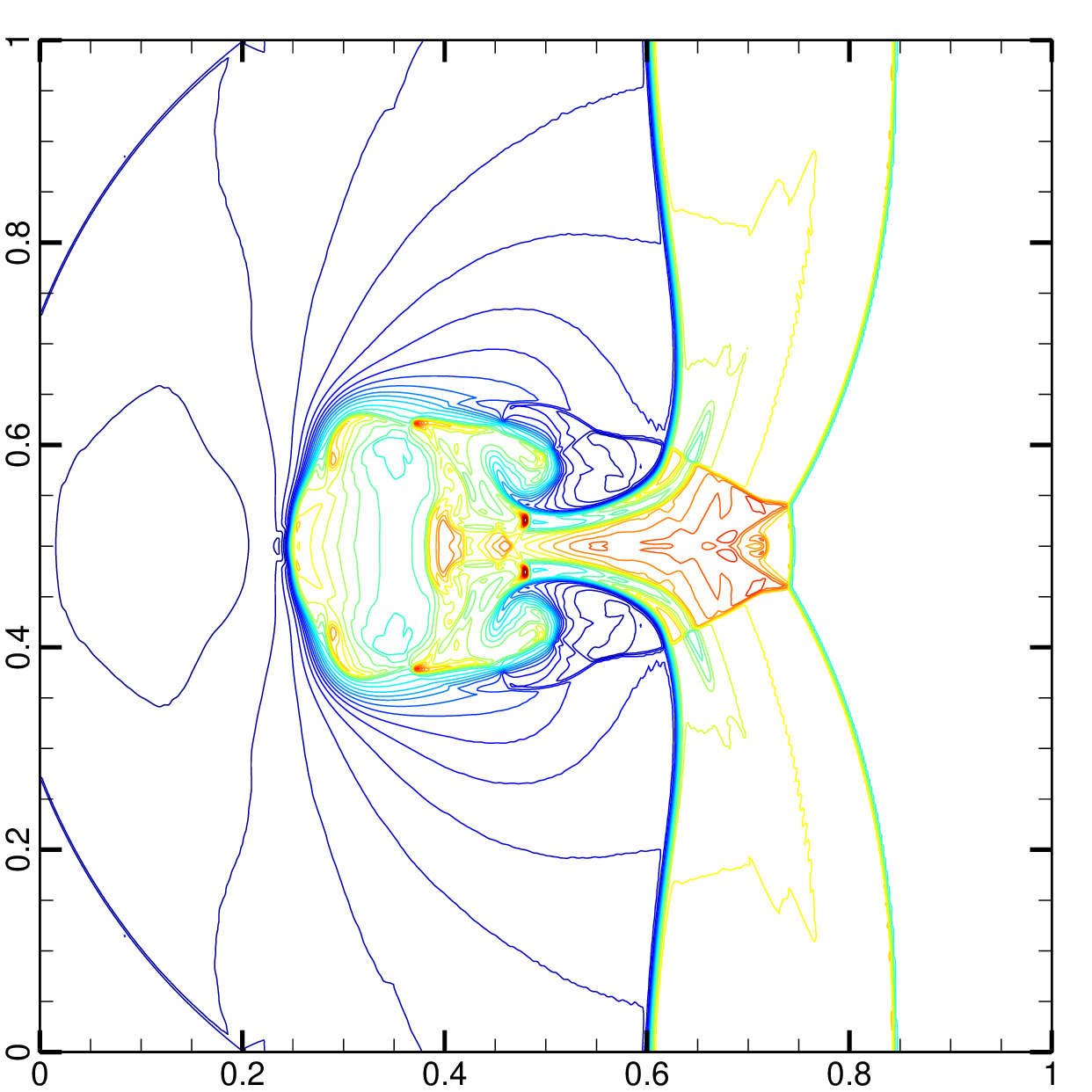}
		\caption{\(B_3\), \(N=512\)}
	\end{subfigure}

	\caption{Results of the cloud-shock interaction problem. For \(\ln\rho\), continuous coloring is used from \(0\) to \(4\). For \(B_3\), 30 equally spaced contours from \(-3.0\) to \(5.0\) are used.}
	\label{fig: cloud-shock}
\end{figure}

This problem, first proposed by Dai and Woodward \cite{10.1006/jcph.1998.5944}, describes the interaction between a dense cloud and a shock. Initially, a shock discontinuity is propagating along the \(x\)-direction with a Mach number 10, which is calculated as the ratio between the shock speed (by Rankine-Hugoniot jump condition) and the fast magnetosonic wave speed of the stationary ambient flow. The discontinuity is a combination of a fast shock wave and a rotational discontinuity in \(B_z\). The rotational discontinuity has no effect on all other flow variables. We use the same setup as \cite{10.1006/jcph.2000.6519}. The domain is \([0,1]^2\) with inflow boundary condition on the right boundary and outflow on the rest. Let \(r = \sqrt{(x-0.8)^2 + (y-0.5)^2}\). The initial condition is
\begin{equation*}
	(\rho, u_1, u_2, u_3, B_1, B_2, B_3, p) = \begin{cases}
		(3.86859, 0, 0, 0, 0, 2.1826182, -2.1826182, 167.345), & x < 0.6, \\
		(10, -11.2536, 0, 0, 0, 0.56418958, 0.56418958, 1), & \text{\(x > 0.6\) and \(r < 0.15\)}, \\
		(1, -11.2536, 0, 0, 0, 0.56418958, 0.56418958, 1), & \text{otherwise}.
	\end{cases}
\end{equation*}
We solve the problem up to \(T = 0.06\) with \(N \times N\) grids for \(N = 128, 256, 512\). The results shown in \Cref{fig: cloud-shock} match those in \cite{10.1137/18M1168042, 10.1007/s00211-019-01042-w, 10.1016/j.jcp.2025.113795} well.

\subsection{Astrophysical jet}

\begin{figure}[htbp]
	\centering

	\begin{subfigure}{0.3\linewidth}
		\centering
		\includegraphics[width=\linewidth]{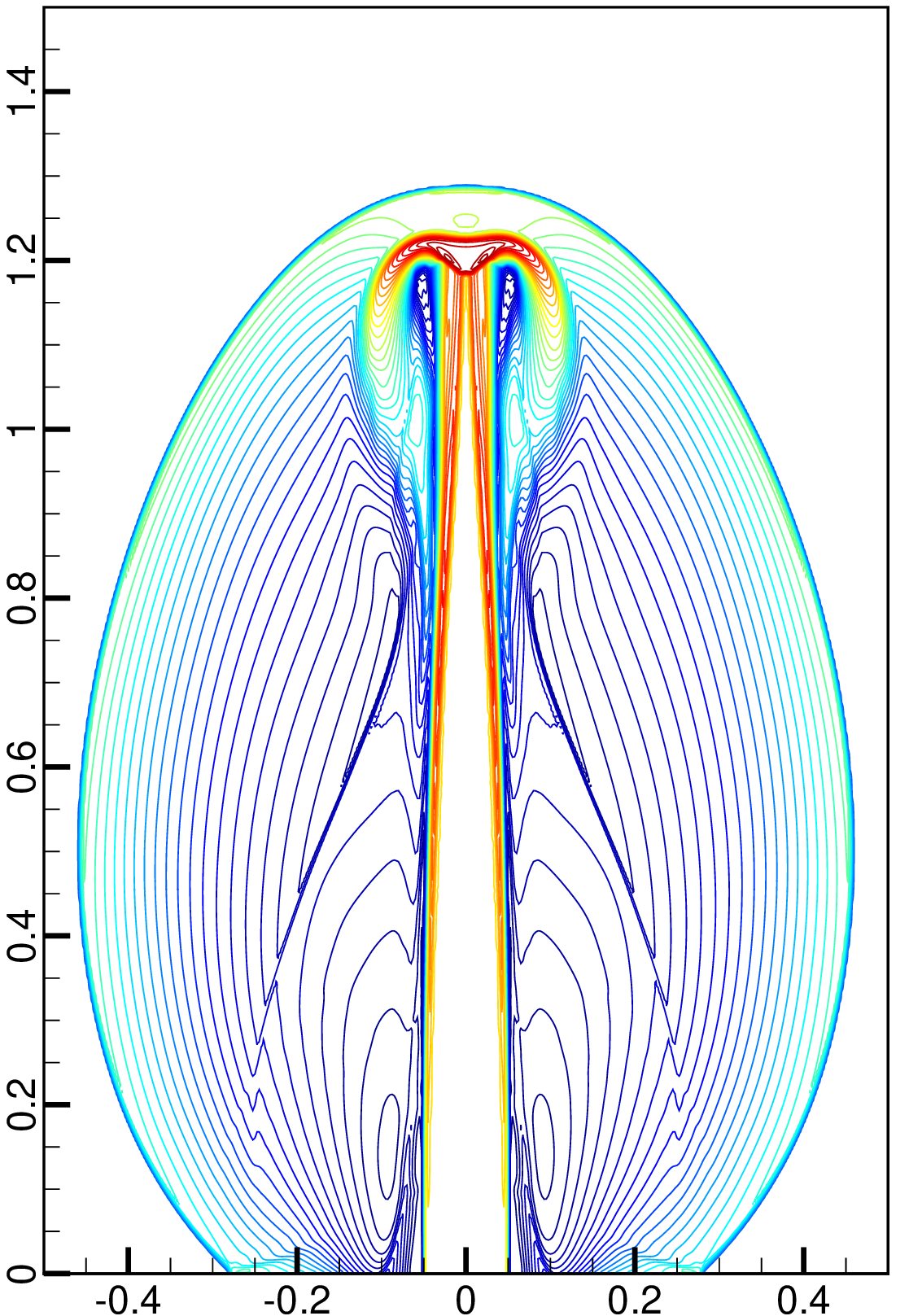}
		\caption{\(\log_{10}\rho\)}
	\end{subfigure}
	\begin{subfigure}{0.3\linewidth}
		\centering
		\includegraphics[width=\linewidth]{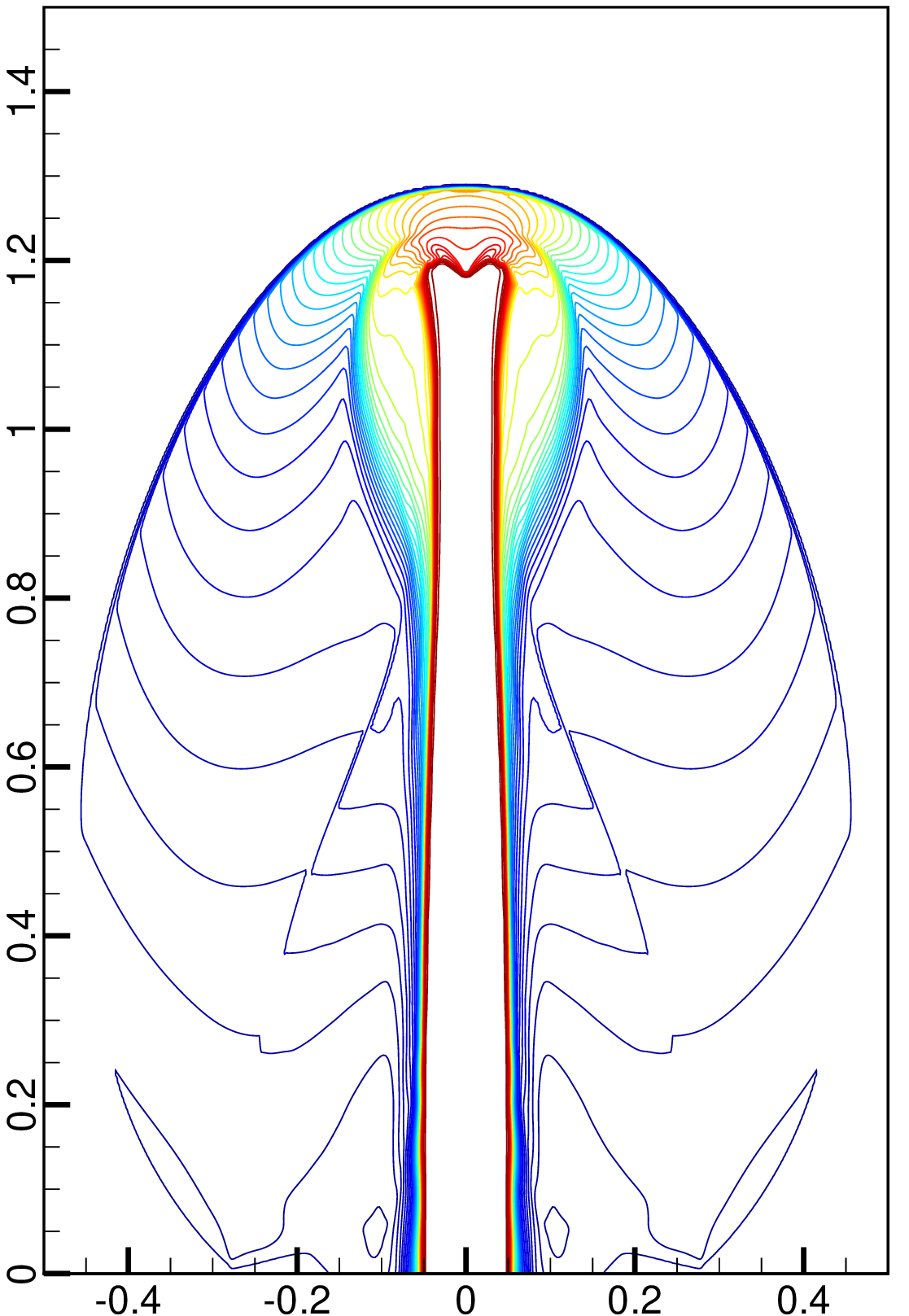}
		\caption{\(u_2\)}
	\end{subfigure}
	\begin{subfigure}{0.3\linewidth}
		\centering
		\includegraphics[width=\linewidth]{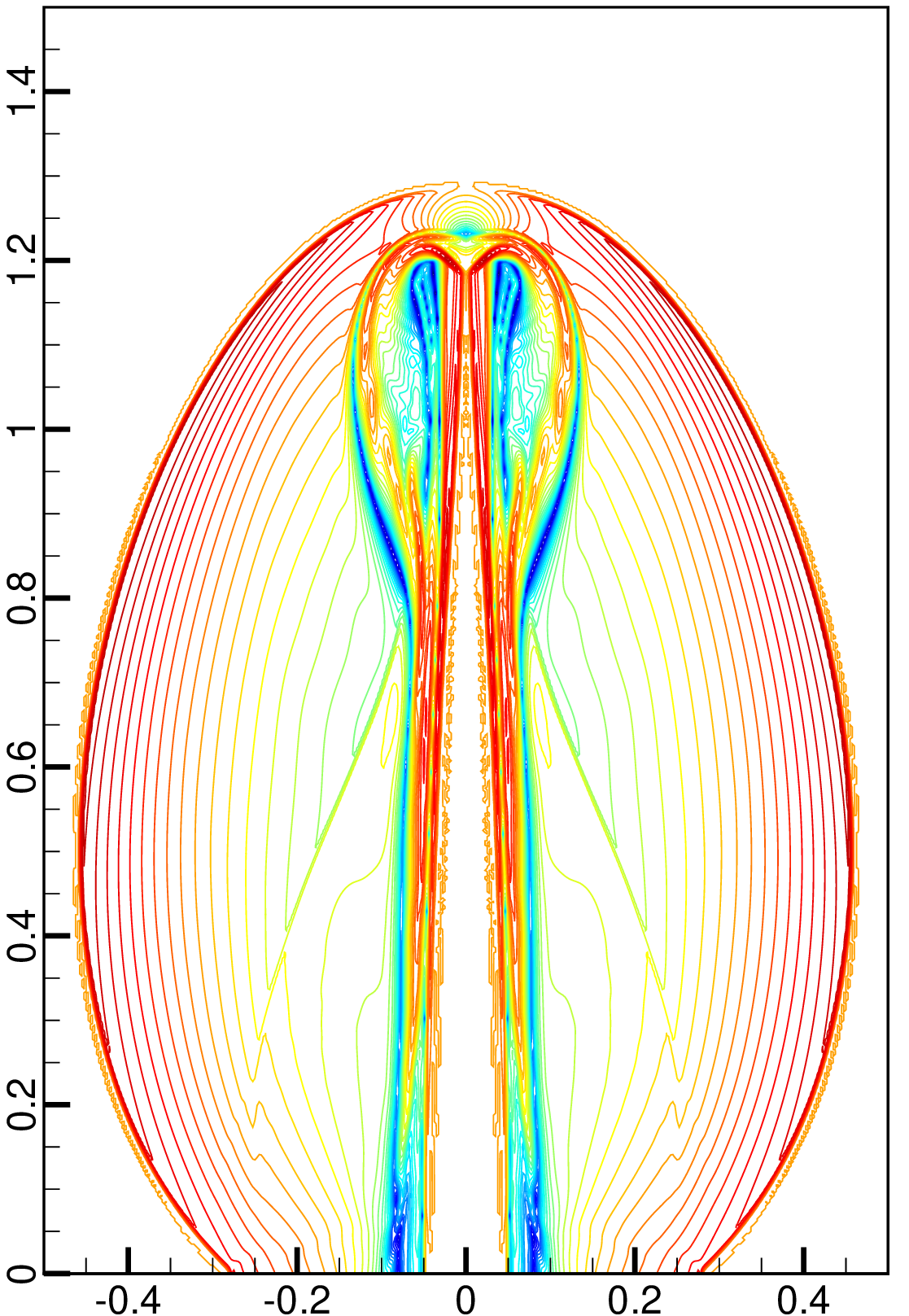}
		\caption{\(\log_{10}\frac{\abs{\bm B}^2}{2}\)}
	\end{subfigure}

	\caption{Contour plots of the astrophysical jet example. 40 equally spaced contours for \(\log_{10}\rho \in [-1.58, 1.12]\), \(v_2 \in [-40, 790]\), and \(\log_{10}\frac{\abs{\bm B}^2}{2} \in [-1.6, 3.6]\).}
	\label{fig: jet}
\end{figure}

This test is taken from \cite{10.1137/18M1168042}, where a magnetic field is added to a Mach 800 dense jet of \cite{10.1016/j.jcp.2012.01.032}, which is a tough case due to wide pressure range. Initially, the domain \([-0.5, 0.5] \times [0, 1.5]\) is filled with static gas with \(\rho = 0.1 \gamma\) and \(p = 1\). A uniform magnetic field \(\bm B = (0, \sqrt{200}, 0)^T\) is added. An inflow jet is imposed at the bottom boundary \(y = 0\) and \(\abs{x} < 0.05\), with \(\rho = \gamma\), \(p = 1\), and \(\bm v = (0, 800, 0)^T\). The specific heat ratio is set to \(\gamma = 1.4\). We use \(400 \times 600\) grids and compute to \(T = 2 \times 10^{-3}\). As shown in \Cref{fig: jet}, our scheme captures the head and beam/cocoon interface, and the internal structure is consistent to the one presented in \cite{10.1016/j.jcp.2012.01.032}.

\section{Conclusions}\label{se:con}

We have developed a nodal discontinuous Galerkin (DG) scheme for ideal MHD equations with the Godunov--Powell source term that is both positivity preserving (PP) and semi-discretely entropy stable (ES). Preservation of positivity keeps the simulation from breaking down, and entropy stability will help the numerical solution to converge to the unique entropy solution. We start from the original flux differencing based nodal DG scheme and incorporate the positivity preserving HLL numerical flux. For each element, entropy conservation in volume is reached via the two-point entropy conservative flux function and entropy dissipation at interface is achieved by our carefully selected signal speed estimates. Local divergence-free projection is applied so that the updated cell average is admissible and thus the PP limiter is applicable. Future work may include the extension to unstructured meshes and an upgrade to fully discretely PP and ES scheme.

\begin{paragraph}{Declaration of conflict of interest}
	The authors declare that there is no conflict of interest associated with this work.
\end{paragraph}

\begin{paragraph}{Funding}
	Research is supported in part by NSF grant DMS-2309249.
\end{paragraph}



\bibliographystyle{amsplain}
\bibliography{refs}

\providecommand{\bysame}{\leavevmode\hbox to3em{\hrulefill}\thinspace}
\providecommand{\MR}{\relax\ifhmode\unskip\space\fi MR }
\providecommand{\MRhref}[2]{%
  \href{http://www.ams.org/mathscinet-getitem?mr=#1}{#2}
}
\providecommand{\href}[2]{#2}
\begin{thebibliography}{10}

\bibitem{10.1016/j.jcp.2018.06.031}
R\'{e}mi Abgrall, \emph{A general framework to construct schemes satisfying additional conservation relations. application to entropy conservative and entropy dissipative schemes}, Journal of Computational Physics \textbf{372} (2018), 640--666.

\bibitem{10.1016/j.jcp.2012.01.032}
Dinshaw~S. Balsara, \emph{Self-adjusting, positivity preserving high order schemes for hydrodynamics and magnetohydrodynamics}, Journal of Computational Physics \textbf{231} (2012), no.~22, 7504--7517.

\bibitem{10.2140/camcos.2021.16.59}
Dinshaw~S. Balsara, Rakesh Kumar, and Praveen Chandrashekar, \emph{Globally divergence-free {DG} scheme for ideal compressible {MHD}}, Communications in Applied Mathematics and Computational Science \textbf{16} (2021), no.~1, 59--98.

\bibitem{10.1006/jcph.1998.6153}
Dinshaw~S. Balsara and Daniel~S. Spicer, \emph{A staggered mesh algorithm using high order {G}odunov fluxes to ensure solenoidal magnetic fields in magnetohydrodynamic simulations}, Journal of Computational Physics \textbf{149} (1999), no.~2, 270--292.

\bibitem{10.1007/978-3-642-58535-7_5}
Timothy~John Barth, \emph{Numerical methods for gasdynamic systems on unstructured meshes}, An Introduction to Recent Developments in Theory and Numerics for Conservation Laws (Dietmar Kr\"{o}ner, Mario Ohlberger, and Christian Rohde, eds.), Lecture Notes in Computational Science and Engineering, vol.~5, Springer, Berlin, Heidelberg, 12 1998, pp.~195--285.

\bibitem{10.1007/0-387-38034-5_4}
\bysame, \emph{On the role of involutions in the discontinuous {G}alerkin discretization of {M}axwell and magnetohydrodynamic systems}, Compatible Spatial Discretizations (Douglas~Norman Arnold, Pavel~Blagovestov Bochev, Richard~Bruno Lehoucq, Roy~A. Nicolaides, and Mikhail~Jurievich Shashkov, eds.), The IMA Volumes in Mathematics and its Applications, vol. 142, Springer, New York, 07 2006, pp.~69--88.

\bibitem{10.1016/j.jcp.2018.06.027}
Marvin Bohm, Andrew~R. Winters, Gregor~J. Gassner, Dominik Derigs, Florian~J. Hindenlang, and Joachim Saur, \emph{An entropy stable nodal discontinuous {G}alerkin method for the resistive mhd equations. part {I}: Theory and numerical verification}, Journal of Computational Physics \textbf{422} (2020), 108076.

\bibitem{10.1007/s00211-007-0108-8}
Fran\c{c}ois Bouchut, Christian Klingenberg, and Knut Waagan, \emph{A multiwave approximate {R}iemann solver for ideal {MHD} based on relaxation. {I}: theoretical framework}, Numerische Mathematik \textbf{108} (2007), no.~1, 7--42.

\bibitem{10.1007/s00211-010-0289-4}
\bysame, \emph{A multiwave approximate {R}iemann solver for ideal {MHD} based on relaxation {II}: numerical implementation with 3 and 5 waves}, Numerische Mathematik \textbf{115} (2010), no.~4, 647--679.

\bibitem{10.1016/0021-99918090079-0}
Jeremiah~Uhler Brackbill and Daniel~Charles Barnes, \emph{The effect of nonzero $\nabla\cdot\bm{B}$ on the numerical solution of the magnetohydrodynamic equations}, Journal of Computational Physics \textbf{35} (1980), no.~3, 426--430.

\bibitem{10.1016/j.jcp.2018.02.033}
Jesse Chan, \emph{On discretely entropy conservative and entropy stable discontinuous {G}alerkin methods}, Journal of Computational Physics \textbf{362} (2018), 346--374.

\bibitem{10.1137/18M1209234}
Jesse Chan, David~Cesar Del Rey~Fern\'{a}ndez, and Mark~Huitt Carpenter, \emph{Efficient entropy stable {G}auss collocation methods}, SIAM Journal on Scientific Computing \textbf{41} (2019), no.~5, A2938--A2966.

\bibitem{10.1137/15M1013626}
Praveen Chandrashekar and Christian Klingenberg, \emph{Entropy stable finite volume scheme for ideal compressible {MHD} on 2-{D} {C}artesian meshes}, SIAM Journal on Numerical Analysis \textbf{54} (2016), no.~2, 1313--1340.

\bibitem{10.1007/s10915-020-01289-8}
Praveen Chandrashekar and Rakesh Kumar, \emph{Constraint preserving discontinuous {G}alerkin method for ideal compressible {MHD} on 2-{D} {C}artesian grids}, Journal of Scientific Computing \textbf{84} (2020), no.~2, 39.

\bibitem{10.1016/j.jcp.2017.05.025}
Tianheng Chen and Chi-Wang Shu, \emph{Entropy stable high order discontinuous {G}alerkin methods with suitable quadrature rules for hyperbolic conservation laws}, Journal of Computational Physics \textbf{345} (2017), 427--461.

\bibitem{10.4208/csiam-am.2020-0003}
\bysame, \emph{Review of entropy stable discontinuous {G}alerkin methods for systems of conservation laws on unstructured simplex meshes}, CSIAM Transactions on Applied Mathematics \textbf{1} (2020), no.~1, 1--52.

\bibitem{10.1016/j.jcp.2012.12.019}
Yue Cheng, Fengyan Li, Jianxian Qiu, and Liwei Xu, \emph{Positivity-preserving {DG} and central {DG} methods for ideal {MHD} equations}, Journal of Computational Physics \textbf{238} (2013), 255--280.

\bibitem{10.1007/978-3-642-59721-3}
Bernardo Cockburn, George~Em Karniadakis, and Chi-Wang Shu, \emph{Discontinuous {Galerkin} methods}, 1 ed., Lecture Notes in Computational Science and Engineering, Springer, Berlin, Heidelberg, 2000.

\bibitem{0021999189901836}
Bernardo Cockburn, San-Yih Lin, and Chi-Wang Shu, \emph{{TVB} {R}unge--{K}utta local projection discontinuous {G}alerkin finite element method for conservation laws {III}: One-dimensional systems}, Journal of Computational Physics \textbf{84} (1989), no.~1, 90--113.

\bibitem{10.1006/jcph.1998.5892}
Bernardo Cockburn and Chi-Wang Shu, \emph{The {R}unge--{K}utta discontinuous {G}alerkin method for conservation laws {V}: Multidimensional systems}, Journal of Computational Physics \textbf{141} (1998), no.~2, 199--224.

\bibitem{10.1016/j.jcp.2017.12.015}
Jared Crean, Jason~Edward Hicken, David~Cesar Del Rey~Fern\'{a}ndez, David~Walter Zingg, and Mark~Huitt Carpenter, \emph{Entropy-stable summation-by-parts discretization of the {E}uler equations on general curved elements}, Journal of Computational Physics \textbf{356} (2018), 410--438.

\bibitem{10.1007/BF00280911}
Constantine~Michael Dafermos, \emph{Quasilinear hyperbolic systems with involutions}, Archive for Rational Mechanics and Analysis \textbf{94} (1986), no.~4, 373--389.

\bibitem{Hyperbolic_Dafermos_4th}
\bysame, \emph{Hyperbolic conservation laws in continuum physics}, 5 ed., Grundlehren der mathematischen Wissenschaften, vol. 325, Springer, Berlin, Heidelberg, 03 2026.

\bibitem{10.1063/1.859081}
Russel~B. Dahlburg and J.~Michael Picone, \emph{Evolution of the {O}rszag--{T}ang vortex system in a compressible medium. {I}. initial average subsonic flow}, Physics of Fluids B: Plasma Physics \textbf{1} (1989), no.~11, 2153--2171.

\bibitem{10.1006/jcph.1998.5944}
Wenlong Dai and Paul~R. Woodward, \emph{A simple finite difference scheme for multidimensional magnetohydrodynamical equations}, Journal of Computational Physics \textbf{142} (1998), no.~2, 331--369.

\bibitem{10.1006/jcph.2001.6961}
Andreas~S. Dedner, Friedemann Kemm, Dietmar Kr\"{o}ner, Claus-Dieter Munz, T.~Schnitzer, and Matthias Wesenberg, \emph{Hyperbolic divergence cleaning for the {MHD} equations}, Journal of Computational Physics \textbf{175} (2002), no.~2, 645--673.

\bibitem{10.1016/j.compfluid.2014.02.016}
David~Cesar Del Rey~Fern\'{a}ndez, Jason~Edward Hicken, and David~Walter Zingg, \emph{Review of summation-by-parts operators with simultaneous approximation terms for the numerical solution of partial differential equations}, Computers and Fluids \textbf{95} (2014), 171--196.

\bibitem{10.1016/j.jcp.2018.03.002}
Dominik Derigs, Andrew~R. Winters, Gregor~J. Gassner, Stefanie Walch, and Marvin Bohm, \emph{Ideal {GLM}-{MHD}: About the entropy consistent nine-wave magnetic field divergence diminishing ideal magnetohydrodynamics equations}, Journal of Computational Physics \textbf{364} (2018), 420--467.

\bibitem{10.1007/978-3-030-56341-7}
Alexandre Ern and Jean-Luc Guermond, \emph{Finite elements {I}}, Texts in Applied Mathematics, vol.~72, Springer, Cham, 02 2021.

\bibitem{10.1086/166684}
Charles~R. Evans and John~F. Hawley, \emph{Simulation of magnetohydrodynamic flows: A constrained transport method}, The Astrophysical journal \textbf{332} (1988), 659--677.

\bibitem{10.1137/110836961}
Ulrik~S Fjordholm, Siddhartha Mishra, and Eitan Tadmor, \emph{Arbitrarily high-order accurate entropy stable essentially nonoscillatory schemes for systems of conservation laws}, SIAM Journal on Numerical Analysis \textbf{50} (2012), no.~2, 544--573.

\bibitem{10.1007/s10915-018-0750-6}
Pei Fu, Fengyan Li, and Yan Xu, \emph{Globally divergence-free discontinuous {G}alerkin methods for ideal magnetohydrodynamic equations}, Journal of Scientific Computing \textbf{77} (2018), no.~3, 1621--1659.

\bibitem{10.1016/j.jcp.2004.11.016}
Thomas~A. Gardiner and James~M. Stone, \emph{An unsplit {Godunov} method for ideal {MHD} via constrained transport}, Journal of Computational Physics \textbf{205} (2005), no.~2, 509--539.

\bibitem{Godunov1972}
Sergei~Konstantinovich Godunov, \emph{Symmetric form of the equations of magnetohydrodynamics}, Numerical Methods of Continuum Mechanics \textbf{3} (1972), no.~1, 26--34.

\bibitem{10.1090/S0025-5718-98-00913-2}
Sigal Gottlieb and Chi-Wang Shu, \emph{Total variation diminishing {R}unge-{K}utta schemes}, Mathematics of Computation \textbf{67} (1998), no.~221, 73--85.

\bibitem{10.1137/S003614450036757X}
Sigal Gottlieb, Chi-Wang Shu, and Eitan Tadmor, \emph{Strong stability-preserving high-order time discretization methods}, SIAM Review \textbf{43} (2001), no.~1, 89--112.

\bibitem{10.1016/j.jcp.2016.05.054}
Jean-Luc Guermond and Bojan Popov, \emph{Fast estimation from above of the maximum wave speed in the {R}iemann problem for the {E}uler equations}, Journal of Computational Physics \textbf{321} (2016), 908--926.

\bibitem{10.1016/j.jcp.2011.02.009}
Christiane Helzel, Rossmanith James~Alexander, and Bertram Taetz, \emph{An unstaggered constrained transport method for the {3D} ideal magnetohydrodynamic equations}, Journal of Computational Physics \textbf{230} (2011), no.~10, 3803--3829.

\bibitem{hindenlang2019new}
Florian~J. Hindenlang and Gregor~J. Gassner, \emph{A new entropy conservative two-point flux for ideal {MHD} equations derived from first principles}, Talk presented at HONOM, 04 2019.

\bibitem{10.1137/050627022}
Rossmanith James~Alexander, \emph{An unstaggered, high-resolution constrained transport method for magnetohydrodynamic flows}, SIAM Journal on Scientific Computing \textbf{28} (2006), no.~5, 1766--1797.

\bibitem{jameson2006eigenvalues}
Antony Jameson, \emph{Eigenvalues, eigenvectors and symmetrization of the magneto-hydrodynamic ({MHD}) equations}, AFOSR Grantees and Contractors Meeting, 06 2006.

\bibitem{10.1006/jcph.2000.6479}
Pekka Janhunen, \emph{A positive conservative method for magnetohydrodynamics based on {HLL} and {R}oe methods}, Journal of Computational Physics \textbf{160} (2000), no.~2, 649--661.

\bibitem{10.1006/jcph.1996.0130}
Guang-Shan Jiang and Chi-Wang Shu, \emph{Efficient implementation of weighted {ENO} schemes}, Journal of Computational Physics \textbf{126} (1996), no.~1, 202--228.

\bibitem{10.2140/camcos.2013.8.1}
Friedemann Kemm, \emph{On the origin of divergence errors in {MHD} simulations and consequences for numerical schemes}, Communications in Applied Mathematics and Computational Science \textbf{8} (2013), no.~1, 1--38.

\bibitem{10.4208/cicp.180515.230616a}
Christian Klingenberg, Frank P\"{o}rner, and Yinhua Xia, \emph{An efficient implementation of the divergence free constraint in a discontinuous {G}alerkin method for magnetohydrodynamics on unstructured meshes}, Communications in Computational Physics \textbf{21} (2017), no.~2, 423--442.

\bibitem{10.1007/s10915-004-4146-4}
Fengyan Li and Chi-Wang Shu, \emph{Locally divergence-free discontinuous galerkin methods for {MHD} equations}, Journal of Scientific Computing \textbf{22} (2005), no.~1, 413--442.

\bibitem{10.1016/j.jcp.2025.113795}
Mengqing Liu and Kailiang Wu, \emph{Structure-preserving oscillation-eliminating discontinuous galerkin schemes for ideal {MHD} equations: Locally divergence-free and positivity-preserving}, Journal of Computational Physics \textbf{527} (2025), 113795.

\bibitem{10.1137/21M140835X}
Yong Liu, Jianfang Lu, and Chi-Wang Shu, \emph{An essentially oscillation-free discontinuous {Galerkin} method for hyperbolic systems}, SIAM Journal on Scientific Computing \textbf{44} (2022), no.~1, A230--A259.

\bibitem{10.1016/j.jcp.2025.113911}
\bysame, \emph{An entropy stable essentially oscillation-free discontinuous {G}alerkin method for solving ideal magnetohydrodynamic equations}, Journal of Computational Physics \textbf{530} (2025), 113911.

\bibitem{10.1016/j.jcp.2017.10.043}
Yong Liu, Chi-Wang Shu, and Mengping Zhang, \emph{Entropy stable high order discontinuous {G}alerkin methods for ideal compressible {MHD} on structured meshes}, Journal of Computational Physics \textbf{354} (2018), 163--178.

\bibitem{10.1090/mcom/3998}
Manting Peng, Zheng Sun, and Kailiang Wu, \emph{{OEDG}: oscillation eliminating discontinuous {G}alerkin method for hyperbolic conservation laws}, Mathematics of Computation \textbf{94} (2025), no.~353, 1147--1198.

\bibitem{Powell1994NASA}
Kenneth~Grant Powell, \emph{An approximate {R}iemann solver for magnetohydrodynamics}, Tech. Report ICASE-94-24, NASA Langley Research Center, Hampton, VA, 04 1994.

\bibitem{10.1006/jcph.1999.6299}
Kenneth~Grant Powell, Philip~L. Roe, Timur~J. Linde, Tamas~I. Gombosi, and Darren~L. De~Zeeuw, \emph{A solution-adaptive upwind scheme for ideal magnetohydrodynamics}, Journal of Computational Physics \textbf{154} (1999), no.~2, 284--309.

\bibitem{10.1137/S1064827503425298}
Jianxian Qiu and Chi-Wang Shu, \emph{{Runge}--{Kutta} discontinuous {Galerkin} method using {WENO} limiters}, SIAM Journal on Scientific Computing \textbf{26} (2005), no.~3, 907--929.

\bibitem{10.1016/j.jcp.2023.112607}
Andr\'{e}s~Mauricio Rueda-Ram\'{\i}rez and Gregor~J. Gassner, \emph{A flux-differencing formula for split-form summation by parts discretizations of non-conservative systems: Applications to subcell limiting for magneto-hydrodynamics}, Journal of Computational Physics \textbf{496} (2024), 112607.

\bibitem{10.1016/j.jcp.2022.111851}
Andr\'{e}s~Mauricio Rueda-Ram\'{\i}rez, Florian~J. Hindenlang, Jesse Chan, and Gregor~J. Gassner, \emph{Entropy-stable {G}auss collocation methods for ideal magneto-hydrodynamics}, Journal of Computational Physics \textbf{475} (2023), 111851.

\bibitem{10.1017/S0962492920000057}
Chi-Wang Shu, \emph{Essentially non-oscillatory and weighted essentially non-oscillatory schemes}, Acta Numerica \textbf{29} (2020), 701–762.

\bibitem{10.1017/S0962492902000156}
Eitan Tadmor, \emph{Entropy stability theory for difference approximations of nonlinear conservation laws and related time-dependent problems}, Acta Numerica \textbf{12} (2003), 451--512.

\bibitem{10.1007/BF01414629}
Eleuterio~Francisco Toro, Michael Spruce, and William Speares, \emph{Restoration of the contact surface in the {HLL}-{R}iemann solver}, Shock Waves \textbf{4} (1994), no.~1, 25--34.

\bibitem{10.1006/jcph.2000.6519}
G\'{a}bor T\'{o}th, \emph{The $\nabla\cdot\bm{B}=0$ constraint in shock-capturing magnetohydrodynamics codes}, Journal of Computational Physics \textbf{161} (2000), no.~2, 605--652.

\bibitem{10.1016/j.jcp.2024.113498}
Lei Wei, Lingling Zhou, and Yinhua Xia, \emph{The jump filter in the discontinuous {Galerkin} method for hyperbolic conservation laws}, Journal of Computational Physics \textbf{520} (2025), 113498.

\bibitem{10.1016/j.jcp.2015.09.055}
Andrew~R. Winters and Gregor~J. Gassner, \emph{Affordable, entropy conserving and entropy stable flux functions for the ideal {MHD} equations}, Journal of Computational Physics \textbf{304} (2016), 72--108.

\bibitem{10.1137/18M1168017}
Kailiang Wu, \emph{Positivity-preserving analysis of numerical schemes for ideal magnetohydrodynamics}, SIAM Journal on Numerical Analysis \textbf{56} (2018), no.~4, 2124--2147.

\bibitem{10.1137/18M1168042}
Kailiang Wu and Chi-Wang Shu, \emph{A provably positive discontinuous {G}alerkin method for multidimensional ideal magnetohydrodynamics}, SIAM Journal on Scientific Computing \textbf{40} (2018), no.~5, B1302--B1329.

\bibitem{10.1007/s00211-019-01042-w}
\bysame, \emph{Provably positive high-order schemes for ideal magnetohydrodynamics: analysis on general meshes}, Numerische Mathematik \textbf{142} (2019), no.~4, 995--1047.

\bibitem{10.1137/21M1458247}
\bysame, \emph{Geometric quasilinearization framework for analysis and design of bound-preserving schemes}, SIAM Review \textbf{65} (2023), no.~4, 1031--1073.

\bibitem{10.1007/3-540-27170-8_42}
Helen~C. Yee and Bj\"orn Sj\"ogreen, \emph{Divergence free high order filter methods for the compressible {MHD} equations}, Modeling, Simulation and Optimization of Complex Processes (Berlin, Heidelberg) (Hans~Georg Bock, Hoang~Xuan Phu, Ekaterina Kostina, and Rolf Rannacher, eds.), Springer, Berlin, Heidelberg, 2005, pp.~559--575.

\bibitem{10.1007/s10915-005-9004-5}
\bysame, \emph{Efficient low dissipative high order schemes for multiscale {MHD} flows, {II}: Minimization of $\nabla\cdot\bm{B}$ numerical error}, Journal of Scientific Computing \textbf{29} (2006), no.~1, 115--164.

\bibitem{10.1016/j.jcp.2010.08.016}
Xiangxiong Zhang and Chi-Wang Shu, \emph{On positivity-preserving high order discontinuous {G}alerkin schemes for compressible {E}uler equations on rectangular meshes}, Journal of Computational Physics \textbf{229} (2010), no.~23, 8918--8934.

\bibitem{10.1016/j.jcp.2012.08.028}
Xinghui Zhong and Chi-Wang Shu, \emph{A simple weighted essentially nonoscillatory limiter for {Runge}--{Kutta} discontinuous {Galerkin} methods}, Journal of Computational Physics \textbf{232} (2013), no.~1, 397--415.

\end{thebibliography}


\end{document}